\input amstex
\documentstyle{amsppt}
\magnification=1200 \hsize6.5truein \vsize9truein
\catcode`\@=11
\let\logo@\relax
\catcode`\@=\active
\NoBlackBoxes
\TagsOnLeft
\define\SV{\operatorname{Sing}V}
\define\7{\operatorname{Supp}\,}
\define\lE{\operatorname{log}E}
\define\ind{\indent}
\define\C{\Cal}
\define\dom{\operatorname{dom}\,}
\define\im{\operatorname{im}\,}

\define\da{\partial}
\define\db{\bar\partial}
\define\ti{\tilde}

\define\ite#1{\item"{\it #1.}"}
\define\f#1#2{\frac{#1}{#2}}
\define \p{\phi}
\define \pr{\prime}
\define\T{\operatorname{Torsion}}
\define\SW{\operatorname{Sing}W}
\define\inc{\hookrightarrow}
\define\di{\operatorname{div}}
\nologo
\topmatter

\title
Pure Hodge structure on the $L_2$-cohomology of varieties with isolated singularities
\endtitle

\rightheadtext {Pure Hodge structure on $L_2$-cohomology.}
\author
William Pardon and Mark Stern
\endauthor

\date
31 October 1997
\enddate

\address {
\newline
Departament of Mathematics \newline
Duke University\newline
Durham, NC 27708-0320 \newline
USA \newline
e-mail:wlp\@math.duke.edu, stern\@math.duke.edu }
\endaddress

\endtopmatter

\document

\footnote""{Partially supported by NSF grants DMS 95-04900 (Pardon) and DMS 9505040 (Stern)\hfill }

\heading Introduction.\endheading
Cheeger, Goresky, and MacPherson 
conjectured [CGM] an $L_2$-de Rham theorem: that the intersection cohomology of a projective 
variety $V$ is naturally 
isomorphic to the $L_2$-cohomology of the incomplete manifold
$V-\SV$, with metric induced by a projective embedding. 
The early interest in this conjecture 
was motivated in large part by the hope that one could then put a Hodge structure on the intersection cohomology of $V$ and even extend the rest of the ``K\"ahler package" (\cite{CGM}) to 
this context. Saito [S1,S2] eventually established the Kahler package for intersection cohomology without recourse to $L_2$-cohomology techniques. Interest 
in $L_2$-cohomology did not disappear with this result, since, among other things, $L_2$-cohomology provides intrinsic geometric invariants of an arbitrary complex projective variety which are not apparent from the point of view of $D$-modules. For instance, $L_2-\db$-coholomology groups depend on boundary conditions ([PS]), which, as we show here, must be treated carefully in order to give the correct Hodge components for the $L_2$-cohomology of a singular variety.  

 It was quickly realized, however, that for incomplete 
manifolds the Hodge and Lefschetz decompositions are not direct consequences of the K\"ahler
condition as they are in the complete case.  The primary obstruction 
to obtaining a Hodge structure on the $L_2$-cohomology is 
an apparent technicality: on an incomplete K\"ahler manifold there are 
several potentially  distinct definitions of a square integrable harmonic form. For example, a form $h$ might be considered harmonic if 
 $dh = 0 =  \delta h$, or if $\db h = 0 = \vartheta h$, or 
simply if $\Delta h = 0$.  Moreover there are further 
domain considerations: one imposes boundary conditions, which turn out to have no effect on cohomology in the case of $d$, but are crucial for $\db$-cohomology.  

On a complete manifold all these definitions of harmonics coincide, and one obtains the 
Hodge decomposition by  decomposing harmonic forms into 
their $(p,q)$ components. The $(p,q)$ components are harmonic in the weakest sense - they are in the kernel of $\Delta$. The equality of the different 
notions of harmonic then allows one to realize these $(p,q)$ 
components as representing both $\db$ and $d$ cohomology classes. 
The equivalence of the different definitions of harmonic is also 
required in order to obtain the Lefschetz decomposition.  Interior 
product with the Kahler form preserves the kernel of $\Delta$ by 
local computation. One requires the equivalence to see that it 
also preserves the kernel  of $d$.  
  
   Ohsawa [O2] proved the conjectured $L_2$-de Rham theorem under the extra assumption that $V$ has only isolated 
singularities. Strangely, the $L_2$-cohomology in the incomplete metric 
played almost no role in his proof. It entered only as a limit of 
cohomology groups with respect to a family of auxiliary complete metrics, 
 which degenerate to the incomplete metric. The proof relies on 
earlier work of Saper [Sa] where $V-\SV$
(under the isolated singularities assumption) is endowed with 
  a complete Kahler metric whose associated 
 $L_2$-cohomology is isomorphic to the intersection cohomology of $V$. 
Saper's result also provides the intersection cohomology of varieties 
with isolated singularities with a (non-canonical) Hodge 
decomposition. 

   This paper began as an attempt to compute those $L_2$-$\db$ cohomology groups for surfaces which had not been computed in [PS] and to show that
they give the same Hodge structure as Hain and Zucker obtain in 
[HZ] using resolution of singulaties. In order to show that the $\db$ cohomology groups actually 
gave the components of the Hodge structure, we were forced to 
overcome the above technical difficulties and to understand the relations 
between the domains of $\Delta,d,\db,$ etc. Ultimately we were led 
to establish a "good" harmonic theory for varieties with isolated 
singularities. To do so, we show that the harmonic forms 
satisfy certain growth estimates near the singular points. With these estimates we can manipulate the harmonic 
forms as though they were on a complete manifold. For example, 
we show that in degrees other than $n-1$, $n$, and $n+1$, $n = dim_\Bbb C V$, 
the $L_2$ kernel of $\Delta$ is contained in the kernels of $d$, $\delta$, 
$\db$, and $\vartheta$. With this and related results we obtain the pure Hodge decomposition and 
the Lefschetz decomposition for the $L_2$-cohomology, in the same manner 
as in the complete case.

   Some of the estimates we derive could also be obtained 
by appropriately elaborating arguments of [O1] and [O2]. They are proved here, however, using the incomplete metric itself
rather than families of auxiliary complete metrics, because we hope to develop the tools for working directly with the $L_2$-complex in the incomplete metric. We have not, 
however, reproved all the results that we need from [O1] and [O2]. In particular, we do not reprove the isomorphism between
 intersection cohomology and $L_2$-cohomology. It is clear that our estimates 
do not yet imply the requisite vanishing in middle degree, which 
in Ohsawa's proof ultimately relies on a computation of Saper [Sa].

As we mentioned above, this line of investigation began with an attempt to calculate those $L_2-\db$-cohomology groups of an algebraic surface which had not already been calculated in [P] and [PS]. Since the $L_2-d$-cohomology had been computed locally, it seemed natural to put our calculations and their relation to a proposed Hodge filtration into a local context using the derived category,  which would globalize to the required Hodge structure. Thus, the second main result of this paper is to show that $(\C L^\cdot_{N/D},\C F^\cdot)$ admits the structure cohomological Hodge complex (\cite{D}), where $\C L^\cdot_{N/D}$ is a complex of sheaves of $L_2$ forms (with mixed Neumann and Dirichlet bounday conditions) and $\C F^\cdot$ is the filtration by holomorphic degree. We have proved this result only for complex surfaces.  

Related work has been done by Fox and Haskell [FH] and Nagase [N] 
in which a Hodge decomposition was stated for the $L_2$-cohomology 
of normal singular surfaces. Their work relies implicitly 
on the inclusion in  degree $0$  
of the $L_2$ kernel of $\Delta$ in the kernel of $d$. The K\"ahler package for curves was proved by Br\"uning and Lesch in [BL].

 In the next section we give a precise statement of our main results.

\newpage

\heading \S 1: Statement of the main results.\endheading

We begin with a review of some of the basic ideas of $L^2$-cohomology for complex varieties (\cite{P, PS}). Let $V$ be a projective variety with singular set $\SV$. From any projective imbedding $V\hookrightarrow M$, where $M$ is a compact K\"ahler manifold, $V-\SV$ inherits a K\"ahler metric $g$, which we call an {\it ambient metric}; or {\it Fubini-Study metric}, if $M=\Bbb P^N$. It is incomplete if $\SV\neq\emptyset$. The pointwise inner product of $k$-forms $\omega_1$ and $\omega_2$ with measurable coefficients will be denoted $<\omega_1,\omega_2>_g$, and the global inner product is
$$ 
(\omega_1,\omega_2)_g:=\int_{V-\SV} <\omega_1,\omega_2>_g\,dV_g
$$
where $dV_g$ denotes the volume form of the ambient metric; subscript $g$'s will be dropped in general and the norm of a form $\omega$ will be denoted $\|\omega\|$. Since the  quasi-isometry class of an ambient metric is independent of the choice of $M$ or its metric, there is,  for each nonnegative integer $k$, a well-defined sheaf $\C L^k$ on $V$ of {\it locally} $L^2$ {\it k-forms}: if $\Cal M^k$ denotes the sheaf on $V-\SV$ of $k$-forms with measurable coefficients, then for each open set $U\subseteq V$ 
$$
 \C L^k(U):=\{\omega\in \C M^k(U-\SV)\,|\,\|\,\omega|K\,\|<\infty,\,\text{for all compact}\,K\subseteq U\}\tag1.1
$$
Now since $V$ is compact the space of global sections $\C L^k(V)$ is a Hilbert space with respect to the inner product $(\cdot,\cdot)$ and the subpaces $\C A^k_c(V-\SV)$ of smooth compactly supported forms and $\C A^k(V-\SV)$ of smooth forms are dense. The exterior derivative $d_{cpt}:\C A^k_c(V-\SV)\to \C A^{k+1}_c(V-\SV)$ admits (at least) two closed extensions to operators $d_N$ and $d_D$ from $\C L^k(V)$ to $\C L^{k+1}(V)$, the {\it Neumann} and {\it Dirichlet}
extensions, so named because of the analogy with classical boundary conditions: these are, respectively, the graph closures of $d$ restricted to $\C A^k(V-\SV)$ and of $d_{cpt}$ (\cite{PS, p,606}). The cohomology groups of the resulting complexes  
$$
(L^{\cdot}_N,d_N):=(L^{\cdot}\cap d_N^{-1} L^{\cdot},d_N)\tag1.2
$$
and 
$$
(L^{\cdot}_D,d_D):=(L^{\cdot}\cap d_D^{-1} L^{\cdot},d_D)\tag1.3
$$
are denoted 
$$
H^*_N(V)\qquad\text{and}\qquad H^*_D(V)\tag1.4
$$
respectively, and are called  $L_2$-{\it de Rham}-{\it cohomology groups}. The operators  $d_N$ and $d_D$ are the maximal and minimal closed extensions of $d_{cpt}$. Others are possible, but it turns out (\cite{O}) that $H^*_N(V)$ and $H^*_D(V)$ are canonically isomorphic to the intersection cohomology groups $IH^*(V)$ when $V$ has only isolated singularities, so that in this case (and probably in general) all choices of boundary conditions (closed extensions of $d_{cpt}$) yield the same de Rham $L_2$-cohomology groups, which we denote
$$
H^*_2(V)
$$ 

Now the definitions (1.2)-(1.4) for de Rham complexes and cohomology work just as well for the corresponding $\bar\partial$-complexes and cohomology, giving for each $p$, $0\leq p\leq\dim V$, complexes 
$$
(L^{p,\cdot}_N,\bar\partial_N),\quad (L_D^{p,\cdot},\bar\partial_D)\tag1.5
$$
and $L_2-\bar\partial$-{\it cohomology groups}
$$
H^{p,*}_N(V),\quad H^{p,*}_D(V)\tag1.6
$$
However, as was already noted in \cite{P, (4.13)},\cite{PS, Theorem B}, unlike the de Rham groups, these {\it are} sensitive to changes in the boundary conditions. In order to state our first main result about Hodge structures, we introduce another $\db$-domplex, which mixes Dirichlet and Neumann boundary conditions, and has cohomology in general different from that of either of the $\db$-complexes above. Namely, we define for each $p$, $0\leq p\leq\dim V$,
$$
(\hat\C L^{p,\cdot},\Hat {\Bar\partial}):=(\C L^{\cdot}\cap \Hat {\Bar\partial}^{-1}\C L^{\cdot},\Hat {\Bar\partial})\tag1.7
$$
where
$$
\Hat{\Bar\partial}^{p,q}:=\left\{
\aligned \bar\partial_D^{p,q}\,,\,\,&p+q<n,\\
         \bar\partial_N^{p,q}\,,\,\,&p+q\geq n \endaligned
\right.\tag1.8
$$
and let 
$$
H^{p,*}_{D/N}(V)\tag1.9
$$
denote its cohomology groups. 

We can now state our first result concerning Hodge structure; for completeness we first give the standard definitions:
\medskip
{\rm 1.10.} {\bf Definition} Let $A$ be a subring of $\Bbb R$ such that
$A\otimes \Bbb Q$ is a field. An $A$-{\it Hodge structure of weight k} is a quadruple $(H_A;H_{\Bbb C},F^{\cdot};i)$, where $H_A$ is a finitely generated $A$-module, $H_{\Bbb C}$ is a $\Bbb C$-vector space, $i:H_A\otimes \Bbb C @>\cong>> H_{\Bbb C}$ is an isomorphism, and $F^{\cdot}$ is a decreasing filtration of $H_{\Bbb C}$ such that $F^0H_{\Bbb C}=H_{\Bbb C}$, $F^{k+1}H_{\Bbb C}=0$ and $H_{\Bbb C}=F^pH_{\Bbb C}\oplus \overline{F^{k-p+1}H_{\Bbb C}}$, for $p=0,1,\dots,k$. We call the above data {\it an A-Hodge structure on} $H_{\Bbb C}$ and say $H_{\Bbb C}$ {\it has an A-Hodge structure.} We define the {\it Hodge} $(p,q)$-{\it component} $H^{p,q}$ to be $F^pH_{\Bbb C}\cap\overline{F^{k-p}H_{\Bbb C}}$; it then follows that $H_{\Bbb C}=\oplus H^{p,q}$, the direct sum of its $(p,q)$-components, $H^{p,q}@>\cong>>F^p/F^{p+1}$, and $\overline{H^{p,q}}=H^{q,p}$ for all $p$. \medskip
Equivalently, a Hodge structure is such a direct sum decomposition; and the filtration is recovered by setting $F^i=\oplus_{p\geq i} H^{p,q}$.
\proclaim {\rm 1.11}{\bf Theorem A} Let $V$ be an $n$-dimensional complex projective variety with isolated singularities. Then for each $k$, $0\leq k \leq n$, there is a canonical isomorphism
$$
H^k_2(V)\cong\oplus_{p+q=k}H^{p,q}_{D/N}(V)
$$
so that $(IH^k(V);H^k_2(V),i;F^\cdot)$ is a $\Bbb Z$-Hodge structure of weight $k$,
where $i:IH^k(V)\otimes\Bbb C\to H^k_2(V)$ is the canonical isomorphism (\cite{O2}) and $F^\cdot$ is filtration by holomorphic degree.
\endproclaim
This theorem is proved by decomposing the space $\C H^k_2(V)$ of harmonic representatives of $H^k_2(V)$ into a direct sum of spaces $\C H^{p,q}_{D/N}(V)$ of harmonic representatives for $H^{p,q}_{D/N}(V)$. In the case of {\it complete } (in particular, compact) K\"ahler manifolds, the $d$-Laplace operator itself decomposes into a sum of $\db$-Laplace operators(\cite{Z1}); this is the usual route to the Hodge decomposition above. In our incomplete case, however, it will turn out that only in degrees $k$, $|n-k|\geq 2$, does such an 
operator decomposition hold. This is connected to the sensitivity of $\db$ to boundary conditions and will be discussed in detail in $\S 2$, where Theorem A is proved. In addition we will prove there that other parts of the "K\"ahler package'' hold, as conjectured in \cite{CGM}: the Lefschetz decomposition (``Hard Lefschetz") and the polarization of the Hodge structure on the primitive subspaces of $L_2$-cohomology.

 The second main theorem of this paper works in a context which is a sheafification of the definition of Hodge strucure above and is proved by a combination of algebraic and analytic methods. Here we get the Hodge strucure by identifying the filtered complex of {\it sheaves} of $L_2$-forms with a filtered complex of sheaves whose global cohomology is known to have a Hodge structure in this sheaf-theoretic sense. The sheaf-theoretic definition of a Hodge structure is given next.

 Let $R$ be a commutative ring and let $X$ be a topological space. Let $\C D^b_R(X)$ denote the derived category of complexes of $R$-sheaves on $X$ that are bounded below, and $\C D\C F^b_R(X)$, the corresponding derived category of filtered complexes.
\medskip
1.12. {\bf Definition} (\cite {D, (8.1.2)}) Let $A$ be a subring of $\Bbb R$ such that
$A\otimes \Bbb Q$ is a field. Let $X$ be a topological space. An $A$-{\it cohomological Hodge complex} is a quadruple $(\C K_A;\C K_{\Bbb C},\C F^{\cdot};\alpha)$ where
\roster

\ite a  $\C K_A$ is an object in $\C D^b_A(X)$,

\ite b  $(\C K_{\Bbb C},F^{\cdot})$ is an object in $\C D\C F^b_{\Bbb C}(X)$,

\ite c $\alpha$ is an isomorphism $ \C K_A\otimes\Bbb C @>\cong>> \C K_{\Bbb C}$ in $\C D^b_{\Bbb C}(X)$ and 

\ite d For each $k$ and $p$ the map $H^k(R\Gamma(X,\C F^p\C K_{\Bbb C}))\to H^k(R\Gamma(X,\C K_{\Bbb C}))$ is injective and the quadruple $(H^k(R\Gamma(X,\C K_A));H^k(R\Gamma(X,\C K_{\Bbb C})),F^{\cdot};H^k(R\Gamma(X,\alpha)))$ is an $A$-Hodge structure of weight $k$ on $H^k(R\Gamma(X,\C K_{\Bbb C}))$, where 
$$
F^pH^k(R\Gamma(X,\C K_{\Bbb C})):=\im (H^k(R\Gamma(X,\C F^p\C K_{\Bbb C}))\to H^k(R\Gamma(X,\C K_{\Bbb C}))).
$$
\endroster
\medskip
We call the above data an $A$-{\it cohomological Hodge structure on} $\C K_{\Bbb C}$. It follows from condition {\it d.} that the spectral sequence of the filtered complex $(R\Gamma(X,\C K_{\Bbb C})),F^{\cdot})$ collapses at the $E^1$-term and that the induced map $$
H^k(R\Gamma(X,gr^p_{\C F}\C K_{\Bbb C})\to F^p/F^{p+1}@<<\cong< H^{p,q}
$$
is an isomorphism.  
\medskip
1.13. {\bf Remark:} The vector spaces in {\it d.} above are hypercohomology. But since all the sheaves we use are fine, this is the same as the cohomology of the global sections. 
\medskip
1.14. {\bf Example:} Let $\pi:({\ti V},E)\to(V,\SV)$ be a resolution of singularities of a complex projective variety $V$ with isolated singularities $\SV$, where $E$ is a divisor with normal crossings. Define a complex of sheaves $\hat\C A^{\cdot}$ on $V$ by
$$\gather
\hat{\C A}^{\cdot} =\pi_*\Cal A_{\ti V}^0(\log E) @>d>> \pi_*\Cal A_{\ti V}^1(\log E) @>d>> 
\cdots @>d>> \pi_*\Cal A_{\ti V}^{n-2}(\log E) \to  \\
\to\pi_*\{\phi\in\Cal A_{\ti V}^{n-1}(\log E)|
d\phi\in\pi_*(\Cal I_E\Cal A_{\ti V}^n(\log E))\} @>d>> \\
\to\pi_*(\Cal I_E\Cal A_{\ti V}^n(\log E)) @>d>> 
\pi_*(\Cal I_E\Cal A_{\ti V}^{n+1}(\log E))\to\cdots @>d>> \pi_*(\Cal I_E\Cal A_{\ti V}^{2n}(\log E))
\endgather
$$
which we abbreviate to
$$
\hat{\C A}^{\cdot}:=\left\{
\aligned &\pi_*(\C A^k_{\ti V}(\lE),\,\,k<n,\\
         &\pi_*(\C I_E\C A^k_{\ti V}(\lE),\,\,k\geq n \endaligned
\right.\tag1.15
$$
where $\C A_{\ti V}^k(\lE)$ is the sheaf of $k$-forms on $\tilde V$ with at worst logarithmic poles along $E$ and $\Cal I_E$ is the ideal sheaf of $E$. In \cite{HZ, Z2}, Hain and Zucker noticed that $\hat \C  A^{\cdot}$ satisfies the axioms (\cite{GM}) for intersection cohomology: for small $U$ containing a singular point $v$ of $V$, the complex of vector spaces $\Gamma (U;\pi_*(\C A^\cdot_{\ti V}(\lE))$ computes the cohomology of $U-v$ while $\Gamma (U;\pi_*(\C I_E\C A^\cdot_{\ti V}(\lE))$ computes that of $U$. which matches the local computation of intersection cohomology. Hence, if $\C I\C H_{\Bbb Z}$ denotes the complex of $\Bbb Z$-sheaves of intersection cochains on $V$ with middle perversity, there is a unique isomorphism (\cite{GM}) $\alpha:\C I\C H_{\Bbb Z}\otimes \Bbb C @>\cong>> \hat\C A^{\cdot}$ in $\C D^b_{\Bbb C}(X)$. It was known that the mixed Hodge structure on $IH^*(V)$ was pure, and they also showed that the usual filtration by holomorphic degree induces a filtration $\C F^{\cdot}$ on $\hat \C  A^{\cdot}$ so that the quadruple $(\C I\C H_{\Bbb Z}; \hat \C A^{\cdot},\C F^{\cdot};\alpha)$ is a $\Bbb Z$-cohomological Hodge complex. In particular, we get a $\Bbb Z$-Hodge structure of weight $k$ on $IH^k(V;\Bbb C)=H^k(V;\hat \C A^{\cdot})$ for all $k$.
\medskip
1.16 {\bf Remark} One may also define for each $p$ the corresponding $\bar\partial$-complex
$$\gather
\hat{\C A}^{p,\cdot} =\pi_*\Cal A_{\ti V}^{p,0}(\log E) @>\bar\partial>> \pi_*\Cal A_{\ti V}^{p,1}(\log E) @>\bar\partial>> 
\cdots @>\bar\partial>> \pi_*\Cal A_{\ti V}^{p,n-p-2}(\log E) \to  \\
\to\pi_*\{\phi\in\Cal A_{\ti V}^{p,n-p-1}(\log E)|
\bar\partial\phi\in\pi_*(\Cal I_E\Cal A_{\ti V}^{p,n-p}(\log E))\} @>\bar\partial>> \\
\to\pi_*(\Cal I_E\Cal A_{\ti V}^{p,n-p}(\log E)) @>\bar\partial>> 
\pi_*(\Cal I_E\Cal A_{\ti V}^{p,n-p+1}(\log E))\to\cdots @>\bar\partial>> \pi_*(\Cal I_E\Cal A_{\ti V}^{p,2n-p}(\log E))
\tag1.17
\endgather
$$
and one expects that the $(p,q)$ components of the Hodge structure on $H^k(V;\hat\C A^{\cdot})$, the spaces  $H^k(V;gr^q_{\C F}\hat\C A^{\cdot})$, will be $H^q(V;\hat\C A^{p,\cdot})$ for $p+q=k$. This is indeed the case when $n=2$, but requires some proof, because the complexes $gr^p_{\C F}\hat\C A^{\cdot}$ and $\hat\C A^{p,\cdot}$ are not the same: for instance, a form $\phi$ of type $(p,n-p-1)$ in $\hat\C A^{n-1}$ must satisfy $\partial\phi\in \pi_*(\Cal I_E\Cal A_{\ti V}^{p+1,n-p-1}(\log E))$, whereas no such condition is required for membership in $\hat\C A^{p,n-p-1}$. So it must be shown that the canonical map of complexes $\hat\kappa^p:gr^p_{\C F}\hat\C A^{\cdot} \to \hat\C A^{p,\cdot}$ is a quasi-isomorphism. This issue is treated in $\S 4$.
\medskip

As with the sheafification of the notion of a Hodge structure in (1.12) above, we must sheafify $d_N$ and $d_D$ by defining operators $d_N(U)$ and $d_D(U)$, for each open set $U\subseteq V$,  which equal those defined above 
in case $U=V$. For $d_N$ this was done in \cite{PS, p.606} and for $d_D$ the definition is similar: for each open set $U$, set
$$
d_N(U):=d_w(U-U\cap \SV)\,\,\text{and} \,\,d_D(U):=d_w(U)\tag1.18
$$
where $d_w(U-U\cap\SV)$ (resp. $d_w(U)$) denotes the weak derivative with
respect to compact subsets of $U-U\cap\SV$ (resp., compact subsets of $U$). Finally, generalizing (1.2) and (1.3), we define complexes of sheaves on $V$,
$$
(\C L^{\cdot}_N,d_N):=(\C L^{\cdot}\cap d_N^{-1}\C L^{\cdot},d_N)\tag1.19
$$
and 
$$
(\C L^{\cdot}_D,d_D):=(\C L^{\cdot}\cap d_D^{-1}\C L^{\cdot},d_D)\tag1.20
$$

Now the results of Ohsawa cited above are actually local, so that each of these complexes of sheaves is isomorphic in the derived category $\C D^b_{\Bbb C}(V)$ to the middle-perversity intersection complex $\C I\C H_{\Bbb C}$. Also, each admits the standard filtration $\C F^{\cdot}$ by holomorphic degree; but neither of these filtered complexes of sheaves will easily produce the associated gradeds we obtained in Theorem A. Rather, another incarnation of $\C I\C H_{\Bbb C}$ will be used, one which mixes the Neumann and Dirichlet boundary conditions as was done in Theorem A. Namely, we define
$$
\hat d^k:=\left\{
\aligned d_D\,\,&k<n,\\
         d_N\,\,&k\geq n \endaligned
\right.\tag1.21
$$   
and then
$$
(\hat\C L^{\cdot},\hat d):=(\C L^{\cdot}\cap \hat d^{-1}\C L^{\cdot},\hat d)\tag1.22
$$
Similarly we have the corresponding $\bar\partial$-complexes, $(\C L^{p,\cdot},\bar\partial_N)$, $(\C L^{p,\cdot},\bar\partial_D)$ and 
$$
(\hat\C L^{p,\cdot},\Hat {\Bar\partial}):=(\C L^{\cdot}\cap \Hat {\Bar\partial}^{-1}\C L^{\cdot},\Hat {\Bar\partial})\tag1.23
$$
defined for each $p$ by
$$
\Hat{\Bar\partial}^{p,q}:=\left\{
\aligned \bar\partial_D^{p,q}\,,\,\,&p+q<n,\\
         \bar\partial_N^{p,q}\,,\,\,&p+q\geq n; \endaligned
\right.\tag1.24
$$

We can now state the second main theorem of this paper. 
\proclaim {\rm 1.25}{\bf Theorem B} Let $V$ be a complex projective variety of dimension two with at most isolated singularities. Then there is
\roster
\ite a a filtered isomorphism 
$$
\hat\Lambda:(\hat\C L^{\cdot},\C F^{\cdot}) @>\cong>> (\hat\C A^{\cdot},\C F^{\cdot})
$$
in $\C D\C F^b_{\Bbb C}(V)$ and
\ite b for each $p\geq 0$, a canonical isomorphism in $\C D^b_{\Bbb C}(V)$
$$
\hat\kappa:gr^p_{\C F}\hat\C L^{\cdot} @>\cong>> \hat\C L^{p,\cdot}$$
\endroster
\endproclaim

\proclaim {\rm 1.26}{\bf Corollary} For $V$ as above, the quadruple $(\C I\C H_{\Bbb Z};\hat\C L^{\cdot},\C F^{\cdot};\beta)$ is a $\Bbb Z$-cohomological Hodge complex, isomorphic to the Hain-Zucker $\Bbb Z$-cohomological Hodge complex $(\C I\C H_{\Bbb Z};\hat\C A^{\cdot},\C F^{\cdot};\alpha)$ and for each $k$, the Hodge components of $H^k(\hat\C L^{\cdot})$ are canonically isomorphic to $H^q(\hat\C L^{p,\cdot})$, $p+q=k$. In particular, the isomorphism
$H^k(V;\hat\C  A^{\cdot}) @>\cong>> H^k(\hat\C L^{\cdot})$, induced by the canonical isomorphism $\Lambda$ preserves the respective $(p,q)$-components of the Hodge structures and induces isomorphisms
$$
H^q(V;\hat\C A^{p,\cdot}) @>\cong>> H^q(V;\hat\C L^{p,\cdot})
$$
\endproclaim
{\bf Proof} By uniqueness of $\C I\C H_{\Bbb C}$ in $\C D^b_{\Bbb C}$, there is a commutative diagram of isomorphisms in $\C D^b_{\Bbb C}(V)$
$$
\CD
\C I\C H_{\Bbb C} @>\alpha>> \hat\C L^{\cdot}\\
@V=VV                         @VV{\hat\Lambda}V   \\
\C I\C H_{\Bbb C} @>\beta>> \hat\C A^{\cdot}
\endCD\tag1.27
$$
in $\C D^b_{\Bbb C}$. It is immediate from the the hypercohomology spectral sequences of
$ \hat\C  A^{\cdot},\,\C F^p\hat\C  A^{\cdot},\, \hat\C L^{\cdot}$ and $\C F^p\hat\C L^{\cdot}$ that $\Lambda$ induces vertical isomorphisms in a commutative diagram
$$
\CD
H^k(V;\C F^p\hat\C  A^{\cdot}) @>>> H^k(V;\hat\C  A^{\cdot})\\
@V\cong VV                               @VV\cong V  \\
H^k(V;\C F^p\hat\C L^{\cdot}) @>>> H^k(V;\hat\C L^{\cdot})
\endCD
$$
for each $p=0,\dots ,n$ and $k=0,\dots,2n$. By the result of Hain and Zucker, the top horizontal is injective, so the bottom one is as well. Finally. the commutativity of (1.27) shows that the right vertical isomorphism in the above diagram preserves the underlying real structures coming from $\alpha$ and $\beta$, so we are done.

\newpage

\heading \S 2. The pure Hodge structure for varieties with isolated singularities\endheading

 Let $V$ be a variety with isolated singularities; to simplify notation, we assume there is only one singular point. Fix an embedding of $V$ in $\Bbb P^N$ and coordinates 
$(z_1,\dots,z_N)$ on the complement $\Bbb C^N$ of a hyperplane so that the image of the singular point is the origin. Let $U=V\cap\Bbb C^N$. Then the restriction to $U-\{0\}$ of the Fubini-Study metric on $\Bbb P^N$ has K\"ahler form 
$$\omega : = i\da\db \log(1+r^2)/2\tag2.1$$
where $r^2=\sum |z_i|^2$.    
Unless otherwise stated, the pointwise inner product and norm
$$
\langle\xi,\eta\rangle\,\,\text{and}\,\,|\xi|:=\langle\xi,\xi\rangle
^\frac{1}{2}
$$
and the (global) $L_2$-inner product and norm
$$
(\xi,\eta):=\int_{U-\{0\}}\langle\xi,\eta\rangle\,dU \quad\text{and}\quad\|\xi\|:=\left (\int_{U-\{0\}}\langle\xi,\xi\rangle\,dU \right)^{1/2}
$$
of $(p,q)$-forms $\xi$ and $\eta$ on $V-\{0\}$ will be with respect to the Hermitian inner product defined by this metric, where $dU$ denotes its volume form: since $V-U$ has measure zero and is supported away from the singular point, these integrals equal their counterparts over $V$.  We say $\xi$ {\it is} $L_2$ when $\|\xi\|<\infty$ ; the space of locally $L_2$ forms of type $(p,q)$ on any subspace $W\subseteq V$ is denoted $L^{p,q}(W)$; similarly, that of $k$-forms is denoted $L^k(W)$. 
   %Set 
   %$$
   %U_{\f12}:=\{z\in V|  \|z\|\leq c\}.
   %$$

   %By \cite{Milnor, Cor. 2.9}, 
   %$dr$ is nonzero on $U_{\f12}-\{0\}$ for sufficiently small $c>0$; we fix 
   %such a $c$ with $c < 1.$ 

Since $dr$, viewed as a 1-form on $\Bbb C^N$, has norm $1+r^2$ with respect to the Fubini-Study metric, its restriction to $U$ has norm 
$$
|dr|\leq 1+r^2.\tag2.2
$$
Hence if we set 
$$u = \tan^{-1}(r),\tag2.3$$
then 
$$
|du|\leq 1\tag2.4$$ 
on $U$ and of course $|\da u| =|\db u|=2^{-\frac{1}{2}}|du|$.

For any form $\alpha$ we will denote by $e(\alpha)$ (resp., $e^*(\alpha)$) the operation of exterior (resp., interior) multiplication on the left by $\alpha$; in case $\alpha$ is the K\"ahler form $\omega$, we denote these operators as usual by $L$ and $\Lambda$. We will often use the fact that that if $\phi$ is a differential form, then
$$
|e(\alpha)\phi|\leq|\alpha|\cdot|\phi|\,\,\text{and}\,\,|e^*(\alpha)\phi|\leq|\alpha|\cdot|\phi|\tag2.5
$$ 
%The following standard facts will also be used frequently. For any $k$-forms %$\alpha$ and $\beta$, \roster
%\item $e^*(\alpha\wedge\beta)=e^*(\beta)e^*(\alpha)$, 
%\item $e^*(f\alpha)=\bar fe^*(\alpha)$ for any function $f$,
%\item $e^*(\alpha)e(\beta)+e(\beta)e^*(\alpha)=<\alpha,\beta>$, if $\alpha$ and %$\beta$ are 1-forms.
%\endroster
On $U$ we may uniquely express any form $\phi$ as :
$$
\phi = \phi_0 + e(\frac{\db u}{|\db u|})\phi_1 + 
e(\frac{\partial u}{|\da u|})\phi_2 +
 e(\frac{\db u}{|\db u|})e(\frac{\partial u}{|\da u|})\phi_3,\tag2.6
$$
where each $\phi_i$ is in the kernel of $e^*(\db u)$ and $e^*(\da u)$, and 
$\phi_i = 0$, for $i>0$ at any point where $du$ vanishes. Then
$$
(\phi,\phi)= (\phi_0,\phi_0)+(\phi_1,\phi_1)+(\phi_2,\phi_2)+(\phi_3,\phi_3).\tag2.7 
$$
Because the K\"ahler form is $\db$-exact on $U$, we may express $L$ there as
$$
 L = -i\{\db, e(\partial \ln(1+r^2)/2)\}=e(-i\db \partial \ln(1+r^2)/2),
$$
where $\{\cdot, \cdot\}$ denotes the anticommutator.
More generally, given a smooth function $f$, we have 

$$  
-i\{\db, e(\partial f(\ln(1+r^2))/2)\} = 
-i2r^2f''(\ln(1+r^2))e(\db u)e(\partial u) + f'(\ln(1+r^2))L.\tag2.8
$$
  We now use this identity to relate certain weighted $L_2$ norms of 
a form $\phi$ on $U-\{0\}$ to the $L_2$ norms of $\db\phi$,  $\vartheta\phi$, 
 and $\phi|V-U_{\f12}$, with  
$$
U_{\f12}:=\{z\in V| \, 0<r\leq \f12\}.\tag2.9
$$
Because of the form of the identity (2.8), our weight functions will be continuous, piecewise smooth functions of $t$, where 
$$
t = \ln(1+r^2).\tag2.10
$$
Thus, for small $r$, $t = r^2 + O(r^{4})$; when convenient in subsequent 
computations, we will replace $t$ by $r^2$ and so introduce the $O(r^2)$-error term. In addition, statements like ``$\p/r\log(1/r^2)$ is $L_2$" and  ``$\p/t^\f12\log(1/t)$ is $L_2$" will be used interchangeably. We will also use 
$$
dt=2r\,du\quad\text{and so}\quad|dt|=(2t^{\f12}+O(t))|du|\tag2.11
$$

Let $D$ and $D'$ denote the operators
$$
D:=\db+\vartheta\qquad\qquad D':=\da+\bar\vartheta
$$
\proclaim{\rm 2.12}{\bf Convention} When such operators are used without subscripts $D$ or $N$ indicating boundary conditions, we always mean the weak derivatives.
\endproclaim
The basis for our estimates will be variations of the following 
proposition which is a variation of a computation of \cite{DF}.  

\proclaim{\rm 2.13}{\bf Proposition} Let $\phi$ be form of type $(p,q)$ on $V$, 
and let
 $f:(0,\infty)\to \Bbb R$ be a piecewise smooth (``weight") function with 
$f'$ continuous and vanishing for sufficiently large $r$. Suppose either
that $\phi$ is supported away from $0\in V$ or that $f'$ is supported away from 0. Then we have
 $$
\multline
(n-p-q)(f'(t)\phi,\phi) 
- (i2r^2f''(t)[e^*(\db u)e^*(\partial u)/|\db u|^2,e(\db u) e(\partial u)]\phi,\phi) \\
= -(D\phi,rf'(t)[e(\db u) + e^*(\db u)]\phi) 
 -(D'\phi,rf'(t)[e(\da u) + e^*(\da u)]\phi). 
\endmultline$$
\endproclaim

{\bf Proof:} The proof is a simple computation involving only integration by
 parts and 
K\"ahler identities. Starting from (2.8) above, integrate by parts ($f'$ is supported away from $\infty$ and either $\phi$ or $f'$ is
 supported away from 0) to get the second equality 
below and use the K\"ahler identity $i\da=[\vartheta,L]$ to get the fourth in

$$\aligned     
(-i2&r^2f''(t)e(\db u)e(\partial u)\phi +  f'(t)L\phi,L\phi)\\
&= (-i\{\db, e(\partial f(t)/2)\} \phi,L\phi)  
=(-ie(\partial f(t)/2) \db \phi,L\phi) -(i e(\partial f(t)/2)\phi,\vartheta L\phi)\\  
&=(-ie(\partial f(t)/2) \db \phi,L\phi) -
 (i e(\partial f(t)/2) \phi,[\vartheta,L]\phi) - (i e(\partial f(t)/2) \phi,L\vartheta\phi)  \\
 &=(-ie(\partial f(t)/2) \db \phi,L\phi) -
 (i e(\partial f(t)/2) \phi,i\partial\phi) - 
(ie(\partial f(t)/2) \phi,L\vartheta\phi). \\
\endaligned
$$ 

We reorganize this as 

$$
\multline
(f'(t)L\phi,L\phi) = 
 (-ie(\partial f(t)/2) \db \phi,L\phi) - 
(ie(\partial f(t)/2) \phi,L\vartheta\phi)
 \\ -(e(\partial f(t)/2) \phi,\partial\phi) +
 (i2r^2f''(t)e(\db u)e(\partial u)\phi,L\phi).
\endmultline\tag2.14
$$

Similarly, applying the identity (2.8) to $\Lambda\p$, taking inner product with $\p$, and using the K\"ahler identity $-\bar\vartheta=[\Lambda,\db]$, we have 

$$
\multline
(f'(t)L\Lambda\phi,\phi) \\
=-(ie(\partial f(t)/2) \Lambda\phi,\vartheta\phi) - 
(ie(\partial f(t)/2)\db \Lambda\phi,\phi)
+(i2r^2f''(t)e(\db u)e(\partial u)\Lambda\phi,\phi) \\= 
(-ie(\partial f(t)/2) \Lambda\phi,\vartheta\phi)
 - (ie(\partial f(t)/2) \Lambda\db\phi,\phi) \\
+ (e(\partial f(t)/2) \bar\vartheta\phi,\phi)
+ (i2r^2f''(t)e(\db u)e(\partial u)\Lambda\phi,\phi). 
\endmultline
$$ 

Subtracting this equality from (2.14) and using the equality 
$\da f(t)/2=rf'(t)\da u$ gives
 %e^*(\db u)e^*(\partial u)/|\db u|^2-( %\db \phi,e(\db f(t))/2)\phi)  
%- (e^*(\db f(t))/2\phi,\vartheta\phi) 

$$
\multline
(f'(t)[\Lambda,L]\phi,\phi)  \\
= (irf'(t) [e(\da u),\Lambda]\db\phi,\phi)  
  +(irf'(t) \db[e(\da u),\Lambda]\phi,\phi)\\ 
- (rf'(t)e(\partial u )\bar\vartheta\phi,\phi)
 -  (rf'(t)e(\partial u) \phi,\partial\phi)\\
 + (i2r^2f''((t)[\Lambda,e(\db u) e(\partial u)]\phi,\phi).
\endmultline 
$$
Then using the K\"ahler identity $[\Lambda,L]\phi=(n-p-q)\phi$ completes the proof. 
 
We use this proposition to bound weighted $L_2$ norms of $\phi$ by 
certain weighted norms of $D\phi$ and $D'\phi$.

\proclaim{\rm (2.15)}{\bf Corollary} With notation and assumptions as in Proposition, and for any function $F:(0,\infty)\to\Bbb R$, we have
$$
\multline
(n-p-q)(f'(t)\phi,\phi) + 
(|d u|^2r^2f''(t)\phi_0,\phi_0) 
-(|d u|^2 r^2f''(t)\phi_3,\phi_3)\\
 \leq(\|F(t)D\phi\|+\|F(t)D'\phi\|)\cdot \|2^{-\f12}|du|(rf'(t)/F(t))\phi\|.
\endmultline
$$
\endproclaim
{\bf Proof}  Multiply and divide by $F$ in the first and second terms on the right in (2.13) and use Cauchy-Schwarz. 
\medskip 
For the applications we have in mind, we will need to relax the support hypothesis on $\phi$ or $f'$ in this Proposition. To do this we choose a sequence $\{\phi_j\}$ of compactly supported forms converging to $\phi$ in the $L_2$ norm; or else a sequence $\{f'_j\}$ of compactly supported functions converging to $f'$. Then we let $j\to\infty$ in (2.15) or in estimates derived from it. In practice, we always choose the approximating sequences in the same way:
$$
\phi_j:=\mu_j(t)\phi\,\,\text{or else}\,\,f'_j=\mu_j(t) f'
$$
where $\mu_j$ is a smooth function $(0,1)\to [0,1]$ having the properties (see \cite{PS, 3.5})
$$
\mu_j(t)=\left\{
\aligned
&1,\,t\geq e^{-e^j}\\
&0,\,t\leq e^{-e^{j+1}}
\endaligned
\right.\tag2.16
$$
and 
$$
|\mu_j'(t)|\leq\frac{2\chi_{I(k)}}{t^{\f12}
\log(1/t)},\tag2.17
$$
where  $\chi_{I(k)}$ is the characteristic function of the interval $I(k)=[e^{-e^{k+1}},e^{-e^k}]$. 

For instance, assume $\p$ and $d\p$ are $L_2$. Clearly, $\mu_j\phi$ converges to $\phi$ in $L_2$. To show that $d\mu_j\phi\to d\p$ as well, we need to control $d\mu_j\wedge\p$. For this, and to widen the applicability of estimates like (2.15) to include non-compactly supported forms, the following obvious lemma will be used. 

\proclaim{\rm 2.18}{\bf Lemma} Let $\xi$ be an $m$-form on $U-\{0\}$, $m\geq 0$, and suppose that $\xi/r\log(1/r^2)\in L^m(U_\f12)$. Then 
$$
\| d\mu_j\wedge\xi\|\to 0
$$
as $j\to\infty$. In particular, if $\xi\in\dom d_N$, then $\xi\in\dom d_D$.
\endproclaim

{\bf Proof} By (2.4) and (2.11) $|dt|\leq Kt^{\f12}$ on $U_\f12$ for some $K>0$, so
$$
\| d\mu_j(t)\wedge\xi\| = 
\|t^{\f12}\log(1/t)d\mu_j(t)\wedge(\xi/t^{\f12}\log(1/t))\|\leq 
2K\|\chi_{I(j)}\xi/t^{\f12}
\log(1/t)\|^2\rightarrow 0
$$
by the Lebesgue dominated convergence theorem. The second 
assertion is an immediate consequence of this.
 %and the identity
%$$
%d\mu_j\xi=\mu_jd\xi+d\mu_j\wedge\xi.
%$$
\medskip
\noindent{\rm 2.19} {\bf Remark} Identical results hold where $d$ is replaced by $\db$, $D$, etc. Thus for such $\xi$ boundary conditions are irrelevant and we may drop the subscript $B$ indicating a boundary condition from $d_B\xi$. Moreover, standard consequences of the K\"ahler condition, which hold for forms compactly supported in $V-0$, are valid for such $\xi$. For instance:
\proclaim{\rm 2.20}{\bf Proposition} Let $\xi$ be an $m$-form on $V-\{0\}$ such that $\xi/r\log(1/r^2)\in L^m(U_\f12)$ and $D\xi$ is $L_2$. Then $D'\xi$ is $L_2$ and $\|D\xi\|=\|D'\xi\|$.
\endproclaim 
\proclaim{\rm 2.21}{\bf Proposition} Let $\phi\in L^{p,q}(V)$ where $n-p-q\geq 1$. Assume that $D\phi$ and $D'\phi$ are $L_2$. Then $\phi/r\log(1/r^2)\in L^{p,q}(U_\f12)$.
\endproclaim
{\bf Proof:} Let $k\in \Bbb R$ and in (2.19) take $F(t)\equiv 1$; and for $t < \f12$
$$
f'(t)=\mu_j(t)\log^k(1/t)
$$
where $\mu_j$ is the cut-off function recalled in (2.16) and $f'$ is extended to a smooth bounded function for $t\geq\f12$.
Then we obtain, for some positive constants $C$ and $K$, the inequality
$$
\multline
(\chi_{\f12}(t)\mu_j(t)\log^k(1/t)\phi,\phi) + ([r^2\mu'_j(t)\log^k(1/t)-k\chi_{\f12}(t)\mu_j(t)\log^{k-1}(1/t)]|du|^2\phi_0,\phi_0)\\
 - ([r^2\mu'_j(t)\log^k(1/t)-k\chi_{\f12}(t)\mu_j(t)\log^{k-1}(1/t)]|du|^2\phi_3,\phi_3)\\
 \leq C\|r\chi_{\f12}(t)\mu_j(t)\log^k(1/t)\phi\|+K\|(1-\chi_{\f12}(t))\p\|
\endmultline\tag2.22
$$
where $\chi_{\f12}(t)$ denotes the characteristic function of $U_{\f12}$.   
The right side is bounded in $j$ (non-uniformly in $k$). Assume that $\chi_{\f12}(t)\log^{(k-1)/2}(1/t)\p\in L^{p,q}_2(U_{\f12})$ for some $k$. Then as $j\to\infty$ the integrals $(r^2\mu^\pr_j(t)\log^k(1/t)\phi_i,\phi_i)$, $i$ = 0, 3, tend to zero by dominated convergence. Hence the integrals $(\chi_{\f12}(t)\mu_j(t)\log^k(1/t)\phi,\phi)$ are bounded as $j\to\infty$, so $\chi_{\f12}(t)\log^{k/2}(1/t)\p$ is $L_2$. Hence, beginning with negative $k$, we conclude by induction that $\chi_{\f12}(t)\log^{k}(1/t)\p$ is $L_2$ for all $k$.

Return to (2.15) and this time take for $t \leq \f12$
$$
f'(t)=\mu_j(t)\log^k(1/t),\qquad F(t)=r\log^{\f {k+1}2}(1/t)
$$
and extend them in a bounded fashion for $t\geq \f12$, and so that $F$ is nowhere zero.  
Now using $\chi_{\f12}(t)\log^k(1/t)\p\in L^{p,q}_2(U_{\f12})$, we argue as above that we may discard the
integrals $(r^2\mu^\pr_j(t)\log^k(1/t)\phi_i,\phi_i)$ and that we may replace $\mu_j(t)$ by 1. Then applying Cauchy-Schwarz to the first term on the right side and discarding a positive term on the left, we get, for any $\epsilon>0$,
$$
\multline
(\chi_{\f12}(t)\log^k(1/t)\phi,\phi)  -k\chi_{\f12}(t)\log^{k-1}(1/t)\phi_0,\phi_0)\\
 \leq \f 1{2\epsilon}(\|r\chi_{\f12}(t)\log^{\f {k+1}2}(1/t)D\phi\|^2+\|r\chi_{\f12}(t)\log^{\f {k+1}2}(1/t)D'\phi\|^2)+\frac{\epsilon}{2}\|\chi_{\f12}(t)\log^{\f {k-1}{2}}(1/t)\phi\|^2\\
+K(\|(1-\chi_{\f12}(t))\p\|+\|(1-\chi_{\f12}(t))D\p\|+\|(1-\chi_{\f12}(t))D'\p\|)
\endmultline\tag2.23
$$
We need the following simple telescoping series calculation:
$$
\multline
-1+\sum_{k=1}^{N+1}\frac{\log^{k-1}(1/t)}{(k+1)!} = \sum_{k=1}^N\bigl(\frac{\log^{k}(1/t)}{(k+1)!}-\frac{\log^{k-1}(1/t)}{k!}+\frac{\log^{k-1}(1/t)}{(k+1)!}\bigr)\\
=\sum_{k=1}^N(\frac{\log^{k}(1/t)}{(k+1)!}-
\frac{k\log^{k-1}(1/t)}{(k+1)!}\bigr).
\endmultline\tag2.24
$$

Now the left side of (2.23) is no larger if $\p_0$ is replaced by $\p$ in the second term; then dividing both sides by $(k+1)!$, summing over all $k\geq 1$ gives, according to (2.23),
$$
\multline
-(\chi_{\f12}(t)\phi,\phi)+(\frac{\chi_{\f12}(t)}{\log^2(1/t)}(\frac{1}{t}-1-\log(1/t))\phi,\phi)\\
\leq
\f 1{2\epsilon}(r^2(\frac{1}{t}-1-\log(1/t))\chi_{\f12}(t)D\p,D\p)+\f 1{2\epsilon}(r^2(\frac{1}{t}-1-\log(1/t))\chi_{\f12}(t)D'\p,D'\p)\\
+\frac{\epsilon}{2}(\frac{\chi_{\f12}(t)}{\log^2(1/t)}(\frac{1}{t}-1-
\log(1/t))\phi,\phi)+K'
\endmultline
$$
for some positive $K'$.
Now take $\epsilon$ small and reorganize to get
$$
\multline
(1-\epsilon/2)\|(\chi_{\f12}(t)/t^{\f12}\log(1/t))\p\|^2\\
\leq \f 1{2\epsilon}(r^2(\frac{1}{t}-1-\log(1/t))\chi_{\f12}(t)D\p,D\p)+\f 1{2\epsilon}(r^2(\frac{1}{t}-1-\log(1/t))\chi_{\f12}(t)D'\p,D'\p)+(\phi,\phi)+K'.
\endmultline\tag2.25
$$
So $(1/r\log(1/t))\phi\in L^{p,q}(U_\f12)$.
\medskip
\noindent{\rm 2.26} {\bf Remark} We will use again below the device of summing estimates like (2.23) but will leave computations like (2.23) to the reader.
\medskip

When one considers forms of degree $k$ with $n-k\geq 2,$ 
difficulties with boundary conditions largely disappear as we will 
show with the following proposition.

\proclaim{\rm 2.27}{\bf Proposition} Let $\phi\in L^{p,q}(V)$ where $n-p-q\geq 2$. Assume 
that $D\phi$ is $L_2$. Then $\phi/r$ is $ L_2$, and 
$D'\phi$ is $ L_2$. Similarly, if $\phi\in L^{k}(V)$ where $n-k\geq 2$ and $(d+\delta)\phi$ is $ L_2$, then $\phi/r$ is $ L_2$. 
\endproclaim
{\bf Proof:} We prove the first assertion. The proof of the 
second is identical.
The idea of the proof is to use the K\"ahler identities so that we can bound 
$D'\phi$ and thus reduce to the preceding corollary. This of course requires integration by parts. The difficulties arise in justifying this integration.

Reflecting our less restrictive hypotheses on the exterior derivatives of $\p$, we 
 recast our basic estimate (2.15) by applying Cauchy-Schwarz only to the first term on the right side of (2.13): for any compactly supprted form $\psi$, we get
$$
\multline
(n-p-q)(f'\psi,\psi) 
+ (r^2|d u|^2f''\psi_0,\psi_0) - (r^2|d u|^2f''\psi_3,\psi_3) \\
\leq \|FD\psi\| \|rf'|\da u|\psi/F\|
 - (FD'\psi,rf'[e(\da u) + e^*(\da u)]\psi/F).\endmultline\tag2.28
$$
where, to shorten notation, we write $F$ for $F(t)$, $f'$ for $f'(t)$, etc.  
Now estimate the last term by 
$$
\multline
-(FD'\psi,rf'[e(\da u) + e^*(\da u)]\psi/F)  \\ 
 =  -(D'F\psi,rf'[e(\da u) + e^*(\da u)]\psi/F) + \\
(2rF'[e(\da u) - e^*(\da u)]\psi, rf'[e(\da u) + e^*(\da u)]\psi/F)\\  \leq
\|D'F\psi\| \||\da u|rf'\psi/F\| +
(r|du|^2F' \psi_0, rf' \psi_0/F)
+ (r|du|^2 F'\psi_1, rf' \psi_1/F)\\- (r|du|^2F'\psi_2, rf'\psi_2/F) 
 - (r|du|^2F'\psi_3, rf'\psi_3/F).
\endmultline \tag2.29 
$$
Using the fact that $\|DF\psi\|=\|D'F\psi\|$ (since the metric is K\"ahler and $\psi$ is compactly supported) and 
some elementary manipulations gives 
$$
\multline
-(FD'\psi,rf'/F(e(\da u) + e^*(\da u))\psi) \leq 
\\
(\|FD\psi\| + \|2rF'|\da u|\psi\|) \||\da u|rf'\psi/F\|  
+(r|du|^2F' \psi_0, rf' \psi_0/F)
+ (r|du|^2F' \psi_1, rf' \psi_1/F)\\- (r|du|^2F'\psi_2, rf'\psi_2/F) 
 - (r|du|^2F'\psi_3, rf'\psi_3/F).
\endmultline\tag2.30
$$

Let us now make the simplifying and useful choice  
$$f' = F^2/r^2,$$
%2FF' = f''r^2 + (1+r^2)f' = r^2f'' + (1+r^2)F^2/r^2.
Inserting this into our earlier inequality (2.28) and using (2.30) and (2.10) we get  
$$ 
\multline
([(n-p-q - |du|^2(r^2+1))F^2/r^2 + |du|^2FF']\psi_0,\psi_0) 
+ ([(n-p-q)F^2/r^2 - |du|^2FF']\psi_1,\psi_1)\\
+ ([(n-p-q)F^2/r^2 + |du|^2FF']\psi_2,\psi_2)
+ ([(n-p-q + |du|^2(r^2+1))F^2/r^2 - |du|^2FF']\psi_3,\psi_3) \\
\leq (2^\f12\|FD\psi\| + \|rF'\psi\|)(\|F\psi/r\|  
\endmultline 
$$

Set now for $t<\f12$, 
$F(t) = F_k(t) = t^{\f12}\log^{\frac{k}{2}}(1/t)$, and extend in a bounded fashion for $t\geq \f12$, with $F_k(t)/\log^{\f k2}(2)$ uniformly bounded as $k\to\infty$.

%FF' = [\log^k(1/t) - kln(1/t)^{k-1}]/2
%F^2 = t\log^k(1/t)
%t^{\f12}F' = [\log^{\frac{k}{2}}(1/t) - kln(1/t)^{k/2-1}]/2
Plugging this in and using $t=r^2+O(r^4)$, gives 
for some positive constants $C$ and $K$,

%([((n-p-q - |du|^2(r^2+1 - r^2/2t))t\log^k(1/t)/r^2) - %|du|^2[kln(1/t)^{k-1}]/2]\psi_0,\psi_0)\\ 
%+ ([((n-p-q)t/r^2 - |du|^2)\log^k(1/t) + %|du|^2[kln(1/t)^{k-1}]/2]\psi_1,\psi_1)\\
%+ ([((n-p-q)t/r^2 + |du|^2)\log^k(1/t) - |du|^2k\log^{\frac{k-1}{2}}(1/t)\psi_2,\psi_2)
%+ ([(n-p-q + |du|^2(r^2+1 r^2/t))t\log^k(1/t)/r^2 + %|du|^2k\log^{\frac{k-1}{2}}(1/t)\psi_3,\psi_3) \\
%\leq \|FD\psi\| 2^{-1/2}\|F/r\psi\| + \\
%(\|FD\psi\|2^{-1/2}\|F/r\psi\| + \|rF'\psi\|\|F/r\psi\|, 
%\endmultline 

$$
\multline
([(n-p-q -\f12- Cr^2)\log^k(1/t) - k\log^{\frac{k-1}{2}}(1/t)/2]\chi_{\f12}(t)\psi,\psi)\\ 
\leq (2^{\f12}\|F_k D\psi\| + \|rF'_k\psi\|)\|F_k\psi/r\|+K\|(1-\chi_{\f12})F_k\psi\| .
\endmultline\tag2.31
$$
 
Applying Cauchy-Schwarz to two terms on the right gives  
$$
\multline
([(n-p-q - 1-Cr^2  - 1/2M)\log^k(1/t)-k^2\log^{k-2}(1/t)/4]\chi_{\f12}(t)\psi,\psi)\\ 
\leq  M\|F_kD\psi\|^2 +K\|(1-\chi_{\f12})F_k\psi\|.
\endmultline\tag2.32 
$$
For $M$ large, this is the estimate we need to prove the proposition. We may apply it 
inductively as follows. 

Suppose that $\log^{\frac{k}{2}-1}(1/t)\phi\in L_2$ 
for some $k$.   
Then $\|F_kD\mu_j\phi\|\rightarrow \|F_kD\phi\|$ (see the proof of Lemma (2.18)). Take $\psi = \mu_j\phi$ in (2.32); then the right side is bounded as $j\rightarrow \infty$ and therefore 
so too is the left. From this we deduce  
$\log^{\frac{k}{2}}(1/t)\phi\in L_2(U_{\f12})$, and by induction we obtain the result 
for all $k$. We may now divide (2.32) by $k!$, 
sum over $k\geq 0$ and, obtain for $M$ large, 
$$\|\chi_{\frac{1}{2}}\phi/t^\f12\|^2 \leq 2M \|D\p\|^2 + 4K\|(1-\chi_{\frac{1}{2}})\p\|^2+C\|rt^{-\f12}\chi_{\frac{1}{2}}\phi\|^2+2^{\f32}\|\log^{-\f12}(1/t)\chi_{\frac{1}{2}}\phi\|^2.$$
In particular, $\phi/r \in L_2(U_{\f12})$. 

Finally, by (2.20) $D'\p$ is $L_2$.

\proclaim{\rm 2.33}{\bf Corollary} Let $\phi\in L^{p,q}(V)$ where $n-p-q\geq 2$. Assume 
that $D\phi = 0$. Then $\phi/r^{n-p-q}\log(1/r^2)\in L_2(U_\f12)$.  
\endproclaim
{\bf Proof:} 
We know $\phi/r \in L_2(U_{\f12})$, so $D'\p=0$ by (2.20); and $\|\log^k(1/t)\mu_j\p\|\to\|\log^k(1/t)\p\|$ for any real $k$. Further $D\mu_j\phi$ and $D'\mu_j\phi$ both tend to zero in $L_2$ (see (2.18)), so we can take the limit over $j$ of the estimates (2.29) where $\psi=\mu_j\p$ to get 

$$
(n-p-q)(f'(t)\phi,\phi) + (r^2|d u|^2f''(t)\phi_0,\phi_0) 
- (r^2|du|^2f''(t)\phi_3,\phi_3) 
\leq 0.
$$   
As in the Proposition (see (2.31)), we obtain for some $C$ and $K$ the estimate  $$
([(n-p-q - Cr^2)\log^k(1/t) - k\log(1/t)^{k-1}]\chi_{\f12}(t)\phi,\phi)
\leq K\|(1-\chi_{\f12}(t))F_k(t)\phi\|^2.
$$
Now, instead of multiplying by $1/k!$, we multiply by 
$A^k/(k+1)!$, $A$ to be determined.  
This gives 
$$
\multline
([(n-p-q - Cr^2)\log(1/t^A)^{k}/(k+1)! - Ak\log(1/t^A)^{k-1}/(k+1)!]\chi_{\f12}(t)\phi,\phi)\\ 
\leq K\|(1-\chi_{\f12}(t))F_k(t)\phi\|^2A^k/k!. 
\endmultline
$$
Summing over $k\geq -1$ gives 
$$
\multline (n-p-q-A)\|\chi_{\f12}(t)\phi/t^{A/2}\log^{1/2}(1/t)\|^2 + 
A\|\chi_{\f12}(t)\phi/t^{A/2}\log(1/t)\|^2\\
\leq C \|r\chi_{\f12}(t)\phi/t^{A/2}\log^{1/2}(1/t)\|^2 +
2^AK\|(1-\chi_{\f12}(t))\phi\|^2 .  
\endmultline 
$$
This inequality implies  
$\phi/r^{A}\log(1/r^2)\in L_2$ if $n-p-q\geq A.$ 
\bigskip  
Now let $\Delta_w$ denote the {\it weak Laplacian} with respect to smooth compactly supported forms: $\Delta_w\phi=\psi$ if and only if
$$
(\phi,\Delta\tau)=(\psi,\tau)
$$
for all smooth compactly supported $\tau$; if $\Delta_w\phi=0$, then $\phi$ is called {\it weakly harmonic}. Here is a useful identity [Ag]. 

\proclaim{\rm 2.34}{\bf Lemma} If $g:V\to\Bbb R$ is a smoooth function supported away from 0 and $\p,\,\,\Delta_w\phi\in L_2$, then
$$
\|Dg\p\|^2=\|[D,g]\p\|^2+(\Delta_w\phi,g^2\p)
$$
\endproclaim
{\bf Proof} Since $g$ is supported away from 0, we use integration by parts and the identity $dg\p=gD\p+[D,g]\p$ to get
$$
(\Delta_w\phi,g^2\p)=(D^2\p,g^2\p)=(D\p,Dg^2\p)=(gD\p,gD\p)+2(gD\p,[D,g]\p)
$$
and 
$$
\|Dg\p\|^2=\|[D,g]\p\|^2+2([D,g]\p,gD\p)+\|gD\p\|^2
$$
which together give the result.

\proclaim{\rm 2.35}{\bf Lemma} Let $\phi\in L^{p,q}_2(V)$ be weakly harmonic where $p+q\leq n-2$. Then $\log^k(1/t)\phi\in L_2(U_{\f12})$ for all real $k$.
\endproclaim
{\bf Proof} Suppose we can show that for each $k$ the following expression is bounded:
$$
\multline
([2k\log^{k-1}(1/t)|du|^2 -  k^2\log^{k-2}(1/t)|du|^2]\chi_{\f12}(t)\phi,\phi) - \\ (|du|^2k\log^{k-1}(1/t)\chi_{\f12}(t)\phi_0,\phi_0)+ 
(|du|^2k\log^{k-1}(1/t)\chi_{\f12}(t)\phi_3,\phi_3).
\endmultline\tag2.36
$$
Then by induction on $k$ we conclude that $\chi_{\f12}(t)\log^k(1/t)\phi$ is $L_2$ for all $k$ as claimed.

We begin the proof of boundedness by making special choices of $f'$ and $F$ in Corollary 2.15. For any $T$ such that $0\leq T\leq \f12$, we set 
$$
f'(t):=f'_T(t)=\left\{
\aligned
&\log^k(1/t),\, \f12\geq t\geq T\\
&\log^k(1/T),\, t\leq T
\endaligned
\right.
$$
and extend in a bounded fashion for $t\geq \f12$, with $F_k(t)/\log^{\f k2}(2)$ uniformly bounded as $k\to\infty$ as we did in the proof of (2.27).
 Then set
$$
F^2(t):=F_T^2(t)=r^2f'_T.
$$
Now using the cut-off $\mu_j=\mu_j(t)$ above and omitting the $t$ variable as usual to simplify notation, we have from (2.15)
$$
\multline
(n-p-q)(f_T'\mu_j\phi,\mu_j\phi) + 
(|du|^2 r^2f_T''\mu_j\phi_0,\mu_j\phi_0) - 
(|du|^2 r^2f_T''\mu_j\phi_3,\mu_j\phi_3)\\
 \leq(\|F_TD\mu_j\phi\|+\|F_TD'\mu_j\phi\|)\cdot \|(2^{-1/2}|du|(rf_T'/F_T)\mu_j\phi\|
\endmultline\tag2.37
$$
Now take $g=F_T\mu_j$ in (2.34) above. Since $\Delta_w\phi=0$ and $[D,F_T\mu_j]=F_T[D,\mu_j]+[D,F_T]\mu_j$, we have 
$$
\multline
\|F_TD\mu_j\p\|\leq\|DF_T\mu_j\p\|+\|[D,F_T]\mu_j\p\|=\|[D,F_T\mu_j]\p\|+\|[D,F_T]\mu_j\p\|\\ \leq \|F_T[D,\mu_j]\p\|+2\|[D,F_T]\mu_j\p\|=\|F_T[D,\mu_j]\p\|+2\|[D,F_T]\mu_j\p\|
\endmultline
$$
Plugging this and its $D'$ analogue into (2.37) and using Cauchy-Schwarz, $[D,F_T]=2rF'_T[e(\db u)+e^*(\db u)]$ and $|\db u|=2^{-1/2}|du|$, we get
$$
\multline
(n-p-q)(f_T'\mu_j\phi,\mu_j\phi) + (|du|^2r^2f_T''\mu_j\phi_0,\mu_j\phi_0) -(|du|^2 r^2f_T''\mu_j\phi_3,\mu_j\phi_3)\\
 \leq 4\|2^{-1/2}|du|2rF'_T\mu_j\p\|\cdot \|2^{-1/2}|du|(rf_T'(t)/F_T)\mu_j\phi\|+\|F_T[D,\mu_j]\p\|\cdot\|(2^{-1/2}|du|rf_T'(t)/F_T)\mu_j\phi\|\\ \leq \||du|2rF'_T\mu_j\p\|^2+\|(|du|rf_T'(t)/F_T)\mu_j\phi\|^2+\|F_T[D,\mu_j]\p\|\cdot\|(2^{-1/2}|du|rf_T'(t)/F_T)\mu_j\phi\|\\
\endmultline
$$
Now the part of the integral over $V-U_\f12$ in the first term on the left side of this inequality is positive, so it can be discarded. So with obvious notation we have, for some positive $C$ and $K$
$$
\multline
(n-p-q)\bigl(\chi_{\f12}\left\{
\aligned
&\log^k(1/t)\\
&\log^k(1/T)
\endaligned
\right\}\mu_j\phi,\mu_j\chi_\f12\phi\bigr) - \bigl( |du|^2r^2\chi_{\f12}\left\{
\aligned
&k\log^{k-1}(1/t)\\
&0
\endaligned
\right\}\mu_j\phi_0,\mu_j\phi_0\bigr)\\ + \bigl( |du|^2r^2\chi_{\f12}\left\{
\aligned
&k\log^{k-1}(1/t)\\
&0
\endaligned
\right\}\mu_j\phi_3,\mu_j\phi_3\bigr)\\
\leq 
\bigl(\chi_{\f12}\left\{
\aligned
&(\log^{\f k2}(1/t)-k\log^{{\f k2}-1}(1/t))^2\\
&\log^k(1/T)
\endaligned
\right\}\mu_j\phi,\mu_j\phi\bigr)+\bigl(\chi_{\f12}\left\{
\aligned
&\log^k(1/t)\\
&\log^k(1/T)
\endaligned
\right\}\mu_j\phi,\mu_j\phi\bigr)\\+\|F_T[D,\mu_j]\p\|\cdot\|2^{-1/2}|du|(rf_T'/F_T)\mu_j\phi\|+C\|(r\chi_{\f12}\log^{\frac{k}{2}}(1/t)\mu_j\phi\|^2 +K\|(1-\chi_{\f12}(t))\log^{\frac{k}{2}}(2)\p\|,
\endmultline
$$
where the term containing $C$ arises from the relation $t=r^2+O(r^2)$. 
The terms of the integrals on the left and right which involve $\log^k(1/T)$ and $\log^k(1/t)$, except for the term containing $C$, are nonnegative when moved to the left side ($n-p-q\geq 2$), so may be discarded. We rewrite the inequality as

$$
\multline
([ 2k\log^{k-1}(1/t)|du|^2 -  k^2\log^{k-2}(1/t)|du|^2]\chi_{[T,{\f12}]}\mu_j\phi,\mu_j\phi) \\- (|du|^2k\log^{k-1}(1/t)\chi_{[T,{\f12}]}\mu_j\phi_0,\mu_j\phi_0) + 
(|du|^2k\log^{k-1}(1/t)\chi_{[T,{\f12}]}\mu_j\phi_3,\mu_j\phi_3)\\
\leq 
\|F_T[D,\mu_j]\p\|\cdot\|(2^{-1/2}|du|(rf_T'/F_T)\mu_j\phi\|+
C\|(r\chi_{\f12}\log^{\frac{k}{2}}(1/t)\mu_j\phi\|^2+K\|(1-\chi_{\f12}(t))\log^{\frac{k}{2}}(2)\p\|
\endmultline
$$

Obviously the second term on the right is bounded uniformly in $j$. We claim that the first term on the right vanishes as $j\to\infty$. Once this latter point is verified, we will take the limit of both sides of this last inequality as $j\to\infty$ and then the limit as $T\rightarrow\infty$ to get 
the estimate 
$$
\multline
([2k\log^{k-1}(1/t)|du|^2 -  k^2\log^{k-2}(1/t)|du|^2]\chi_{\f12}(t)\phi,\phi) - \\ (|du|^2k\log^{k-1}(1/t)\chi_{\f12}(t)\phi_0,\phi_0)+ 
(|du|^2k\log^{k-1}(1/t)\chi_{\f12}(t)\phi_3,\phi_3)\\
\leq C\|(r\chi_{\f12}\log^{\frac{k}{2}}(1/t)\phi\|^2+K\|(1-\chi_{\f12}(t))\log^{\frac{k}{2}}(2)\p\|
\endmultline\tag2.38
$$
which proves the boundedness of (2.36).

It remains to show that for fixed $T>0$ the term
$$
\|F_T[D,\mu_j]\p\|\cdot\|2^{-1/2}|du|(rf_T'(t)/F_T)\mu_j\phi\|
$$
 tends to zero as $j$ tends to infinity. First note that for $t<T$, $rf_T'(t)/F_T(t)=\log^{\f k2}(1/T)$; hence 
$\|(rf_T'(t)/F_T)\mu_j\phi\|$ is uniformly bounded as $j\to\infty$.
 Again for $t<T$, the pointwise norm squared
$$
|F_T[D,\mu_j]\p|^2\leq \frac{2\chi_{I(j)}(t)r^2\log^k(1/T)}{t\log^2(t)}|\p|^2
$$ 
which is bounded and tends pointwise to zero as $j\to\infty$. This completes
 the proof.

%{\bf Proof:} Since $\db f(r^2)/2=rf'(r^2)\db r=r\log^{k/2}(1/r^2)$ and %$\db(\mu_j\phi=\db\mu_j\wedge\phi+\mu_j\db\phi$, we have
%$$
%(\db \phi_j,e(\db f(r^2)/2)\phi_j)=(\db\mu_j\wedge r\log^s(1/r^2)\phi,e(\db %r)\phi_j)+
%(\mu_j\db\phi,e(\db f(r^2)/2)\phi_j)
%$$
%Now take $\xi=r\log^s(1/r^2)\phi$ in the Lemma to see that the first term in ? %tends to 0; and use the first property of $\mu_j$ to see that the second tends %to $(\db\phi,e(\db f(r^2)/2)\phi)$. This takes care of the first term in ?; the %econd is handled similarly. The third term is zero-th order in $\phi$ and so %needs no comment. As for the fourth, note first that by our choices of $f$ and %$F$, $([e^*(\partial f(r^2)) - (e(\partial f(r^2))]\phi/2F)=[e^*(\partial r) - %(e(\partial r)]\phi/)$, whiich is $L_2$; and since %$\phi=F(r^2)\phi/r\log(1/r^2)\in L^{p,q}$, we handle this term as we did the %first two. The last two terms are zero-th order in $\phi$; and $\da %F(r^2)=2r[s\log^{s-1}(1/r^2)+\log^s(1/r^2)]\da r$ so that $e^*(\partial F)\phi$ %and $e(\partial F)\phi$ are $L^2$, so we are done.

 Let $\Delta_{d_D}$ denote the ``strong" Dirichlet Laplacian in the usual sense of functional analysis: $\Delta_{d_D}\phi=\psi$  if and only if $\phi\in\dom d_D\cap\dom\delta_N$, $d_D\phi\in\dom\delta_N$, $\delta_N\phi\in\dom d_D$ and 
$$
\delta_N d_D\phi+d_D\delta_N\phi=\psi.
$$
(Note that $\delta_N$ is the Hilbert space adjoint of $d_D$.) The kernel of $\Delta_{d_D}$ in degree $k$ is denoted $\C H^k_{d_D}$ and its elements are called {\it strongly} $d_D$-{\it harmonic}. Analogous definitions can be made for $\Delta_{\da_D}$ and $\Delta_{\db_D}$ and also for Neumann boundary conditions. It follows from Stokes' theorem in the usual way that a form $\p$ is strongly $d_D$-harmonic if and only if $d_D\p=0=\delta_N\p$; again there is an analogous statement for the other Laplacians. 
Evidently, if $\p$ is strongly harmonic in any sense, it is weakly harmonic. Outside the middle three degrees the converse holds:

\proclaim{\rm 2.39}{\bf Theorem} Suppose $\phi$ is a weakly harmonic $(p,q)$-form such that  $n-p-q\geq 2$. Then  $\phi$ is strongly harmonic in any of the above senses, 
and hence 
$$
\phi/r^{n-p-q}\log(1/r^2)\in L_2.
$$
\endproclaim

{\bf Proof:} Divide both sides of the estimate (2.38) by $(k+2)!$ and expand to obtain
$$
\multline
([((n-p-q)-2|du|^2)\log^k(1/t)/(k+2)! + 2\log^{k-1}(1/t)/(k+1)!|du|^2 \\
- 4\log^{k-1}(1/t)/(k+2)!|du|^2
 -  \log^{k-2}(1/t)/k!|du|^2  \\
+ 3\log^{k-2}(1/t)/(k+1)!|du|^2 - 8\log^{k-2}(1/t)/(k+2)!|du|^2]
\chi_{\f12}(t)\phi,\phi)\\ - 
(|du|^2[(\log^{k-1}(1/t)/(k+1)! -2\log^{k-1}(1/t)/(k+2)!]\chi_{\f12}(t)\phi_0,\phi_0)
\\ + 
(|du|^2[\log^{k-1}(1/t)/(k+1)!  - 2\log^{k-1}(1/t)/(k+2)!]\chi_{\f12}(t)\phi_3,\phi_3)\\
\leq  C\|r\chi_{\f12}\log^{\frac{k}{2}}(1/t)\phi\|^2/(k+2)!+K\|(1-\chi_{\f12}(t))\log^{\frac{k}{2}}(2)\p\|/(k+2)!
\endmultline
$$
Sum over $k\geq -2$  to obtain 

$$
\multline
((n-p-q-|du|^2)/t\log^2(1/t) - |du|^2/t\log^3(1/t) - 8|du|^2/t\log^4(1/t)]
\phi,\chi_\f12\phi)\\
 + 
(|du|^2[(-1/t\log^2(1/t) + 2/t\log^3(1/t)]\phi_0,\chi_\f12\phi_0)
\\ + 
(|du|^2[1/t\log^2(1/t) - 2/t\log^3(1/t)]\phi_3,\chi_\f12\phi_3)\\
\leq C\|rt^{-\f12}\log^{-1}(1/t)\chi_{\f12}\phi\|^2+2\log^{-2}(2)K\|(1-\chi_{\f12}(t))\p\|
+\|2^{-\f12}\log^{-2}(1/t)\chi_{\f12}\phi\|^2
\endmultline
$$

>From this inequality, we deduce that $\phi_i/r\log(1/r^2)\in L_2$ 
for $i>0$, and $\phi_0/r\log^{3/2}(1/r^2)\in L_2$.  If 
$n-p-q>2$, we get $\phi/r\log(1/r^2)\in L_2$; to get this for $n-p-q=2$, we must argue more carefully.  
To begin, we choose for $t\leq\f12$
$$
F(t)=\log^{\f k2}(1/t),
$$
extend in a bounded fashion for larger $t$ and get from (2.34) the inequality
$$
\|\mu_jFD\p\|\leq 2\|[D,\mu_jF]\p\|.
$$
If we assume that $t^{-\f12}\log^{\frac{k}{2}-1}(1/t)\p$ is $L_2,$ then 
$[D,\mu_jF]\p$ is $L_2$ and, as $j\rightarrow\infty$, converges 
 in $L_2$ to 
$[D,F]\p.$ Hence taking the limit, we obtain the following lemma.

\proclaim{\rm 2.40}{\bf Lemma} If $t^{-\f12}\log^{\frac{k-1}{2}}(1/t)\p$ is $L_2$ on $U_\f12$, then so is $\log^{\f k2}(1/t)D\p$.
\endproclaim

Now return to our basic estimate in Proposition (2.15) and take for $t\leq \f12$
$$
f'(t)=t^{-1}\log^k(1/t)\qquad\text{and}\qquad F(t)=\log^{\frac{k}{2}}(1/t)
$$
so that
$$
r^2f''(t)= -t^{-1}\log^k(1/t)-t^{-1}k\log^{k-1}(1/t)+O(r^2)\log^k(1/t)
$$
there; and extend both $f'$ and $F$ in a bounded fashion as usual for $t\geq\f12$. Using $n-p-q\geq 2$ and discarding some positive terms from the left side, we get, for some positive constants $C$ and $K$ 
$$
\multline
\|t^{-\f12}\log^{\f k2}(1/t)\mu_j\p\chi_\f12\|^2 - \|kt^{-\f12}\log^\f {k-1}{2}(1/t)\mu_j\phi_0\chi_\f12\|^2\\
 \leq(\|\chi_\f12\log^{\f k2}(1/t)D\mu_j\phi\|+\|\chi_\f12\log^{\f k2}(1/t)D'\mu_j\phi\|)\cdot \|\chi_\f12(\log^{\f k2}(1/t)/r)\mu_j\phi\|+\\
+C\|r\log^{\f k2}(1/t)\mu_j\p\chi_\f12\|+K(\|(1-\chi_{\f12}(t))\p\|+\|(1-\chi_{\f12}(t))D\p\|+\|(1-\chi_{\f12}(t))D'\p\|)
\endmultline
$$
Now if we assume that $\log^{\f k2}(1/t)D\phi$ is $L_2$, then we may  use (2.34) as we did above to bring the cut-off $\mu_j$ past $D$ and $D'$ on the right side and get
$$ 
\multline
\|t^{-\f12}\log^{\f k2}(1/t)\mu_j\p\chi_\f12\|^2 - \|kt^{-\f12}\log^\f {k-1}{2}(1/t)\mu_j\phi_0\chi_\f12\|^2\\
 \leq(\|\chi_\f12\log^{\f k2}(1/t)\mu_jD\phi\|+\|\chi_\f12\log^{\f k2}(1/t)\mu_jD'\phi\|+\ell)\cdot \|\chi_\f12(\log^{\f k2}(1/t)/r)\mu_j\phi\|\\
+C\|r\log^{\f k2}(1/t)\mu_j\p\chi_\f12\|+K(\|(1-\chi_{\f12}(t))\p\|+\|(1-\chi_{\f12}(t))D\p\|+\|(1-\chi_{\f12}(t))D'\p\|)
\endmultline
$$
for some constant $\ell$. Let us now assume in addition that 
$\chi_{\f12}t^{-\f12}\log^\f {k-1}{2}(1/t)\phi$ is $L_2$ and divide both sides by $\|(\log^{\f k2}(1/t)/r)\mu_j\p\|$. Remembering that $t^{\f12}\sim r$, we see that the right side is then bounded as $j\to\infty$, so the left is as well. Hence we obtain
\proclaim{\rm 2.41}{\bf Lemma} If $\log^{\f k2}(1/t)D\phi$ and $t^{-\f12}\log^\f {k-1}{2}(1/t)\phi$ are $L_2$ on $U_\f12$, so is $t^{-\f12}\log^{\f k2}(1/t)\p$.
\endproclaim
Now induction on the hypothesis of this last lemma, beginning with $k=-2$, together with (2.40) and (2.41), proves that $\chi_\f12 t^{-\f12}\log^k(1/t)\phi$ and $\chi_{\f12}\log^k(1/t)D\phi$ are  $L_2$ for all $k$. In particular, $D\phi$ is $L_2$, so $D'\phi$ is as well by (2.20).

To prove that $\p$ is strongly $d_D$-harmonic, observe first that $\p\in\dom d_D\cap\dom \delta_D$ by (2.18).  
Hence 
$$\|d\p\|^2 + \|\delta\p\|^2 = 
\lim_{j\rightarrow\infty} (d\p,d\mu_j\p) + (\delta\p,\delta\mu_j\p) $$
$$= (\Delta\p,\mu_j\p) = 0.$$

That $\p$ is strongly harmonic in the other senses is proved similarly; this and Proposition (2.33) complete the proof of (2.39).
\bigskip

\proclaim{\rm 2.42}{\bf Corollary} For $|n-k|\geq 2$ or $|n-p-q|\geq 2$, 
$$
\Cal H^k_w(V)=\Cal H^k_B(V)\quad\text{and}\quad\Cal H^{p,q}_w(V)=\Cal H^{p,q}_B(V)
$$
for all boundary conditions $B$.
\endproclaim
{\bf Proof} As remarked before the statement of the theorem, it is easy to see the inclusions $\supseteq$; and for $(p,q)$-forms with $n-p-q>0$, the opposite inclusion is part of the Theorem. If $n-p-q<-2$, the result follows from this and \cite{PS, 1.3}. The $(p,q)$-components of a weakly harmonic $k$-form are weakly $\db$- and $\da$-harmonic, so we are done.

The preceding estimates are not quite strong enough to carry over to 
$p+q = n-1,$ where one can obtain estimates in the complete case.  We need to 
introduce a variational argument to handle $p+q = n-1$.

\proclaim{\rm 2.43}{\bf Theorem} Let $\p\in \ker d_D$ and suppose $\deg\p = n-1$. Then the $d_D$-harmonic representative $h$ of $\p$ 
satisfies
$$
\delta_Dh=d^c_Dh=\delta^c_Dh=0,
$$
and $h/(r\log(1/r^2))\in L_2$.
\endproclaim
\noindent

{\bf Proof.} There is a sequence $\p_j\in A^{n-1}_c(V)$ such that $\p_j\to\p$ in $L^{{n-1}}(V)$ and $d\p_j\to 0$ in $L^{n}(V)$. Let 
$$
Q:L^{{n-1}}(V)\times L^{{n-1}}(V)\to \Bbb C
$$
be the unbounded, densely defined Hermitian form 
$$
Q(\alpha,\beta)=(\delta\alpha,\delta\beta)+(d^c\alpha,d^c\beta)+(\delta^c\alpha,\delta^c\beta).
$$
Let $q(\alpha)=Q(\alpha,\alpha)$ denote the corresponding quadratic form. Then $q$ is nonnegative and $\p_j\in\dom q$ since it is compactly supported; so we can define 
$$
m_j=\inf\{q(\p_j+d_D\beta)|\beta\in\dom d_D\}
$$
We now show by a standard argument that this infimum is realized.

\proclaim{\rm 2.44}{\bf Lemma} There is $\beta_j\in\dom d_D$ such that 
$q(\p_j+d_D\beta_j)=m_j$
\endproclaim
\noindent
{\bf Proof.} For the proof, let us drop the subscript $j$'s on $\p_j$and $m_j$. 
Let $S$ denote the Hilbert space closure of $dA_c^{n-2}(V)$ 
with respect to the norm induced by $Q$.  The finite dimensionality 
of the $L_2$ cohomology (and the closed graph theorem) imply this 
norm dominates a multiple of the $L_2$ norm.  Moreover, by the 
ellipticity of $d+ \delta$, the norm is equivalent on compact subsets 
to the Sobolev norm of forms with one $L_2$ derivative. 

Choose a minimizing sequence $db_i\in S$ with
 $q(\p+d_Db_i)\downarrow m$. Because this is bounded in $S$, we 
may extract a weakly convergent subsequence, which we also label 
$db_i$, converging to some limit $Z$, which is clearly in the range of 
$d_D$, 
$$
Z = d_D\beta_j,
$$
for some $\beta_j$.  Then 
$$
m=\lim_{i\rightarrow\infty}q(\p + db_i) = 
q(\p) + 2\text{Re}\lim_{i\rightarrow\infty}Q(\p,db_i) + \lim_{i\rightarrow\infty}q(db_i)$$
$$ 
= q(\p) + 2\text{Re}\, Q(\p,Z) + \lim_{i\rightarrow\infty}q(db_i).
$$
Recalling that for weak limits, 
$\lim_{i\rightarrow\infty}q(\p+db_i) \geq q(\p+Z)$, gives 
$q(\p+Z)\leq m$. By hypothesis, we also have $q(\p+Z) \geq m$ and hence 
$q(\p+Z) = m$. The infimum is achieved.

Let us now denote a minimizing form constructed above
$$
\psi_j:=\p_j+d_D\beta_j.
$$
By construction, 
$$d\psi_j=d\p_j \text{  and  } \psi_j\in\dom d_N^c\cap\dom \delta_N^c.$$

\proclaim{\rm 2.45}{\bf Lemma} $\delta\psi_j=0$, and $\psi_j\in\dom d_D^c\cap\dom \delta_D^c$.
\endproclaim
{\bf Proof.} We apply the usual variational argument: by the minimality of $q(\psi_j)$, we have for smooth compactly supported $w$,
$$0 = \frac{d}{dt}q(\psi_j + tdw)_{|t=0} = $$
$$(\delta\psi_j,\delta dw) + (d^c\psi_j,d^c dw) + (\delta^c\psi_j,\delta^c dw).$$
Integrating by parts and using the K\"ahler identities gives 
 $\Delta\delta\psi_j = 0$.  Thus $\delta\psi_j$ is weakly 
harmonic and of degree $n-2$. 
By (2.39) this implies that it is strongly harmonic and therefore perpendicular 
to the image of $\delta$. Hence
$\delta\psi_j = 0$.  Now the $L_2$ boundedness of $D\psi_j$ and $D'\psi_j$ allows us to use 
(2.21) to conclude $\psi_j/r\log(1/r^2)$ is $L_2$. This then implies  
$\psi_j\in\dom d_D^c\cap\dom \delta_D^c$.  
\medskip

Returning to the proof of the theorem, we have constructed a sequence 
$\psi_j=\p+d_D\beta_j$ with $d\psi_j\rightarrow 0$ and  $\delta\psi_j = 0$ and 
$\|D'\psi_j\| \leq \|d\psi_j\|$. Hence as $j\rightarrow\infty,$ 
$\psi_j$ converges to the harmonic representative $h$ of $\phi.$  
Moreover, according to (2.25), $\|\psi_j/r\log(1/r^2)\|$ is bounded in terms of 
$\|D\psi_j\|$, $\|D'\psi_j\|$ and 
 $\|\psi_j\|$. Hence, the harmonic limit $h$ also satisfies 
$h/r\log(1/r^2)\in L_2$. From this we may immediately deduce the claims of 
the theorem by applying the K\"ahler identities (cf. (2.19) and (2.20)).

\proclaim{\rm 2.46}{\bf Corollary} Let $\p\in L^{n-1}(V)$ and suppose $d_D\p=\xi\in L^{n}(V)$. Then there is $\psi\in L^{n-1}(V)$ such that $\psi/(r\log(1/r^2))\in L^{n-1}(U_\f12)$, $d_D\psi=\xi$ and $\delta \psi = 0$.
\endproclaim
\noindent

{\bf Proof} Let $\phi_j$ be a sequence of compactly supported forms so that 
 $\phi_j\to\phi$ and $d\phi_j\to\xi$; construct $\psi_j=\p_j+d_D\beta_j$ as in the preceding theorem.  
We can then write for some $(d_D)$-harmonic $h_j$
$$\psi_j = h_j + \delta_N\alpha_j.$$ 
Since $d_D\psi=d_D\p\to\xi$ and $d_D$ has closed range, $\delta_N\alpha_j$ converges; $h_j$ also converges (it is the harmonic component of $\p_j$ too), so $\psi_j$ does. By the preceding theorem, $h_j/r\log(1/r^2)\in L_2$ and by the 
proof of the preceding theorem $\psi_j/r\log(1/r^2)\in L_2.$ 
Hence $\delta_N\alpha_j/r\log(1/r^2)\in L_2$.  Taking the limit as $j\rightarrow\infty$, 
we obtain the desired result.

Here is a variation on the previous result in which we must keep track of $(p,q)$-type in reaching a slightly different conclusion. It will be used in \S4.

\proclaim{\rm 2.47}{\bf Theorem} Let $\p\in L^{p,q}$ and suppose $\db_D\p=\xi\in L^{p,q+1}$. Then if $p+q<n$, there is $\psi\in L^{p,q}$ such that $\psi/(r\log(1/r^2))\in L^{p,q}$, $\db_D\psi=\xi$ and $\da_D\psi\in L^{p+1,q}$.
\endproclaim

{\bf Proof} The proof is quite similar to that of the previous theorem and will only be sketched. By hypothesis there is a sequence $\p_j\in A^{p,q}_c(V-\{0\})$ such that $\p_j\to\p$ and $\db\p_j\to\xi$. Using the Hermitian form  
$$
Q:L^{p,q}(V)\times L^{p,q}(V)\to \Bbb C
$$
where 
$$
Q(\alpha,\beta)=(\vartheta\alpha,\vartheta\beta)+(\da\alpha,\da\beta)+(\bar\vartheta\alpha,\bar\vartheta\beta)
$$
as we did above, we get a sequence $\psi_j=\p_j+\db_D\beta_j$ such that

\proclaim{\bf Lemma} $\db\psi_j=\db\p_j$, $\vartheta\psi_j=0$, and $\psi_j\in\dom\db_D\cap\dom\bar\vartheta_D$.
\endproclaim

Since $\vartheta\psi_j=0$, $\db_D\beta_j$ is the $\db_D$-exact component of $\p_j$ in its Hodge decomposition, so $\db_D\beta_j$ converges because $\p_j$ does. Hence $\psi_j$ converges, say $\psi_j\to\psi\in L^{p,q}(V)$. 
Again using (2.25), we get $\psi_j/r\log(1/r^2)\in L_2$ and from this 
$$
\|\da\psi_j\|^2+\|\bar\vartheta\psi_j\|^2=\|\db\psi_j\|^2
$$
Now $\db\psi_j$ is Cauchy, so this shows $\da\psi_j$ is as well. Since $\psi_j/r\log(1/r^2)\in L_2$ (uniformly) implies $\psi/r\log(1/r^2)\in L_2$, we have $\psi\in\dom \da_D$
and are done.
\bigskip

Let $\Cal H^k_{D/N}(V)$ denote the space of harmonic forms on $V$ in degree $k$ with respect to the operator
$$
\hat d^k:=\left\{
\aligned d_D\,\,&k<n,\\
         d_N\,\,&k\geq n \endaligned
\right.
$$ 
so that
$$
 \Cal H^k_{D/N}(V)
=\left\{
\aligned \ker d_D\cap\ker\delta_N,\,\,&k<n,\\
         \ker d_N\cap\ker\delta_N,\,\,&k= n \\
         \ker d_N\cap\ker\delta_D,\,\,&k>n\endaligned\right.\tag2.48
$$ 
We now verify that the operators $L$ and $\Lambda$ act on $\Cal H^\ast_{D/N}(V):=\oplus_{k\geq 0} \Cal H^k_{D/N}(V)$, satisfying the standard K\"ahler identities. In an unfortunate convergence of notation, we 
 let $L^{\ast}(V):=\oplus_{k\geq 0}L^k(V)$, let $\Pi_k:L^{\ast}(V)\to L^{\ast}(V)$ be the projection to $L^k(V)$, and let $H=\oplus_{k\geq0}(n-k)\Pi_k:L^{\ast}(V)\to L^{\ast}(V)$. Here 
$L^*$ is not to be confused with the adjoint of $L$.

\proclaim{\rm 2.49}{\bf Theorem} The operators $L$, $\Lambda$, and $H$ preserve $\Cal H^\ast_{D/N}(V)$ and satisfy
$$
[\Lambda,L]=H,\qquad [H,L]=-2L,\qquad [H,\Lambda]=2\Lambda
$$
In particular, $\Cal H^\ast_{D/N}(V)$ is an $sl_2(\Bbb C)$-module and so the Lefschetz decomposition theorem holds for $L_2$-cohomology.
\endproclaim
{\bf Proof} It is sufficient to verify that the three operators preserve $\Cal H^\ast_{D/N}(V)$. Since this is clear for $H$ and since $\Lambda$ is dual to $L$, it is enough to show (omitting $V$ from the notation now)
$$
L\Cal H^{k}_{D/N}\subseteq\Cal H^{k+2}_{D/N}
$$
For all $k$, $L\Cal H^{k}_w\subseteq\Cal H^{k+2}_w$; and 
$\Cal H^{k+2}_{D/N}\subseteq\Cal H^{k+2}_w$ with equality except possibly for $k=n-3,n-2,n-1$ by (2.42). In these cases, for $\phi\in\Cal H^{k}_{D/N}$, we have $\p/(r\log(1/r^2))\in L^k(U_{\f12})$ by (2.39) and (2.43), so that boundary conditions on $\phi$ are irrelevant and $\p\in\ker d\cap\ker\delta\cap\ker d^c\cap\ker\delta^c$. Since $dL\p=Ld\p=0$ and $\delta L\p=L\delta\p-4\pi d^c\p=0$ and  $L\p/(r\log(1/r^2))\in L^{k+2}(V)$, we get $d_D L\p=0$ and $\delta_D L\p=0$, so the proof is complete.
\bigskip

We are now ready to put a Hodge structure on $L_2$-cohomology. 
Let $\Cal H^{p,q}_{D/N}(V)$ denote the space of harmonic forms 
in degree $(p,q)$ on $V$ with respect to the operator
$$
\Hat{\Bar\partial}^{p,q}:=\left\{
\aligned \bar\partial_D^{p,q}\,,\,\,&p+q<n,\\
         \bar\partial_N^{p,q}\,,\,\,&p+q\geq n; \endaligned
\right.
$$
so that
$$
\Cal H^{p,q}_{D/N}(V)
=\left\{
\aligned \ker \bar\partial_D\cap\ker\vartheta_N,\,\,&p+q<n,\\
         \ker \bar\partial_N\cap\ker\vartheta_N,\,\,&p+q=n \\
         \ker \bar\partial_N\cap\ker\vartheta_D,\,\,&p+q\geq n\endaligned\right.
$$
\proclaim{\rm 2.50}{\bf Theorem} Let $V$ be a complex projective variety of dimension $n$ with at most isolated singularities. Then 
$$
L(\Cal H^{p,q}_{D/N}(V))\subseteq \Cal H^{p+1,q+1}_{D/N}(V)
$$
and for each $k=0\dots 2n$, we have the equality of subspaces of $L^k(V)$
$$
\Cal H^k_{D/N}(V)=\oplus_{p+q=n} \Cal H^{p,q}_{D/N}(V)
$$
where the summands on the right side are the $(p,q)$-components of the left side.
\endproclaim 
{\bf Proof} The first assertion is proved in the same way as the last theorem. The K\"ahler identities imply that
$$
\Cal H^k_w(V)=\oplus_{p+q=n} \Cal H^{p,q}_w(V)
$$
so for $|n-k|\geq 2$, the Theorem follows from Cor. (2.42). If $\p\in\Cal H^{n-1}_{D/N}$, then Theorem (2.43) and (2.19), (2.20) say we have 
$$
\|d\p\|^2+\|\delta\p\|^2=\|\db\p\|^2+\|\vartheta\p\|^2 
$$
which verifies the theorem in this case; the equality in case $k=n-1$ follows from this and duality (\cite{PS, 1.3}.  In the remaining case,
$$
\Cal H^n_{D/N}=\ker L\oplus\im L
$$
by the Lefschetz decomposition and $\phi\in\Cal H^n_{D/N}$ if and only if $d_N\phi=0=\delta_N\phi$. If $\p=L\psi$, where $\psi \in \Cal H^{n-2}_{D/N}$, then $\psi=\sum\psi^{p,q}$, where $\psi^{p,q}\in\Cal H^{p,q}_{D/N}(V)$, so we are done by the first assertion of the Theorem. In case $L\p=0$, $\Lambda\p=0$ as well, and so we may use the K\"ahler identities
$$
[L,\delta]=4\pi d^c\qquad[\Lambda,d]=-4\pi\delta^c 
$$
to conclude from $d\p=\delta\p=0$ that $d^c\p=\delta^c\p=0$. (Subscript $N$'s are intended on the operators here.) This implies that $\db\phi= \vartheta\phi=\da\phi=\bar\vartheta\phi=0$ so that if $\phi=\sum\phi^{p,q}$, then $\phi^{p,q}\in\Cal H^{p,q}_{D/N}$. This completes the proof.

Finally, we can conclude in the usual way (\cite{Hi, \S15.8}) that our Hodge structure (2.49), with its Lefschetz decomposition (2.50) is polarized, in the following (standard) sense, by the inner product $(\xi,\eta)$.
\medskip
\proclaim{\rm 2.51}{\bf Definition} Let $A$ be a subring of $\Bbb R$ such that
$A\otimes \Bbb Q$ is a field. A {\rm polarized} $A$-{\rm Hodge structure of weight} k is a Hodge structure $(P_A;P_{\Bbb C},F^{\cdot};i)$, together with a symmetric bilinear form $Q:P_A\times P_A\to A$, such that
$$
Q(F^p,F^{k-p+1})=0,\,\,\,\text{for all}\,\, p
$$
and
$$
i^{p-q}Q(v,\bar v)<0,\,\,\,\text{for all}\,\,v\in P^{p,q}
$$
where $P^{p,q}:=F^p\cap \overline{F^{k-p}}$. We say {\rm the Hodge structure} $(P_A;P_{\Bbb C},F^{\cdot};i)$ {\rm is polarized by} Q.
\endproclaim

\newpage
\heading \S3 Hsiang-Pati Coordinates and the Nash bundle\endheading

Let $U$ be a small neighborhood of an isolated singular point $v$ on
a complex algebraic surface $V$ and as usual let $g$ denote the K\"ahler metric on $U-v$ inherited from an imbedding of $(U,v)\subseteq (\Bbb C^N,0)$. Let 
$$
\pi:(\tilde U,E) \to (U,v)\tag3.1
$$ 
be a resolution of the singularity $v$ of $U$. Then 
$$
\gamma:=\pi^*g\tag3.2
$$
is a K\"ahler metric on $U-E$. Hsiang and Pati showed ([HP]) that when $\pi$ is a  sufficently fine resolution, then, up to quasi-isometry, $\gamma$ assumes near $E$ a normal form in appropriate coordinates.\newline
(3.3)  Specifically, they showed that $\pi:(\tilde U,E) \to (U,v)$ can be chosen so that $E=\cup E_i$ is a divisor with normal crossings and has the following properties:
\roster
\ite a For each point $e \in E$, there is a neighborhood $W$ of $e$ in $\tilde U$ and linear functions $k,l:{\Bbb C}^N \to {\Bbb C}$ such that
$$
\gamma \sim d{\phi}d{\bar \phi}+d{\psi}d{\bar \psi}
$$
on $W$, where ${\phi}=l\circ{\pi}$ and ${\psi}=k\circ{\pi}$. 

This means that the linear projection $(l,k)
:\Bbb C^N\to\Bbb C^2$ is such that $(l,k)\circ\pi|W$ pulls back the Euclidean metric on $\Bbb C^2$ to one on $W$ which is quasi-isometric to $\gamma$.

\ite b {\it (Local description of $\phi$)}  $\phi$ locally defines the scheme-theoretic inverse image $\pi^{-1}(\frak m_v)$ of $v$ (where $\frak m_v$ is the maximal ideal of $v$). This means that if $w_1,\,w_2,\,\dots,\,w_N$ are coordinates on $\Bbb C^N$, then the restriction of the ideal $(w_1\circ\pi,\cdots,\,w_N\circ\pi):=\pi^{-1}(\frak m_v)$ in $\C O_{\tilde U}$ to $\C O_W$ is principal and is generated by $\phi$. Hence, $\pi^{-1}(\frak m_v)$ may be identified with its divisor $Z=:\sum m_iE_i$ and there are coordinates $u, v$ on $W$ such that if $e\in E_i \cap E_j$, then $E_i=\{u=0\}$, $E_j=\{v=0\}$ and  $\phi/u^{m_i}v^{m_j}$ is non-vanishing holomorphic in $W$; and if $e\in E_i$ is a simple point of $E$, then $E_i=\{u=0\}$ and $\phi/u^{m_i}$ is non-vanishing holomorphic in $W$. 

\ite c{\it (Local description of $\psi$)} There are integers $n_i\ge m_i$ such that $n_im_j-n_jm_i\neq0$ if $E_i \cap E_j\neq\emptyset$ and $\psi$ is the sum of two holomorphic functions $\psi=f(\phi)+{\psi^\prime}$, where $f=\sum a_jz^{\epsilon_j}$ is a series where the $\epsilon_j$ are rationals $\geq 1$ and ${\psi^\prime}$ defines a divisor $N:=\sum n_iE_i$ in $W$; in fact, with the {\it same} coordinates $u, v$ as in b), $\psi^\prime/u^{n_i}v^{n_j}$ is non-vanishing holomorphic in $W$ if $e\in E_i \cap E_j$; and if $e\in E_i=\{u=0\}$ is a simple point of $E$, then $\psi^\prime/u^{n_i}v$ is non-vanishing holomorphic in $W$. Moreover, $(n_i,n_j)$ (resp. $n_i$) is minimal with this property: if for some linear function $h:\Bbb C^N\to\Bbb C$, $h\circ\pi:=\eta=g(\phi)+\eta'$ with $g$ a series in rational powers $\geq 1$ and $\eta^\prime/u^{p_i}v^{p_j}$  non-vanishing holomorphic in $W$ (resp., $\psi^\prime/u^{p_i}v$ is non-vanishing holomorphic in $W$), where $p_i\ge m_i$, $p_j\ge m_j$ and $p_im_j-p_jm_i\neq0$ (resp., $p_i\ge m_i$), then $p_i\ge n_i$ and $p_j\ge n_j$ (resp., $p_i\ge n_i$).

\ite d On the above neighborhood $W$ of $e\in E_i \cap E_j$ let $\zeta_1=u^{m_i}v^{m_j}$ and $\zeta_2=u^{n_i}v^{n_j}$; or if $e\in E_i$ is a simple point of $E$, let $\zeta_1=u^{m_i}$ and $\zeta_2=u^{n_i}v$. Then in $W$
$$
{\gamma}\sim d\zeta_1d{\bar \zeta_1}+d\zeta_2d{\bar \zeta_2}
$$
So we have 
$$d{\phi}d{\bar \phi}+d{\psi}d{\bar \psi} \sim {\gamma} \sim d\zeta_1d{\bar \zeta_1}+d\zeta_2d{\bar \zeta_2}
$$
\endroster

(3.4) {\bf Remarks} \roster
\ite a Property (3.3d) is an easy consequence of the others (see \cite{HP, p. 401}). In an Appendix to this Chapter we will show, using properties of the Nash blow-up, the existence of a linear projection $(l,k)
:\Bbb C^N\to\Bbb C^2$ satifying the first three properties. 
\ite b It follows from (3.3b) that $|\phi|\sim r\circ\pi$ on $W$, where  $r^2=|w_1|^2+\dots+|w_N|^2$
\ite c It is not in general possible to remove $g(\phi)$ from the expression for $h\circ\pi$ in (3.3c).
\endroster
\medskip

\indent One final point which will prove useful later is that it is possible to choose a linear function $h$ which can be taken to be $l$ in (3.3a) outside a finite set of points of $E$; while near each point of this finite set, it can be taken to be $k$. The proof will also show that $l$ in (3.3a) and (3.3b) is generic among all linear functions $\Bbb  C^N\to \Bbb C$. Before stating the result, we give an example. 

\noindent (3.5) {\bf Example:} Let $V\subset\Bbb C^3(x, y, z)$ be the cone $\{y^2=xz\}$. Then blowing up $V$ at its singular point $(0,0,0)$ produces a resolution $\pi:\tilde V\to V$, where $\tilde V$ is the total space of the line bundle of degree -2 over $\Bbb P^1$, the exceptional divisor $E$, and $\pi$ collapses the zero section to $(0,0,0)$. Let $U$ be the intersection of a small ball about $(0,0,0)\in \Bbb C^3$ with $V$. Then $\tilde U:=\pi^{-1}(U)$ is covered by two open sets $\tilde U_1\subseteq\Bbb C^3(u,v)$ and $\tilde U_2\subseteq\Bbb C^3(u',v')$ which contain the $u$-axis and $u'$-axis respectively, and which are glued by 
$$
u'=uv^2
$$
$$
v'=v^{-1}
$$
Then the linear function $h$ on $U$ is $h(x,y,z)=y$, which is $uv=u'v'$ on $\tilde U$. Hence the proper transform $R\subset \tilde U$ of $h=0$ has two components transversely intersecting $E=\Bbb P^1$ at 0 and $\infty$. Notice that  if $l:=y+\epsilon x$ is a small perturbation of $y=h$, then the pair $\{l,h\}$ satifies the conditions of (3.3) in neighborhoods in $\tilde U$ of $0\in \Bbb P^1$ and $\infty\in \Bbb P^1$; and that in neighborhoods of all other points on $\Bbb P^1$, $h$ itself satifies the conditions of $l$.  

\proclaim
 {\rm (3.6)} {\bf Proposition:} Let ${\pi}:({\tilde U},E) \to (U,v)$ be a resolution of the singularity $v$ and let $(U,u_0)\subset ({\Bbb C}^N, 0)$. Then there is a linear function $h:{\Bbb C}^N \to {\Bbb C}$ such that 
$$
\di (h \circ {\pi})=Z+R
$$
where $R$ is reduced and meets $E$ transversely at smooth points of $E$. Moreover, if $R \cap E_i \neq \emptyset$, then $m_i=n_i$; and if $e \in R \cap E_i$, we may take $k=h$ in (1) above while if $e \notin R \cap E_i$, then we may choose $l=h$ in (1) above.
\endproclaim
{\bf Proof:} The idea is as follows; more details can be found in \cite{GS1, Lemma 2.1}. Let $\tau:B\ell(U)\to U$ denote the proper transform of $U$ in $B\ell(\Bbb C^N)$, the blow-up of $\Bbb C^N$ at the origin, and let $C\subset B\ell(U)$ be the reduced exceptional set, a curve in the fiber of $B\ell(\Bbb C^N)\to\Bbb C^N$over the origin. Then $C\subset \Bbb P^{N-1}$ and Bertini's Theorem says a generic hyperplane meets $C$ transversely in isolated smooth points of $C$. One gets such generic hyperplanes as the intersections of proper transforms $B\ell(H)$ of generic hyperplanes $H$ in $\Bbb C^N$ passing through zero with $\Bbb P^{N-1}$; so the desired $h$ is a linear function vanishing on such an $H$. Since $\frak{m}_v\cdot\C O_{B\ell(U)}$ is locally free of rank one on $B\ell(U)$ (a basic property of blowing up), and $h\circ\tau$ is a global section of it vanishing only along $R$, we have $\di (h\circ\tau)=\di (\tau^*\frak{m}_v)+R$. To begin the passage from $B\ell(U)$ to a resolution $\tilde U$ one must first normalize $B\ell(U)$, which requires a more careful, but still generic choice of $H$. Finally one may complete the resolution of $U$ with modifications away from the intersection points of $C$ with $B\ell(H)$ and then $\di (h\circ\pi)=\di (\pi^*\frak{m}_v)+R=Z+R$.

Since $R$ and $E$ are transverse at such an intersection point $e\in R\cap E$, there are local equations $\{v=0\}$ of $R$ and $\{u=0\}$ of $E$ so that $h=u^mv$ near $e$. A small perturbation $l$ of $h$ has the form $l=\delta u^m$, where $\delta$ is holomorphic and nowhere zero near $e$, since the set of $h$ above was generic. Now extracting an $m$-th root of $\delta$ and replacing $u$ with $u\delta^{1/m}$, we have $l=u^m$ and $h=\delta^{-1}u^mv$. If we replace $v$ by $\delta^{-1}v$ then we have coordinates $\{u,\,v\}$ on a neighborhood of $e$ and linear functions $k$, $h:\Bbb C^N\to\Bbb C$ such that $\phi:=l\circ\pi=u^m$ and $\psi:=h\circ\pi=u^mv$. By (3,3b) and (3.3c) above (in particular, the minimality property in (3.3c)), we are done.
\medskip 

The locus of the vanishing of the determinant of $\gamma$ gives a measure of its degeneracy. Now on $W$, $\gamma\sim d{\phi}d{\bar \phi}+d{\psi}d{\bar \psi}$
and the determinant of $d{\phi}d{\bar \phi}+d{\psi}d{\bar \psi}$ is $|\phi_u\psi_v-\phi_v\psi_u|^2$. A calculation using (2) and (3) above shows that $\phi_u\psi_v-\phi_v\psi_u$ locally defines the divisor
$$ 
D_{\gamma}=\sum (m_i+n_i-1)E_i=Z+N-E,\tag3.7
$$
so we call it the {\it degeneracy divisor} of $\gamma$. This calculation also shows that the volume form of $\gamma$ in $W$ is
$$
d{\tilde U}_{\gamma}\sim |{\phi}_u{\psi}_v-{\phi}_v{\psi}_u|^{2}d{\tilde U}\sim 
|u|^{2(m_i+n_i-1)}|v|^{2(m_j+n_j-1)}d{\tilde U}\tag3.8
$$
{where $d{\tilde U}:=du \wedge d{\bar u} \wedge dv \wedge d{\bar v}$ is the volume form in the Euclidean metric $du d\bar u+ dv d\bar v$. For any differential forms $\omega_1,\,\,\omega_2$ defined a.e. on $\tilde U$, let  
$$
<\omega_1,\omega_2>_{\gamma}
$$
denote as usual the  pointwise-defined inner product and  
$$
\|\omega\|_{\gamma}:=\left (\int_{\tilde U}<\omega,\omega>_{\gamma}\,d{\tilde U}_{\gamma} \right)^{1/2}
$$
the {\it $L^2$-norm of $\omega$}. Unless otherwise specified, we understand the pseudo-metric $\gamma$ on $\tilde U$ and omit the subscripts in such expressions, unless another metric is intended.

 Let $\C L^{p,q}_{\gamma}$ denote the sheaf of measurable forms on $\ti U$ which have locally finite $L_2$-norm. Notice that if $\tau$ is a differential form on $U-v$ and $\omega:=\pi^*\tau$, then
$$
\|\tau\|_g=\|\omega\|_\gamma\tag3.9
$$
so that the norm of a form on $U-v$, measured using the metric coming from the imbedding $(U,u_0)\subset({\Bbb C}^N, 0)$ is the same as the $L_2$-norm of its pullback to $\tilde U$. In particular, we have the equality of sheaves for each $p,q$,
$$
\pi_*\C L^{p,q}_{\gamma}=\C L^{p,q}\tag3.10
$$
where $\C L^{p,q}$ denotes the sheaf of $L_2$-forms on $U$. 
There will be advantages to working in $\tilde U$ rather than in $U$; indeed, some situations require it.

 \proclaim {\rm (3.11)}{\bf Definition:} The {\it Nash sheaf} is defined to be
$$
 \Cal N:=\varOmega^1_{{\tilde U},(2)}\otimes\Cal O_{\tilde U}(-D_{\gamma})
$$
where $\varOmega^1_{\tilde U,(2)}$ is the sheaf of 1-forms which have {\it locally} finite $L^2$-norm on $\tilde U$ and are holomorphic on $U-E$.
\endproclaim

If $i:U-E\inc U$ denotes the inclusion, then $\C N$ is the subsheaf of $i_*\varOmega^1_{U-E}$ defined by the local condition: if $\omega$ is defined near $u\in {\tilde U}$ where $\di (d_{\gamma})=D_{\gamma}$, then $\omega\in\Cal N$ if and only if $\|d^{-1}_\gamma{\omega}\|<\infty$. This implies that actually $\C N\subseteq\varOmega^1_U$: since $\gamma$ degenerates near $E$ with respect to any Hermitian (non-degenerate) metric $\mu$, $\|d^{-1}_\gamma{\omega}\|<\|{\omega}\|_\mu$ for any 1-form $\omega$, so the Laurent expansion of $\omega$ can have no polar part. For the same reason, $dw_i\in\C N$ for any set $w_1,\,w_2,\,\dots,\,w_N$ of coordinates on $\Bbb C^N$. 

B. Youssin independently noticed part {\it d.} of the folowing proposition (\cite{Y}), but for arbitrary varieties. In an appendix to this chapter, we will elaborate this point further, in particular in relation to (3.3), and give the  reason for the name ``Nash sheaf".

\proclaim {\rm (3.12)}{\bf Proposition}  
\roster
\ite a $\Cal N|U-E=\varOmega^1_{U-E}$
\ite b $\Cal N$ is locally free of rank 2: if $W$ is a neighborhood of $e\in \tilde U$ as in (3.3a) above, then $\{d\zeta_1, d\zeta_2\}$ from (3.3d) is an $\Cal O_W$-basis of $\Cal N^*(W)$; $\{d\phi, d\psi\}$ from (3.3a) is likewise a basis.
\ite c $\varOmega^1_{\tilde U}$ is a subsheaf of $\C N(N)$.
\ite d $\C L^{p,q}_{\gamma}=\C M(\Lambda^p\C N\otimes \Lambda^q\bar\C N\otimes \C O(D_\gamma))$, where, for any Hermitian bundle $B$ on $\tilde V$ and $\C B$ its sheaf of sections, $\C M(\C B)$ denotes its sheaf of measurable sections. 
\endroster
\endproclaim
{\bf Proof} Part {\it a.} is obvious and reduces {\it b.} to the case stated there: that $\C N|W$ is $\C O_W$ free with basis $\{d\zeta_1, d\zeta_2\}$; we also assume $W$ is a neighborhood of a crossing point $e\in E_i\cap E_j$ of $E$, the other case being similar. We may use any of the three quasi-isometric metrics in (3.3d) to determine whether a form on $W$ is $L_2$ with respect to $\gamma$, and here we choose $d\zeta_1d\bar  \zeta_1+d\zeta_2d\bar \zeta_2$. Let $d_\gamma=u^{(m_i+n_i-1)}v^{(m_j+n_j-1)}$, a local defining function for $D_\gamma$. To begin, note that $d\zeta_i\in\C N$: $\| d^{-1}_\gamma d\zeta_i\|<\infty$ since $\langle d\zeta_i,d\zeta_i\rangle=1$ and $|d_\gamma|^2 du \wedge d{\bar u} \wedge dv \wedge d{\bar v}\sim$ the volume form of $d\zeta_1d\bar  \zeta_1+d\zeta_2d\bar \zeta_2$. Since $d\zeta_1\wedge d\zeta_2$ vanishes only on $E$, $d\zeta_1$ and $d\zeta_2$ are $\C O_W$-independent. To show they generate, let $\omega\in\C N$ and write
$$
\omega=\alpha_1 d\zeta_1+\alpha_2 d\zeta_2
$$
where $\alpha_1$ and $\alpha_2$ are meromorphic. Then since $\|d^{-1}_\gamma\omega\|<\infty$ and $d\zeta_1$ and $d\zeta_2$ are pointwise orthonormal, the $\alpha_i$ are of finite $L_2$-norm in any non-degenerate metric, hence are holomorphic. That $\{d\phi, d\psi\}$ is also a basis is proved in the same way. Part {\it c.} asserts that $\|u^{-(m_i-1)}v^{-(m_j-1)}\omega\|<\infty$ for any holomorphic 1-form $\omega$, which holds because $ \|u^{n_i}v^{n_j}d\zeta_1\|$ and $ \|u^{n_i}v^{n_j}d\zeta_2\|$ are both finite. The proof of {\it d.}) is similar to that of {\it b.}).

\indent For the coming comparison between $\Cal N$ and ${\varOmega}^1_{\tilde U}(\lE)$, the sheaf of holomorphic 1-forms on $\tilde U$ with logarithmic singularities along $E$, it will be useful to have other local meromorphic sections (with poles along $E$) of $\varOmega^1_{\tilde U}$ with which to express elements of $\varOmega^1_{\tilde U}(\lE)$ and $\Cal N^*$.

Let $W$ be a neighborhood of $e\in E_i \cap E_j$ (resp., of $e\in E_i$, away from the crossings of $E$) with coordinates $u,\,v$ as in (3.3) above. Then we have the {\it logarithmic frame},
$$
\{\frac{du}{u},\,\frac{dv}{v}\}\quad\text{(resp.}\{\frac{du}{u},\,dv\}\text{)}\tag3.13
$$
the standard local basis for ${\varOmega}^1_{\tilde U}(\lE)$.
Referring now to the functions $\zeta_1$ and $\zeta_2$ in (3.3d), we define 
$$
 \zeta'_2:=\left\{
\aligned
&\zeta_2\text{, if $e$ is at a crossing of $E$}\\
&\zeta_2v^{-1}\text{, if $e$ is away from a crossing of $E$}\endaligned
\right.\tag3.14
$$
 and then the {\it logarithmic Nash frame} is
$$
\{\frac{d\zeta_1}{\zeta_1},\,\frac{d\zeta_2}{\zeta'_2}\}\tag3.15
$$
If we write a meromorphic 1-form $\omega$ on $W$ in the logarithmic frame
$$
\omega=f\frac{du}{u}+g\frac{dv}{v}\quad\text{(resp., }\,f\frac{du}{u}+g\,dv\text{)}\tag3.16
$$
then in the logarithmic Nash frame,
$$
\omega=\frac{n_jf-n_ig}{d}\,\frac{d\zeta_1}{\zeta_1}+\frac {m_ig-m_jf}{d}\,\frac{d\zeta_2}{\zeta'_2}\quad\text{(resp. }\frac{f-n_igv}{m_i}\,\frac{d\zeta_1}{\zeta_1}+g\frac{d\zeta_2}{\zeta'_2}\text{)}\tag3.17
$$
where $d=m_1n_2-m_2n_1$. Since $d\neq 0$ (resp., $m_i\neq 0$), it follows from this that
$$
\varOmega^1_{\tilde U}(\lE)(W)=\{k_1\frac{d\zeta_1}{\zeta_1}+k_2\frac{d\zeta_2}{\zeta'_2}\bigm| \text{ $k_1$ and $k_2$ are holomorphic in $W$}\}\tag3.18
$$
so $\{\frac{d\zeta_1}{\zeta_1},\,\frac{d\zeta_2}{\zeta'_2}\}$ is also a local basis for $\varOmega^1_{\tilde U}(\lE)$. And since $Z=\di (\zeta_1)$ and $\zeta'_2/\zeta_1$ is holomorphic in $W$, $\C N(Z)$ is a subsheaf of $\varOmega^1_{\tilde U}(\lE)$; in fact,
$$
\C N(Z)=\{k_1\frac{d\zeta_1}{\zeta_1}+k_2\frac{d\zeta_2}{\zeta'_2}\bigm| \text{ $k_1$ and $\frac{\zeta_1}{\zeta'_2}k_2$ are holomorphic in $W$}\}\tag3.19
$$

\proclaim
{\rm (3.20)}{\bf Proposition:} Let $\C I_E$ denote the ideal sheaf of $E$. There is an exact sequence of sheaves on $\tilde U$
$$
0\rightarrow\C N(Z-E) @>\alpha>>{\Cal I}_E{\varOmega}^1_{\tilde U}(\log\,E)@>\beta>>{\varOmega}^2_{\tilde U}\otimes{\Cal O}_{N-Z}\rightarrow0
$$
\endproclaim
\noindent (3.21) {\bf Remark:} Tensoring  the exact sequence with ${\Cal O}(E-Z)$ gives a description of the dual sheaf in terms of resolution data:
$$
0\rightarrow{\Cal N}\rightarrow{\Cal I}_E{\varOmega}^1_{\tilde U}(\log\,E)\otimes{\Cal O}_{\tilde U}(E-Z)\rightarrow{\varOmega}^2_{\tilde U}\otimes{\Cal O}_{N-Z}(E-Z)\rightarrow0
$$
{\bf Proof}: For the proof ${\varOmega}^i$ and $\Cal O$ will denote ${\varOmega}^i_{\tilde U}$ and ${\Cal O}_{\tilde U}$. The injection $\alpha$ is the tensor product of the inclusion $\C N(Z)\subseteq \varOmega^1(\lE)$ with $\C I_E$. To define $\beta$, recall from Proposition (3.6) the holomorphic function $h$ on $\tilde U$ such that $\di (h)=Z+R$, where $R$ is reduced and meets $E$ transversely and away from the crossings. 
Define
$$
{\tilde \beta}: {\Cal I}_E{\varOmega}^1(\log\,E)\rightarrow{\varOmega}^2\otimes{\Cal O}(R),\qquad {\tilde \beta}(\omega)=\omega\wedge\frac{dh}{h}
$$
{We first show that $\omega\wedge\frac{dh}{h}\in{\varOmega}^2\otimes{\Cal O}(R)$. In a neighborhood $W$ of $e\in E_i \cap E_j$ as in (3.3), we have $h=ku^{m_i}v^{m_j}$ (resp., $h=ku^{m_i}$ if $e\in E_i$ is away from a crossing and away from $R$), where $k$ is a nowhere-vanishing holomorphic function. Change $u$ by mutiplying it by the inverse of an $m_i$-th root of $k$, so that $h=u^{m_i}v^{m_j}$ (resp., $h=u^{m_i}$), for this choice of coordinates $\{u,\,v\}$ on $W$. Now let $\omega=k_1\,\frac{d\zeta_1}{\zeta_1}+k_2\,\frac{d\zeta_2}{\zeta'_2}\in{\Cal I}_E{\varOmega}^1(\lE)$, so that $(uv)^{-1}k_1$ and$(uv)^{-1}k_2$ are holomorphic (resp., $u^{-1}k_1$ and$u^{-1}k_2$ are holomorphic). Then 
$$
\align
{\tilde \beta}(\omega)=\omega\wedge\frac{dh}{h}&=(uv)^{-1}k_2\frac{uv\,d\zeta_2\wedge d\zeta_1}{\zeta'_2\zeta_1}\\
\text{(resp.}&=u^{-1}k_2\frac{uv\,d\zeta_2\wedge d\zeta_1}{\zeta'_2\zeta_1}\text{)}\endalign
$$
which is a holomorphic 2-form in $W$, since $\frac{uv\,d\zeta_2\wedge d\zeta_1}{\zeta'_2\zeta_1}$ is a smooth nowhere-vanishing multiple of $du\wedge dv$. Near a point $e \in E_i \cap R$, we can find coordinates $\{u,\,v\}$ so that $h=u^{m_i}v$ with $R=\{v=0\}$ and $\omega=f\,du+ug\,dv$. Then ${\tilde \beta}(\omega)=\omega\wedge\frac{dh}{h}=(fv^{-1}-m_ig)du \wedge dv$, which is clearly in ${\varOmega}^2\otimes{\Cal O}(R)$. These computations also show that $\tilde \beta$ is surjective.}

\indent Finally, define $\beta$ to be $\tilde \beta$ composed with the quotient map 

$$
{\varOmega}^2\otimes{\Cal O}(R)\rightarrow{\varOmega}^2\otimes{\Cal O}(R)/{\varOmega}^2\otimes{\Cal O}(R-N+Z)
$$

\noindent Since $N=Z$ at points of $R\cap E$, we have

$$
{\varOmega}^2\otimes{\Cal O}(R)/{\varOmega}^2\otimes{\Cal O}(R-N+Z)\cong{\varOmega}^2/{\varOmega}^2\otimes{\Cal O}(-N+Z)\cong{\varOmega}^2\otimes{\Cal O}_{N-Z}
$$

\indent Looking back now to the computation of $\tilde \beta(\omega)$ in the neighborhood $W$ of $e\in E_i \cap E_j$, we see that if $\tilde \beta(\omega)=0$, then $(uv)^{-1}k_2\in \C O(Z-N)$, so  $(uv)^{-1}k_2\frac{\zeta_1}{\zeta'_2}$ is holomorphic in $W$. Consequently, we see from our description of $\C N(Z)$ in (3.19) that $\omega\in\C N(Z-E) $
\medskip
\indent To state the first corollary of the Proposition, we need a definition.
\smallskip
\noindent {\rm (3.22)} {\bf Definition:} Let $\Cal F$ be an ${\Cal O}_{\tilde U}$-module. The Serre dual of $\Cal F$ is $\Cal F\sphat:=\C Hom_{
{\Cal O}_{\tilde U}}({\Cal F},{\varOmega}^2_{\tilde U})$. If ${\alpha} : {\Cal F} \to {\Cal G}$ is an ${\Cal O}_{\tilde U}$-homomorphism, then the Serre dual of $\alpha$ is ${\alpha}\sphat:=\C Hom_{{\Cal O}_{\tilde U}}({\alpha},{\varOmega}^2_{\tilde U})$.
\smallskip
\proclaim
{\rm (3.23)} {\bf Corollary} There is a short exact sequence of sheaves
$$
0\rightarrow{\varOmega}^1_{{\tilde U}}(\log\,E)@>\alpha\sphat>>\C N(N)\rightarrow{\Cal O}_{\tilde U}(N-Z)/{\Cal O}_{\tilde U}\rightarrow0
$$
where $\alpha \sphat$ is the Serre dual of $\alpha$.
\endproclaim

\indent {\bf Proof:} We make the notational conventions of the proof of Proposition (3.20). Wedge product induces  $\Cal O$-bilinear pairings
$$
\varOmega^1(\log\,E)\times\Cal I_E\varOmega^1(\log\,E)\to \varOmega^2\,\,\text{and}\,\,\C N(N)\times\C N(Z-E)\to \varOmega^2
$$
which are easily checked to be nonsingular in the sense that the induced $\Cal O$-homomorphisms are isomorphisms:
$$
\varOmega^1(\log\,E)@>\cong>>\Cal I_E\varOmega^1(\log\,E)\sphat\,\,\text{and}\,\,\C N(N)@>\cong>>\C N(Z-E)\sphat
$$
Thus, taking the Serre dual of the exact sequence of (3.20), we get an exact sequence
$$
0\rightarrow\varOmega^1(\log\,E)\to\C N(N)\to {\Cal Ext}^1_{\Cal O}(\varOmega^2\otimes{\Cal O}_{N-Z},\varOmega^2)\rightarrow0,
$$
and it remains to identify the rightmost term with ${\Cal O}(N-Z)/{\Cal O}$. By [AK, p. 74], we have
$$
\C Ext^1_{\Cal O}(\varOmega^2\otimes{\Cal O}_{N-Z}, \varOmega^2)@>\cong>>\C Ext^1_{\Cal O}(\Cal O_{N-Z},\varOmega^2)\otimes(\varOmega^2)^*
$$
where $(\varOmega^2)^*:={\Cal Hom}(\varOmega^2,{\Cal O})$. By Serre-Grothendieck duality [AK, p. 13]
$$
{\Cal Ext}^1_{\Cal O}(\Cal O_{N-Z},\varOmega^2)\cong{\Cal Hom}_{\Cal O_{N-Z}}({\Cal I_{N-Z}}/{\Cal I^2_{N-Z}},\varOmega^2\otimes{\Cal O_{N-Z}})
$$
where $\Cal I_{N-Z}$ is the ideal sheaf of $N-Z$, so that $\Cal O_{N-Z}\cong\Cal O/\Cal I_{N-Z}$. Hence, since ${\Cal I_{N-Z}}/{\Cal I^2_{N-Z}}$ is locally free and $\varOmega^2$ is locally free of rank one,
$$
\multline
\shoveright{{\Cal Ext}^1_{\Cal O}(\Cal O_{N-Z},\varOmega^2)\otimes(\varOmega^2)^*\cong{\Cal Hom}_{\Cal O_{N-Z}}({\Cal I_{N-Z}}/{\Cal I^2_{N-Z}},\varOmega^2\otimes{\Cal O_{N-Z}})\otimes(\varOmega^2)^*\cong}\\
{\Cal Hom}_{\Cal O_{N-Z}}({\Cal I_{N-Z}}/{\Cal I^2_{N-Z}},{\Cal O_{N-Z}})\otimes\varOmega^2\otimes(\varOmega^2)^*\cong{\Cal Hom}_{\Cal O_{N-Z}}({\Cal I_{N-Z}}/{\Cal I^2_{N-Z}},{\Cal O_{N-Z}})
\endmultline
$$
and this last is isomorphic to $(\Cal I_{N-Z})^{-1}/\Cal O ={\Cal O}(N-Z)/\Cal O$ as claimed.
\smallskip
\proclaim{\rm (3.24)}{\bf Corollary} 
\roster
\ite a The ${\Cal O}_{\tilde U}$-homomorphisms $\alpha$ and $\alpha\otimes{\Cal O_{\tilde U}}(E-Z)$ induce surjections 
$$
\alpha_*:H^1({\tilde U};\C N(Z-E))\twoheadrightarrow H^1({\tilde U};\Cal I_E\varOmega^1_{\tilde U}(\log\,E))
$$
and
$$
H^1({\tilde U};\C N)\twoheadrightarrow H^1({\tilde U};\Cal I_E\varOmega^1_{\tilde U}(\log\,E)\otimes{\Cal O}(-Z))
$$
\ite b The ${\Cal O}_{\tilde U}$-homomorphisms $\alpha\sphat$ and $\alpha\sphat\otimes{\Cal O}_{\tilde U}(Z-E)$ induce injections
$$
\alpha\sphat_*:H^1({\tilde U};\varOmega^1_{\tilde U}(\log\,E))\rightarrowtail H^1({\tilde U};\C N(N))
$$
and
$$
H^1({\tilde U};\varOmega^1_{\tilde U}(\log\,E)\otimes{\Cal O}_{\tilde U}(Z-E))\rightarrowtail H^1({\tilde U};\C N(D_\gamma))
$$
and isomorphisms
$$
\Gamma({\tilde U};\varOmega^1_{\tilde U}(\log\,E))@>\cong>> \Gamma({\tilde U};\C N(N))\,\,\,\text{and}\,\,\,\Gamma({\tilde U};\varOmega^1_{\tilde U}(\log\,E)\otimes{\Cal O}_{\tilde U}(Z-E))\cong\Gamma({\tilde U};\C N(D_\gamma))
$$
\endroster
\endproclaim
\indent {\bf Proof:} Part {\it b.} is evidently equivalent to
$$
\Gamma({\tilde U};{\Cal O}(N-Z)/{\Cal O})=0\quad\text{and}\quad \Gamma({\tilde U};{\Cal O}(N-E)/{\Cal O}(Z-E))=0
$$
which are [P, 4.1] (where $D=N-Z$) and [PS, p.619] (where we take $D=N$ and observe that $H^0({\tilde U};{\Cal O}(Z-E))\to H^0({\tilde U};{\Cal O}(N-E)) $ is surjective since, by [P, 4.2], $H^0({\tilde U};{\Cal O})\to H^0({\tilde U};{\Cal O}(N-E))$ is.)

\indent To prove part {\it a.}, let $\Cal S$ be a coherent $\Cal O_{\tilde U}$-module and let $\Cal I$ denote the ideal sheaf of the complete intersection scheme $N-Z$. Then
$$
H^1({\tilde U}; \varOmega^2\otimes{\Cal O}_{N-Z}\otimes\Cal S)\cong H^1(N-Z;\varOmega^2/{\Cal I}\varOmega^2\otimes{\Cal S}/{\Cal I}{\Cal S})
$$
which by Grothendieck duality on $N-Z$ is isomorphic to
$$ 
Hom_{\Cal O_{N-Z}}(\varOmega^2/{\Cal I}\varOmega^2\otimes{\Cal S}/{\Cal I}{\Cal S},\omega_{N-Z})
$$
where $\omega_{N-Z}$ is the dualizing sheaf of $N-Z$. Since 
$$
\omega_{N-Z}:={\Cal Ext}^1_{\Cal O}({\Cal O}_{N-Z}, \varOmega^2)
$$
which is, by [AK, p. 13]
$$
{\Cal Hom}_{\Cal O_{N-Z}}(\Cal I/{\Cal I}^2, \varOmega^2/{\Cal I}\varOmega^2)
$$
we have
$$
\Cal Hom_{\Cal O_{N-Z}}(\varOmega/{\Cal I}\varOmega^2\otimes{\Cal S}/{\Cal I}{\Cal S},\omega_{N-Z})\cong \Cal Hom_{\Cal O_{N-Z}}(\varOmega^2/{\Cal I}\varOmega^2\otimes{\Cal S}/{\Cal I}{\Cal S}\otimes\Cal I/{\Cal I}^2,\varOmega^2/{\Cal I}\varOmega^2)
$$
Since $\varOmega^2/{\Cal I}\varOmega^2$ is locally ${\Cal O}_{N-Z}$-free, this last is 
$$
\Cal Hom(\Cal S\otimes\Cal I/{\Cal I}^2,{\Cal O}_{N-Z}).
$$
When $\Cal S={\Cal O}$ (resp. $\Cal S={\Cal O}(E-Z)$), we get, since $N-Z$ is a complete intersection scheme,
$$
H^0({\tilde U};{\Cal O}(N-Z)/\Cal O)\,\,\,\,\,\,(\text{resp.}H^0({\tilde U};{\Cal O}(N-E)/{\Cal O}(Z-E)))
$$
which were shown to vanish in the proof of part {\it b.} above.

\ind Now combining parts {\it a.} and {\it b.} of this corollary, we get a commutative diagram, in which the left vertical map is surjective and the right, injective:
$$
\CD
H^1(\tilde U;\C I_E\varOmega^1(\lE)) @>>> H^1(\tilde U;\varOmega^1(\lE))\\
@A{\alpha_*}AA             @VV{\alpha\sphat_*}V\\
H^1(\tilde U;\C N(Z-E)) @>>> H^1(\tilde U;\C N(N))
\endCD\tag3.25
$$
>From this the first statement in the following corollary is immediate, and the second and third are proved similarly.
\proclaim {\rm (3.26)}{\bf Corollary} \roster
\ite a The homomorphisms $\alpha$ and $\alpha\sphat$ induce an isomorphism
$$
\im (H^1(\tilde U;\C I_E\varOmega^1(\lE))\to  H^1(\tilde U;\varOmega^1(\lE)) @>\cong>> \im (H^1(\tilde U;\C N(Z-E))\to H^1(\tilde U;\C N(N));
$$
\ite b $\alpha$ and $\alpha\sphat\otimes\C O(Z-E)$ induce an isomorphism
$$
\im (H^1(\tilde U;\C I_E\varOmega^1(\lE))\to  H^1(\tilde U;\varOmega^1(\lE)\otimes\C O(Z-E)) @>\cong>> \im (H^1(\tilde U;\C N(Z-E))\to H^1(\tilde U;\C N(D_\gamma))
$$
\ite c $\alpha\otimes\C O(E-Z)$ and $\alpha\sphat$ induce an isomorphism
$$
\im (H^1(\tilde U;\C I_E\varOmega^1(\lE)\otimes\C O(E-Z))\to  H^1(\tilde U;\varOmega^1(\lE)) @>\cong>> \im (H^1(\tilde U;\C N)\to H^1(\tilde U;\C N(N))
$$

\endroster
\endproclaim

\ind To conclude this section recall the short exact sequences of sheaves on $\tilde U$:
$$
0\to \C I_E\varOmega^1(\lE) @>>> \varOmega^1 @>>> \oplus i_*\varOmega^1_{E_i}\to 0\tag3.27
$$ 
and
$$
0\to \varOmega^1 @>>> \varOmega^1(\lE) @>>> \oplus i_*\C O_{E_i}\to 0
$$
These give rise to exact sequences
$$
H^1(\tilde U;\C I_E\varOmega^1(\lE)) @>f>> H^1(\tilde U;\varOmega^1) @>g>> \oplus H^1(E_i;\varOmega^1_{E_i})\tag3.28
$$
and
$$
\oplus H^0(E_i;\C O_{E_i}) @>d>> H^1(\tilde U;\varOmega^1) @>h>> H^1(\tilde U;\varOmega^1(\lE))
$$
\proclaim {\rm (3.29)}{\bf Proposition} 
\roster
\ite a $H^0(\tilde U;\varOmega^1) @>\cong>>  H^0(\tilde U;\varOmega^1(\lE))$
\ite b $f(\ker hf)=0$
\endroster
\endproclaim
\ind {\bf Proof} These follow easily from the fact that $gd$ is an isomorphism, since it can be identified with the adjoint of the cup product pairing on $H^2(\ti U;\Bbb C)$, which is well-known to be nonsingular.
\bigskip

{\bf Appendix: The Nash bundle and Hsiang-Pati coordinates}
\smallskip 
\indent Let $M^m$ be a smooth quasiprojective variety and let $Gr(n,TM)$ denote the bundle with fiber $Gr(n,m)$, the Grassmanian of $n$-planes in $\Bbb C^m$, associated to the tangent bundle $TM$ of $M$. Let $i:W\inc M$ be a subvariety (always reduced) of dimension $n$ and let $\hat W$ be the closure in $Gr(n,TM)$ of the image of the section over $W-\SW$ defined by the derivative $di:T(W-\SW)\to TM$. Then the bundle projection $Gr(n,TM)\to M$ restricts to a proper algebraic map $\hat\pi:\hat W\to W$ whose restriction to $\hat W-\hat\pi^{-1}(W-\SW)$ is a biholomorphism onto $W-\SW$. The pair $(\hat W,\hat\pi)$ is called {\it the Nash blow-up of} $W$; it is independent, in the obvious sense, of the choice of imbedding $i$. The canonical $n$-plane bundle over $Gr(n,TM)$ restricts to an $n$-plane bundle $N\hat W$ over $\hat W$, whose restriction to $\hat W-\hat\pi^{-1}(W-\SW)$ is the tangent bundle. $N\hat W$ is called {\it the Nash bundle} and is also intrinsic.

Now one may get a resolution of singularities $\pi:\ti W\to W$ by resolving the singularities of $\hat W$, say $\ti\pi:\ti W\to \hat W$, and setting $\pi:=\hat\pi\circ\ti\pi$. The following result was told to us by R. MacPherson (see also \cite{GS2}, where $N\ti W$ is called the {\it generalized} Nash bundle) It gives bundle data on $\ti W$ equivalent to the existence of such a factorization and allows us to extend the notion of Nash bundle to such $\ti W$.   
\proclaim{\rm (A3.1)}{\bf Proposition} Let $i:W\inc M$ be a subvariety and let $\pi:\ti W\to W$ be a resolution of singularities. 
\roster
\ite a There is a proper algebraic map $\ti\pi:\ti W\to \hat W$ such that $\pi=\hat\pi\circ\ti\pi$ if and only if there is a bundle $N\ti W$ over $\ti W$ and a bundle map $n:T\ti W\to N\ti W$ such that the tangent map $d(i\circ\pi):T\ti W\to TM$ factors
$$
d(i\circ\pi)=n\circ m
$$
where $m:N\ti W\to TM$ is a bundle map which is injective on fibers and covers $i\circ\pi$.
\ite b If the pair $(N\ti W,n)$ exists as in (1), then it is unique.
\endroster
\endproclaim
\proclaim{\rm (A3.2)}{\bf Definition} $N\ti W$ is called the Nash bundle (of $\pi:\ti W\to W$).
\endproclaim

We assume this result for the time being; a proof is given below.

Let now $\pi:\ti W\to W$ be a resolution which factors through $\hat W$. Let $\ti w\in \ti W$ and let $\{w_1,\dots,w_m\}$ be local holomorphic coordinates for $M$ centered at $w=\pi(\ti w)\in M$. Then the 1-forms $\{dw_1,\dots,dw_m\}$  can be pulled back to sections $\{n^*(dw_1),\dots,n^*(dw_m)\}$ of $N\ti W^*$. Since $n^*:TM^*\to N\ti W^*$ is surjective on fibers, the tensor
$$
\gamma_N:=\sum n^*dw_i\otimes n^*d\bar w_i\tag A3.3
$$
defines a (nonsingular) Hermitian metric on $N\ti W$ which restricts to the singular Hermitian metric $\gamma$ on $T\ti W$ pulled up from the metric $\sum dw_i\otimes d\bar w_i$ on a neighborhood of $w\in M$. 

Now suppose that the exceptional set of $\pi$ is a divisor $E$ with normal crossings. Then it is easy to see that, in local holomorphic coordinates near each $\ti w$, the volume form of $\gamma$ is
$$
dW_\gamma=|d_\gamma|^2 dW_\sigma
$$
where $d_\gamma=0$ is the local defining equation for a divisor $D_\gamma$ supported in $E$ and $dW_\sigma$ is a local (nonsingular) Hermitian form on on $\Lambda^nN\ti W$. In case $\dim W=2$, $D_\gamma$ will turn out to be the degeneracy divisor of $\gamma$ defined in (3.7) above. It is now immediate that
$$
\C L^{p,q}_{\gamma}=\C M(\Lambda^p\C N\otimes \Lambda^q\bar\C N\otimes \C O(D_\gamma))
$$
The proof is the same as that of part {\it d.} of Proposition (3.12 ) above. 

Observe that if the pair $(N\ti W,n)$ exists as in {\it a.}, then $\ti\pi:\ti W\to \hat W$ factors $\pi$, where 
$$
\ti\pi(x):=(\pi(x),m(N\ti W_x)\in Gr(n,TM)_x.
$$
To prove the rest of the theorem and make the connection to Hsiang-Pati coordinates, we need to work the context of sheaves. We use without comment the well-known eqivalence between the categories of locally free sheaves and of algebraic vector bundles given in \cite{Ha, pp.128-129}. For example, using standard properties of the sheaf $\varOmega^1$ of differentials, the translation to the language of sheaves of part {\it a.} of the proposition is:

\proclaim  {\rm (A3.4)}A resolution $\pi:\ti W\to W$ factors through the Nash blow-up if and only if there is a pair $(\C N,\nu:\C N\to\varOmega^1_{\ti W})$, where $\C N$ is a locally free sheaf of rank $n$ on $\ti W$ and $\nu|\ti W-\pi^{-1}(\SW)$ is an isomorphism, such that the canonical map 
$$
\delta:\pi^*\varOmega^1_{W}\to\varOmega^1_{\ti W}
$$
factors through $\nu$,
$$
\delta=\nu\circ\mu
$$
where $\mu:\pi^*\varOmega^1_{W}\to\C N$ is is a surjective morphism of sheaves on $\ti W$.
\endproclaim
Notice that our use of $\C N$ here is consistent with that of Definition (3.11): the pair $(\C N, \pi:\ti U\to U)$, where $\C N$ is the sheaf on $\ti U$ defined there, satisfies the condition above: we observed (following (3.11)) that $\C N$ is a subsheaf of $\varOmega^1_{\ti U}$ and that the canonical $\C O_{\ti U}$-morphism $\pi^*\varOmega^1_U\to\varOmega^1_{\ti U}$ factors through that $\C N$. 

We next recall the notion of blowing up a coherent sheaf; our discussion comes directly from \cite{NA}.

Let $\C F$ be a coherent sheaf on a variety $W$ 
\proclaim{\rm (A3.5)}{\bf Definition} We call $\beta:\hat W\to W$ the blow-up of $W$ relative to $\C F$ if
\roster
\ite a $\beta$ is birational and proper,
\ite b $\beta^*\C F/\T\beta^*\C F$ is locally free on $\hat W$ and
\ite c if $\rho: Z\to W$ also satifies (1) and (2), then there is a unique regular map $Z\to \hat W$ which factors through $\beta$:
$$
\rho=\phi\circ\beta 
$$
and $\phi^*(\beta^*\C F/\T\beta^*\C F)\cong\rho^*(\C F/\T\rho^*\C F)$.
\endroster
\endproclaim
Here is a local construction of the blow-up of $W$ relative to $\C F$, which shows that it exists. (In fact, it is also unique, \cite{R}.) Let $U\inc W$ be an open affine subset over which there is an $\C O_U$-surjection 
$$
\C O^N_U\to\C F|U\to 0.\tag A3.6
$$
Let $U^\circ\inc U$ be the open dense subset where $\C F|U^\circ$ is locally free, of rank $r$ say. From this we get a section
$$
\sigma:U^\circ\to Gr(r,\epsilon^N_U)\cong U\times Gr(r,N)\tag A3.7
$$
since the surjection $\C O^N_U\to\C F|U$ of locally free sheaves on $U^\circ$ corresponds to an injection of vector bundles 
$$
\Bbb V(\C F|U^\circ)\inc\epsilon^N
$$
\proclaim{\rm (A3.8)}{\bf Proposition} \cite{Ro} 
Let $\hat U:=\overline{\sigma(U^\circ)}$, $\hat E=\hat U-\sigma(U^\circ)$ and let $\beta:\hat U\to U$ and $\gamma:\hat U\to Gr(r,N)$ be the induced maps. Then $\beta$ is the blow-up of $\C F|U$, $\beta|\hat U-\hat E:\hat U-\hat E\to U^\circ$ is biholomorphic and
$$
\beta^*\C F/\T\beta^*\C F\cong\gamma^*\C Q(r,N)
$$
where $\C Q(r,N)$ is the universal quotient sheaf over $Gr(r,N)$.
\endproclaim
Evidently, these constructions over affine $U$ covering $W$ patch together
to show that the the blow-up of $W$ relative to $\C F$ exists.
We call $\hat E$ the {it exceptional set} of $\beta$.
\medskip
\noindent(A3.9) {\bf Examples} \roster\ite a Let $\C F=\C I$ be an ideal sheaf in $\C O_W$. Here $r=1$ and $\beta$ is the blow-up of the ideal $\C I$(\cite{Ha,Rie}). In particular, if $\C I$ is the maximal ideal $\frak m_w$  of a point $w$ in $W$, then $\beta$ is the blow-up of $W$ at $w$, denoted $\tau:B\ell(W)\to W$ in the proof of (3.6) above: $U^\circ=U-w$, a generator of the stalk of $\frak m_w$ over $x\in U^\circ$ is the linear form whose zero set is perpendicular at $x$ to the secant line $\overarrow{xw}$, so a section sends $x\in U^\circ$ to $(x,\overarrow{xw})\in  U\times Gr(1,N)=\Bbb P^{N-1}$.
\ite b Let $\C F=\varOmega^1_W$ be the sheaf of differentials on $W$. Here $r=\dim W$ and it is clear that $\beta$ is the Nash blow-up $\hat\pi:\hat W\to W$ (\cite{T, III.1.2, Remark 3}). Set
$$
\C N_{\hat W}:=\hat\pi^* \varOmega^1_W/\T\hat\pi^*\varOmega^1_W.
$$
Then it follows from (A3.8) that $\C N_{\hat W}$ is the sheaf of sections of the dual $N\hat W^*$ of the Nash bundle. 
\endroster
\medskip
We can now prove Proposition (A3.1), in the form (A3.4). Suppose that $\pi=\hat\pi\circ\ti\pi:\ti W\to W$ is a resolution of singularities of $W$ which factors through the Nash blow-up $\hat\pi:\hat W\to W$ and set $\C N_{\ti W}:=\ti\pi^*\C N_{\hat W}$. Then the canonical map of sheaves on $\hat W$, $\hat\pi^*\varOmega^1_W\to\C N_{\hat W}$ gives rise to $\pi^*\varOmega^1_W=(\hat\pi\circ\ti\pi)^*\varOmega^1_W\to\ti\pi^*\C N_{\hat W}=\C N_{\ti W}$ which evidently divides out the torsion subsheaf. Since the canonical map $\delta$ kills torsion, $\delta$ factors uniquely through an $\C O_W$ map $\nu:\C N_{\ti W}\to \varOmega^1_{\ti W}$. Conversely if a pair $(\C N,\nu:\C N\to\varOmega^1_{\ti W})$ satifies the condition in (A3.4), then $\mu$ is the morphism which divides out torsion. By Proposition (A3.8) a factorization $\pi=\hat\pi\circ\ti\pi$ exists such that $\ti\pi^*\C N_{\hat W}=\C N_{\ti W}$. This completes the proof of (A3.1), including part {\it b.}, which follows from the uniqueness in (A3.5). 

Now let $D_r\subseteq\Bbb C^m$ be a subspace of codimension $r$. Let 
$$
S(D_r):=\{E^r\in Gr(r,m)|\dim (E^r\cap D_r)\geq 1\}\tag A3.10
$$
This is a {\it Schubert variety}, of codimension 1 in $Gr(r,m)$ and denoted $c_1(D_r)$ in \cite{LeT}; in \cite{NA}, it is $S(D)$ where $s=1$, $D=(D_r)$ and $a=(1)$ in the notation used there. Let $\beta:\hat W\to W$ be the blow-up of $W$ relative to $\C F$ and $\gamma:\hat W\to Gr(r,m)$, the canonical map. The following is a consequence of the transversality theorem of Kleiman (\cite{Kl, (.)}
\proclaim {\rm (A3.11)}{\bf Proposition} With notation as above, we have for generic $D_r$ in $Gr(m-r,m)$
\roster
\ite a $\gamma^{-1}S(D_r)\cap (\hat W-\hat E)$ is smooth and dense in $\gamma^{-1}S(D_r)$,
\ite b $\gamma^{-1}S(D_r)$ is either empty or has codimension $1$ in $\hat W$ and 
\ite c $\gamma^{-1}S(D_r)\cap\hat E$ is either empty or has codimension one in $\hat E$, can be arranged to miss any given finite set of points of $\hat E$, and, if
$\dim \hat E=2$, then $\gamma^{-1}S(D_r)\cap\hat E$ consists of isolated smooth points of $\hat E$.
\endroster\endproclaim
If $p:\Bbb C^m\to\Bbb C^r$ is any linear projection with $\ker p=D_r$, then $p$ induces a trivialization of the universal sheaf restricted to $Gr(r,m)-S(D)$, hence also, by (rossi), of $\beta^*\C F/\T\beta^*\C F$ over  $\hat W-\gamma^{-1}S(D_r)$. We show next that it is this trivialization that is the source of the linear projection $\Bbb C^N\to\Bbb C^2$ in (3.3). 
\medskip
\noindent (A3.12) {\bf Examples} Return to the two examples (A3.9). Since it is all we need in what follows, we assume that $W=U$, the neighborhood of a singular point $v\in V$ considered throughout \S 3.
\roster
\ite a If $\tau:(B\ell(U),C)\to (U,v)$ is the blow-up of $U$ at $v$, then $Gr(r,m)=\Bbb P^{m-1}$, $S(D)$ is a hyperplane $H$ and the trivilization of 
$$
\tau^*\frak m_w/\T\tau^*\frak m_w=\frak m_w\C O_{B\ell(U)}
$$
over $\gamma^{-1}(\Bbb P^{m-1}-S(D))$ is induced by a non-trivial linear function $h:\Bbb C^m\to\Bbb C$ where $H=(h)$ and is just the global section $h\circ\tau$. The reader may now recast the proof of (3.6) using Prop. (A3.11).
In particular, if $\ti\pi:(\ti U,E)\to (B\ell(U),C)$ is a resolution of singularities, $e\in E$, $\pi:=\check\pi\circ\ti\pi$ and $\pi(e)\notin H\cap C$, then $h\circ\pi$ is a generator of $\frak m_v\C O_{\ti U}$ near $e$.

\ite b If $\hat\pi:(\hat U,\hat E)\to (U,v)$ is the Nash blow-up and $k,\,l:\Bbb C^m\to\Bbb C$ are linear functions such that $\ker k\cap\ker l=D_2$, then $\C N_{\hat U}$ is trivialized over $\hat U-\gamma^{-1}S(D_2)$ by the global sections $\{\mu(\hat\pi^*dk),\mu(\hat\pi^*dl)\}$, where $\mu:\hat\pi^*\varOmega^1_U\to\C N$ is the factor of the canonical map $\delta$ in (A3.4). Hence if $\ti\pi:(\ti U,E)\to(\hat U,\hat E)$ is a resolution of singularities, $e\in E$, $\pi:=\hat\pi\circ\ti\pi$ and $\C N_{\ti U}:=\ti\pi^*\C N_{\hat U}$ and $\pi(e)\notin S(D)\cap \hat E$, then $\{d(k\circ\pi), d(l\circ\pi) \}$ is a basis for $\C N_{\ti U}$ near $e$. Remembering that $M$ in Prposition (A3.1) is a neighborhood of the origin in $\Bbb C^N$ so that $TM$ is identified with $\Bbb C^N$, we can equivalently say that $(k,l)\circ m:N\ti U\to \Bbb C^2$ is in isomorphism on the fibers of $N\ti U$ near $e$. 
\endroster

We must now combine these examples. To do this consider the commutative diagram
$$
\CD
(\check U,\check E) @>\check\tau>> (\hat U,\hat E)\\
@V\Check\pi VV                         @VV{\hat\pi} V   \\
(B\ell(U),C) @>\tau>> (U,v)
\endCD\tag A3.13
$$
where $(\check U,\check E)$ is the fiber product of $(B\ell(U),C)$ and $(\hat U,\hat E)$.
In fact, $\tau\circ\Check\pi =\hat\pi\circ\check\tau$ is the blow-up of $U$ corresponding the the sheaf (module in this case) $\frak m_v\oplus\varOmega^1_U$. Then $\check\tau$ is the blow-up of $\hat E$, $\check\pi$ is one-to-one except over a finite subset $B(C)$ of $C$ and $\check\tau$ is finite-to-one outside a finite subset $B(\hat E)$ of $\hat E$. Hence outside $B(C)$, $\check\pi$ factors the normalization of $B\ell(U)$; similarly for $B(\hat E)$ and $\check\tau$. Finally, let 
$$
\bar n:(\bar U,\bar E)\to (\check U,\check E)
$$ 
be the normalization. Then $\bar U$ has only isolated singularities, outside $(\check\pi\circ n)^{-1}B(C)$ is the normalization of $B\ell(U)-B(C)$ and, outside $(\check\tau\circ n)^{-1}B(\check E)$, is the normalization of $\hat U-B(\hat E)$.

Now suppose given $F(C)$, a finite subset of $C$ and $F(\hat E)$, a finite subset of $\hat E$. Choose the linear function $h$ in Example (A3.12a) (or the proof of Proposition (3.6)) so that the corresponding hyperplane also misses $B(C)\cup F(C)$; and choose the linear functionals $k$ and $l$ in part {\it b.} of (A3.12) so that the Schubert variety $S(\ker k\cap\ker l)$ misses $B(\hat E)\cup F(\hat E)$. Such choices are possible and generic, according to (A3.11c). 

\proclaim{\rm (A3.14)}{\bf Proposition} Let $\ti\pi:(\ti U,E)\to (\check U,\check E)$ be a resolution of singularities factoring through the normalization $\bar U $ of $\check U$, let $e\in E$ and set $\pi=(\hat\pi\circ\check\tau)\circ\ti\pi=(\tau\circ\check\pi)\circ\ti\pi$. Then there is a hyperplane $H=\ker h\in Gr(N-2,N)=\Bbb P^{N-1}$ and a codimension two plane $D=\ker k\cap\ker l\in Gr(N-2,N)$ such that 
\roster
\ite a $h\circ\pi$ generates $\frak m_v\C O_{\ti U}$ near $e$,
\ite b $\{d(k\circ\pi),d(l\circ\pi)\}$ generates $\C N_{\ti U}$ near $e$ and 
\ite c $D\subset H$
\endroster
\endproclaim
{\bf Proof} The discussion above, where $F(C)=\check\pi\circ\ti\pi(e)$ and $F(\hat E)=\check\tau\circ\ti\pi(e)$, proves {\it a.} and {\it b.}

To prove {\it c.}, begin with the following simple fact: Let $D\in Gr(N-2,N)$, let $\C G$ be a neighborhood of it in $ Gr(N-2,N)$ and let $D\subset H\in Gr(N-1,N)$. Then any $H'$ sufficiently close to $H$ is contained in some $D'\in \C G$. Now call $D\in Gr(N-2,N)$ {\it good} if it misses $B(\hat E)\cup F(\hat E)$ and $ H\in Gr(N-1,N)$ good if it misses $B(C)\cup F(C)$. The sets of good planes are dense and open in there respective Grassmannians by the discussion above. Let $D$ be a good plane in $Gr(N-2,N)$. Then there is a hyperplane $H\supset D$ and a sequence $\{H_i\}$ of good hyperplanes converging to $H$. Now choose a neighborhood $\C G$ of $D$ consisting of good $D$'s. Then the simple fact says we can find some $D'\in \C G$ and some $H_j$ with $ D'\subset H_j$. 
\medskip
 We can now use these blow-up and transversality considerations to derive (3.3) in its proper context.

\proclaim{\rm (A3.15)}{\bf Corollary} \cite{HP} Let $v\in V$ be an isolated singular point on a complex surface. Then there is a neighbohood $U$ of $v$ and a resolution of singularities $\pi:(\ti U,E)\to(U,v)$ such that for each $e\in E$ there are linear functions $k\,,l:\Bbb C^N\to\Bbb C$ satisfying properties {\it a.}-{\it c.} of (3.3).
\endproclaim
\demo{\bf Proof} It follows from part {\it c.} of Proposition (A3.14)) that we may take $k=h$ in parts {\it a.} and {\it b.} We take $E:=\sum m_iE_i$ to have simple normal crossings and suppose $e\in E_1\cap E_2$ is at a crossing; the proof in case $e$ is a simple point of $E$ is similar. Then we have, for $\phi:=k\circ\pi$ and $\psi:=l\circ\pi$ and suitable local coordinates $\{u,v\}$ with $E_1=\{u=0\}$ and $E_2=\{v=0\}$ near $e$,
$$
\phi=u^{m_1}v^{m_2},\qquad\qquad d\phi\wedge d\psi=d_\gamma du\wedge dv,
$$
where $d_\gamma=0$ is a local defining eequation for the degeneracy divisor of $\gamma$. Let us now for the sake of convenience work with elements of $\hat\C O$, the ring of germs of holomorphic functions at $e$. The fact that $\gamma$ degenerates only along $|E|$ means that
$$
d_\gamma=\mu u^{d_1}v^{d_2},\qquad \mu\in\hat\C O^*
$$
for some non-negative integers $d_1$ and $d_2$.

Since 
$$
d\phi\wedge d\psi=u^{m_1-1}v^{m_2-1}(m_1v\psi_v-m_2u\psi_u)
$$
we will need the following simple lemma whose proof is left to the reader.
\proclaim{\bf Lemma} Let $D_{m_1,m_2}:\hat\C O\to\hat\C O$ be the $\Bbb C$-derivation
$$
D_{m_1,m_2}g=m_1vg_v-m_2ug_u
$$
Then 
$$
\ker D_{m_1,m_2}=\Bbb C\{z\}
$$
the ring of convergent power series in $z:=u^{\frac{m_1}{(m_1,m_2)}}v^{\frac{m_2}{(m_1,m_2)}}$, where $(m_1,m_2)$ is the greatest common divisor of $m_1$ and $m_2$; and $\im  D_{m_1,m_2}$ consists of convergent power sreies 
$$
\sum r_{a,b}u^av^b
$$
where $m_1b-m_2a\neq 0$. 
\endproclaim
Using the Lemma write
$$
\psi=\sum s_iz^i+\sum r_{a,b}u^av^b
$$
Since $\psi$ must vanish along $E$, but to no lower order than $\phi$, we have $i\geq (m_1,m_2)$, $a\geq m_1$ and $b\geq m_2$ in these sums. We can now compute
$$
u^{m_1-1}v^{m_2-1}\sum(m_1b-m_2a)r_{a,b}u^av^b=u^{m_1-1}v^{m_2-1}D_{m_1,m_2}\psi=\mu u^{d_1}v^{d_2}
$$
and hence
$$
\sum r_{a,b}u^av^b=\nu u^{d_1-m_1+1}v^{d_2-m_2+1}
$$
for some $\nu\in\hat\C O^*$. Set
$$
n_i:=d_i-m_i+1
$$
Then since $a\geq m_1$ and $b\geq m_2$ for all $a$ and $b$, $n_i\geq m_i$ and the proof is complete.

\enddemo

\newpage

\heading \S 4: The cohomological Hodge structure in dimension two.\endheading

The main purpose of this section is to prove Theorem B: we construct the filtered quasi-isomorphism  $\Lambda:(\hat\C L^{\cdot},\C F^{\cdot}) @>\cong>> (\hat\C A^{\cdot},\C F^{\cdot})$ and verify that the canonical maps ${\kappa}^p_{\C L}:gr^p\hat\C L^{\cdot}_V\to \hat\C L^{p,\cdot}$ are isomorphisms for each $p$. In addition, we will compute ((4.19)-(4.22)) the $L_2-\db$-comology groups $H^{p,q}_B(V)$, where $\dim V= 2$ and $B=D$, $D/N$ and $N$. Throughout this chapter we keep the notations and conventions of \S1.  

Since morphisms in $\C D^b_{\Bbb C}$ are in general not chain maps, but rather equivalence classes of pairs of them, it is reasonable to expect $\Lambda$ to have this form, and this is indeed the case. Let

$$
\C N^k_{\tilde V}:=\left\{
\aligned &\C A^k(\lE),\,k=0,1\\
         &\C A^0(\varOmega^2)\oplus \C A^1(\C N(Z-E))\oplus \C A^2(\C O(-E)),\,k=2\\
         & 0,\,k=3,\,4
\endaligned\right.\tag4.1
$$
and define 
$$
\hat\C N^{\cdot}:=\pi_*\C N^{\cdot}_{\tilde V}
$$
We use the notation $\C N^k_{\tilde V}$ here because of the important role to be played by the Nash sheaf $\C N$, defined in \S3. 

There is an obvious decreasing filtration (by holomorphic degree) on $\hat\C N^{\cdot}$, and we thus regard it as an element of $\C D\C F^b_{\Bbb C}(V)$. It is easy to see that $\hat\C N^{\cdot}$ is a subcomplex of both $\hat\C L^{\cdot}$ and of $\hat\C A^{\cdot}$, respecting all the filtrations, so there are maps, $\lambda_1$ and $\lambda_2$ in  $\C D\C F^b_{\Bbb C}(V)$, 
$$
\hat\C L^{\cdot} @<\lambda_1<< \hat\C N^{\cdot} @>\lambda_2>> \hat\C A^{\cdot}\tag4.2
$$ 
What we will show is that
$$
\Lambda:=\lambda_1\lambda_2^{-1}:\hat\C L^{\cdot}_V @>>> \hat\C A^{\cdot}_V\tag4.3
$$
is an isomorphism in $\C D\C F^b_{\Bbb C}(V)$.
To do this it is necessary and sufficient (\cite{Il, V.1.2}) to show that for each $p$ the morphisms of the associated graded complexes
$$
gr^p_{\C F}\hat\C L^{\cdot}_V @< gr^p_{\C F}\lambda_1<< gr^p_{\C F}\hat\C N^{\cdot}_V @>gr^p_{\C F}\lambda_2>> gr^p_{\C F}\hat\C A^{\cdot}_V \tag4.4
$$
are isomorphisms in $\C D^b_{\Bbb C}(V)$; {\it i.e.}, that $gr^p_{\C F}\lambda_1$ and $gr^p_{\C F}\lambda_2$ are quasi-isomorphisms. 

Now in a fixed degree $k$, we have for $q:=k-p$
$$
gr^p_{\C F}\hat\C L^k_V\subset \hat\C L^{p,q}\tag4.5
$$
the $L^2$-forms of type $(p,q)$ on $V$. Here is an important observation:

\proclaim{\rm 4.6}{\bf Remark} This is {\rm not} an equality unless $p=\dim V$. In general, to show that a form $\omega\in L^{p,q}$ is in 
the image of  $gr^p_{\C F}\hat\C L^k_V$,  we must show either that $\partial_{B}\omega\in L_2$ (for appropriate $B$) or find a 
form $\alpha$ in $\C F^{p+1}\hat\C L^k_V$ so that 
$\da \omega + d\alpha\in L_2$. 
\endproclaim

Entirely analogous remarks apply to the filtered complexes of sheaves $\hat\C N^{\cdot}$ and $\hat\C A^{\cdot}$: the associated graded complex of each differs, in its imposition of an ``extra" $\partial$-condition, from its naturally associated ``pure $\bar\partial$"-complex, denoted $\hat\C N^{p,\cdot}$ and $\hat\C A^{p,\cdot}$; and there are inclusions of complexes of sheaves
$$
\hat\C L^{p,\cdot} @<\lambda_1^p<< \hat\C N^{p,\cdot} @>\lambda_2^p>> \hat\C A^{\cdot}\tag4.7
$$ 

Now it is cohomology of these latter complexes  $\hat\C L^{p,\cdot}$, $\hat\C N^{p,\cdot}$ and $\hat\C A^{p,\cdot}$ which can be computed most easily, so our strategy is to show that the morphisms which forget the $\partial$-conditions are quasi-isomorphisms and then to show that $\lambda_1^p$ and $\lambda_2^p$ are quasi-isomorphisms.

To be precise, we have the commutative diagram of complexes of sheaves
$$
\matrix\format\c&\quad\c&\quad\c\\
gr^p\hat{\Cal L}\spdot\gets & gr^p\hat{\Cal N}\spdot\to  & gr^p\hat{\Cal A}\spdot\\
\downarrow{\hat\kappa}^p_{\Cal L} &     \downarrow{\hat\kappa}^p_{\Cal N} &    \downarrow{\hat\kappa}^p_{\Cal A}\\
\hat{\Cal L}^{p,\,\cdot}\gets & \hat{\Cal N}^{p,\,\cdot}\to & \hat{\Cal A}^{p,\,\cdot}
\endmatrix\tag4.8
$$
Hence to prove the top horizontals in (4.8) are quasi-isomorphisms, it is enough to prove that the vertical and bottom horizontal morphisms are.

These are local statements which are obvious on the smooth part of $V$, where the complexes in (4.8) are identical. Hence Theorem B will follow from:
\medskip
\proclaim {\rm 4.9}{\bf Theorem} Let $V$ be a complex projective surface and let $\Cal S^{\cdot}=\hat\Cal L^{\cdot}$, $\hat\Cal N^{\cdot}$ or $\hat\Cal A^{\cdot}$.
\roster
\ite a Let $v\in V$ be a singular point. Then the stalk map
 $$
{\kappa}^p_{{\Cal S},v}: gr^p{\Cal S}^{\cdot}_v \to {\Cal S}^{p,\,\cdot}_v
$$
induces isomorphisms on cohomology for all $p$.  
\smallskip
\ite b The local cohomology groups $H^{p,q}({\Cal S}^{\cdot}_v):=H^q({\Cal S}^{p,\,\cdot}_v)$, arranged in the Hodge diamond,
$$
H^{2,2}({\Cal S}^{\cdot}_v)
$$
$$
 H^{2,1}({\Cal S}^{\cdot}_v) \qquad H^{1,2}({\Cal S}^{\cdot}_v)
$$
$$
\hss H^{2,0}({\Cal S}^{\cdot}_v)\qquad H^{1,1}({\Cal S}^{\cdot}_v)\qquad H^{0,2}({\Cal S}^{\cdot}_v)
$$
$$
H^{1,0}({\Cal S}^{\cdot}_v) \qquad H^{0,1}({\Cal S}^{\cdot}_v)
$$
$$
H^{0,0}({\Cal S}^{\cdot}_v)
$$
are isomorphic to the local cohomology groups\pagebreak 
$$
0
$$
$$
 0 \,\,\,\qquad\qquad\,\,\, 0
$$
$$
{\pi}_*{\Omega}^2_{\tilde V,v}\qquad H^1({\hat\Cal A}^{1,\cdot}_v)\qquad R^2{\pi}_*{\Cal O}_{\tilde V,v}
$$
$$
{\pi}_*{\Omega}^1_{\tilde V,v} \qquad R^1{\pi}_*{\Cal O}_{\tilde V,v}
$$
$$
{\pi}_*{\Cal O}_{\tilde V,v}
$$
where $\pi:{\tilde V}\to V$ is any resolution of the singularities of $V$.
Moreover, these isomorphisms are compatible with the maps on stalk cohomology induced by the bottom horizontal maps in (4.8).
\endroster
\endproclaim

\indent{\bf Proof}: The stalk cohomology groups are all direct limits of global section cohomology of open neighborhoods $U$ of $v$. In case $\Cal S^{\cdot}=\hat\Cal N\spdot$ or $\hat\Cal A\spdot$ these in turn are (since $\hat\Cal N\spdot\, \text{and}\, \hat\Cal A\spdot$ are direct image sheaves) cohomology over open neighborhoods ${\tilde U}:={\pi}^{-1}(U)$ of $E:={\pi}^{-1}(v)$, the exceptional divisor in ${\tilde V}$; in case ${\Cal S}\spdot=\hat\Cal L\spdot$, the same is true because of (3.9) and (3.10). So in the proof we will work in such ${\tilde U}$ without further comment, except when it is necessary to choose $U$ (and hence $\tilde U$) to have a psedudoconvex boundary. (This is permissible since such $U$ are cofinal among all neighborhoods of $v$ in $V$.) Also the compatibility of the isomorphisms in {\it b.} will be clear from the proofs and will be left to the reader. Finally, the isomorphisms in {\it a.} will all take the form 
$$
{\kappa}^{p,q}_{{\Cal S},v}: H^{p+q}(gr^p{\Cal S}^{\cdot}_v) \to H^{q}({\Cal S}^{p,\,\cdot}_v)
$$
and will be done case-by-case, identified by a choice of ${\Cal S}^{\cdot}$ and of $(p,q)$. The issue in these arguments will be the same as that described in the introduction to this section, in particular (4.6): the elements of $gr^p{\Cal S}^{\cdot}_v$ must be in $\dom d$, while those in ${\Cal S}^{p,\,\cdot}_v$ satisfy the weaker $\dom \bar\partial$ condition.

\indent$\bullet$ ${\Cal S}\spdot=\hat\Cal L\spdot \,\,\text{and}\,\, (p,q)=(0,0)$: Let $[\phi ]\in H^{0}( U;{\Cal L}^{0,\,\cdot}_v)$. Then $[\phi ]$ has a representative $\phi\in {\Cal L}^{0,0}(U)$ with $\partial_D \phi=0$. This follows from (2.27), but can also be proved more 
readily in this case as follows. Let ${\phi}_i\to \phi$ and $\bar{\partial}{\phi}_i\to \bar{\partial}_D\phi=0$, where ${\phi}_i$ is a sequence of smooth functions on $ U$ supported away from $v$. Let $\eta$ be a smooth compactly supported function on $U$, $\eta\equiv 1$ on a neighborhood $U'$ of $v$. Clearly,
$$
\bar\partial(\eta{\phi}_i)=\bar\partial\eta\land {\phi}_i+\eta\bar\partial{\phi}_i\to \bar\partial\eta\land \phi
$$
In particular, lim $\bar\partial(\eta{\phi}_i)$ exists, and so is Cauchy. Now using equality of the Laplacians $\Delta_{\partial}$ and $\Delta_{\bar\partial}$ on functions of compact support, we get
$$
<\partial(\eta{\phi}_i-\eta{\phi}_j),\partial(\eta{\phi}_i-\eta{\phi}_j)>\,\,=\,\,<\bar\partial(\eta{\phi}_i-\eta{\phi}_j),\bar\partial(\eta{\phi}_i-\eta{\phi}_j)>
$$
for all $i$ and $j$, so that $\partial(\eta{\phi}_i)$ is Cauchy as well. Thus we have the convergence of $\partial(\eta{\phi}_i)$ on $U'$ so that $\phi|U'\in gr^0{\hat\Cal L}^{0.\cdot}(U')$. Thus ${\kappa}^{0.0}_{{\hat\Cal L},v}: H^{0}(gr^0{\hat\Cal L}^{\cdot}_v) \to H^{0}({\Cal L}^{0,\,\cdot}_v)$ is surjective. As it is clearly injective, it is an isomorphism. 

Finally, by [P, 4.7], the natural map ${\Cal O}(Z-E)(\tilde U)\to {\Cal L}^{0,0}_N(\tilde U)$ induces an isomorphism onto $\ker\bar\partial_N$. This implies that $\ker\db_D=\C O(\tilde U)$; for if $\db_D\phi=0$, then 
by (2.27), $\phi/\rho\log(1/\rho)\in\C L^{0,0}(\tilde U)$, where $\rho=r\circ\pi$ and $r$ is the 
distance from the singular point $v$. Since we already have $\phi\in\C O(Z-E)(\tilde U)$, it follows from (3.4b) and (3.3b) that $\phi\in\C O(\tilde U)$. Conversely, if $\phi\in\C O(\tilde U)$, then clearly $\db_N\phi=0$; and $\phi/\rho\in\C L^{0,0}$ by (3.4b), since $Z\leq D_\gamma$. 
Thus by (2.18), $\phi\in\dom\bar\partial^0_D$. This completes the proof of part {\it b.} in this case.
\smallskip
\indent$\bullet$ ${\Cal S}\spdot=\hat\Cal A\spdot\,\, \text{and}\,\, (p,q)=(0,0)$: By definition, $\hat\C A^{0,0}(U)=\C A^{0,0}(\tilde U;\log E)=\C A^{0,0}(\tilde U)$. 
Since $\da\C A^{0,0}(\tilde U)\subseteq \C A^{1,0}(\tilde U)\subseteq \C A^0(\tilde U,\varOmega^1(\lE))=\hat\Cal A^{1,0}(U)$, ${\kappa}^{0,0}_{{\Cal A},v}$ is surjective. Injectivity is obvious.
\smallskip
\indent$\bullet$ ${\Cal S}\spdot=\hat\Cal N\spdot\,\, \text{and}\,\, (p,q)=(0,0)$: In degree $\leq 1$, $\hat\Cal N\spdot=\hat\Cal A\spdot$, so this case is identical to the previous one.
\smallskip
\indent$\bullet$ ${\Cal S}\spdot=\hat\Cal L\spdot\,\, \text{and} \,\,(p,q)=(0,1)$: To begin, observe that we have natural inclusions of complexes of sheaves
$$
{\pi}_*{\Cal A}^{0,\cdot}_{\tilde V}\to \hat{\Cal L}^{0,\cdot}
$$
and
$$
\hat{\Cal L}^{2,\cdot}\to {\pi}_*{\Cal A}^{2,\cdot}_{\tilde V}
$$
These induce the horizontal maps in the commutative diagram
$$
\CD
H^{0,1}(\tilde U)@>>> H^1(\tilde U; {\Cal L_D}^{0,\cdot})\\
@V{\cong}VV     @V{\cong}VV\\
H^{2,1}_c(\tilde U)^*@>\cong>> H^1_c(\tilde U; {\Cal L_N}^{2,\cdot})^*
\endCD
$$
in which we take $\tilde U$ to be a pseudoconvex neighborhood of $E$, so that $H^{0,1}(\tilde U)$ and $H^{2,1}_c(\tilde U)$ are finite -dimensional (\cite{KF, 4.3.2, 5.1.7}). The bottom horizontal map is an isomorphism by the main theorem of [PS], and the verticals are Serre Duality isomorphisms ([KF, 5.1.7]and [PS, 1.3c]). This proves part {\it b.} in this case; and it shows that, given $[\phi]\in H^1(\tilde U, {\Cal L_D}^{0,\cdot})$, we may assume  $\phi$ is a smooth 1-form on $\tilde U$. We now claim that $[\partial\phi]\in H^1(\tilde U;{\varOmega}^1)$ is in the image of the map
$$
H^1(\tilde U;{\Cal N}(Z-E))\to H^1(\tilde U;{\varOmega}^1)
$$
induced by the inclusion ((3.19) or (3.20)) of sheaves ${\Cal N}(Z-E)\hookrightarrow {\varOmega}^1$. If this is assumed, there is $\psi\in A(\tilde U;\C N(Z-E))\subseteq A^{1,0}(\tilde U)$ such that 
$$
\bar\partial\psi-\partial\phi\in \Cal A^1(\tilde U;{\Cal N}(Z-E))
$$
Since ${\Cal N}(Z-E)\subseteq {\Cal L}^{1,1}$ and $\partial\psi\in \Cal A^{2,0}(\tilde U)\subseteq \Cal L^{2,0}(\tilde U)$ we have $d(\psi+\phi)\in \Cal L^1(\tilde U)$, i.e., $\psi+\phi\in \dom\,d_N$. But since $\psi+\phi$ is smooth it follows from (2.18) that $\psi+\phi\in \dom\,d_D$, so ${\kappa}^{0,1}_{{\Cal L},v}$ is surjective. (Compare this argument to Remark (4.6).)

 To prove the above claim, recall the exact sequence of sheaves (3.27)
$$
0\to \Cal I_E\varOmega^1(\lE)\to \varOmega^1\to \oplus\varOmega^1_{E_i} \to 0
$$
This induces the exact sequence of vector spaces (3.28)
$$
H^1(\tilde U;\Cal I_E\varOmega^1(\lE)) @>f>> H^1(\tilde U;\varOmega^1) @>g>> \oplus H^1(E_i;\varOmega^1_{E_i})
$$
Clearly, $g[\partial\phi]=[\partial g\phi]=0$, so $[\partial\phi]\in \operatorname{im}f$. Hence, (3.24a) finishes the proof of the claim.

\indent To show injectivity, suppose $\omega:=\psi+\phi\in \Cal L^{1,0}(\tilde U)\oplus \Cal L^{0,1}(\tilde U)=\Cal L^1(\tilde U)$ is such that $d_D\omega\in F^1\Cal L^2(\tilde U)$ and $\phi=\bar\partial{_D}f$ for some $f\in \Cal L^0(\tilde U)$. By the argument from the case ${\Cal S}\spdot=\hat\Cal L\spdot$ and $(p,q)=(0,0)$ above, $f\in \dom \partial_D$, so in $H^1(gr^0\hat\C L^{\cdot}_v)$, $[\omega]=[\psi]=0$.
\smallskip
\indent$\bullet$ ${\Cal S}\spdot=\hat\Cal A\spdot\,\, \text{and}\,\, (p,q)=(0,1)$: By definition, $H^{0,1}(\tilde U)@>\cong>> H^1(\tilde U; {\Cal A}^{0,\cdot})$. Next, given $\phi\in \hat{\Cal A}^{0,1}(U)=\Cal A^{0,1}(\tilde U)$ with $\bar\partial\phi=0$, the argument in the previous case shows there is $\psi\in A^{1,0}(\tilde U)\subseteq \Cal A^0(\tilde U;{\varOmega}^1(\lE)):= \hat{\Cal A }^{1,0}(U)$ such that $\bar\partial\psi-\partial\phi\in A^1(\tilde U;\Cal I_E{\varOmega}^1(\lE))=\hat{\Cal A }^{1,1}(U)$. Since $\partial\psi\i \Cal A^0(\tilde U;{\varOmega}^2)=\hat{\Cal A }^{2,0}(U)$, this shows ${\kappa}^{0,1}_{{\Cal A},v}$ is surjective. Injectivity is obvious, so ${\kappa}^{0,1}_{\hat\Cal A,v}$ is an isomorphism. 
\smallskip
\indent$\bullet$ ${\Cal S}\spdot=\hat\Cal N\spdot \,\,\text{and} \,\,(p,q)=(0,1)$: Since
 $\hat{\Cal N}^{1,1}(U):=\Cal A^1(\tilde U;{\Cal N}(Z-E))$ and 
 $\hat{\Cal N}^{2,0}= {\Cal A}^{2,0}$, the surjectivity follows from the 
proof of the case ${\Cal S}\spdot=\hat\Cal L\spdot$, except that we may 
immediately assume that our form $\phi$ is in ${\Cal A}^{0,1}$, 
since $\hat{\Cal N}^{0,q}={\Cal A}^{0,q}$. The injectivity follows similarly and $H^{0,1}(\tilde U)@>>> H^1(U, \hat{\Cal N}^{0,\cdot})$ is an isomorphism by the previous case, because $\hat{\Cal N}^{0,\cdot}:=\hat{\Cal A}^{0,\cdot}$.
\smallskip
\indent$\bullet$ ${\Cal S}\spdot=\hat\Cal L\spdot \,\,\text{and}\,\, (p,q)=(0,2)$: First note that $H^{2}({\hat\Cal L}^{0,\,\cdot}_v)=0$ since, for any pseudoconvex neighborhood $\tilde U$ of $E$, \cite{PS,KF. {\it loc.cit.}} shows $H^{0,2}_D(\tilde U)\cong \varOmega^2_c(\tilde U)^*=0$. It remains to show that 
$H^2(gr^0\hat\C L^\cdot_v)=0$

Let $[\phi]\in H^2(U,gr^0\hat\C L^\cdot)$; then $\phi\in\hat\C L^2(U)$ and $d_N\phi\in\hat\C L^3(U)$. Suppose $\kappa^{0,2}_{\hat\C L}[\phi]=0$; then there exists $\xi\in \hat\C L^{0,1}(U)$ such that $\db_D\xi=\phi^{0,2}$. By (2.47), we can choose  this $\xi$ so that $\da_D\xi\in\hat\C L_2$. Thus, $[\psi]=0\in H^2(V;gr^0\hat\C L^\cdot)$. 
\smallskip
\indent$\bullet$ ${\Cal S}\spdot=\hat\Cal A\spdot \,\,\text{and}\,\, (p,q)=(0,2)$: By Malgrange's theorem [M], $H^2(\tilde U;\Cal I_E)=0$, so $H^{2}({\hat\Cal A}^{0,\,\cdot}_v)=0$. We will show that 
$$
H^2(\hat\Cal A^{\cdot}_v)\to  H^2(gr^0\hat\Cal A^{\cdot}_v)\tag4.11
$$
is surjective. Then since $\hat\C A^\cdot$ is quasi-isomorphic to $\C I\C C^\cdot$, the intersection cohomology complex on $V$ (see (1.  )), we have $H^2(\hat\C A^\cdot_v)=0$ so  that $H^2(gr^0\hat\Cal A^{\cdot}_v)=0$ as well.

\indent So let $\xi=\xi^{2,0}+\xi^{1,1}+\xi^{0,2}\in \hat\Cal A^2(U)$ satisfy $d\xi\in F^1\hat\Cal A^3(U)=\hat\Cal A^3(U)$. By Malgrange again, $H^2(\tilde U;\hat\C A^{1,\cdot})=H^2(\tilde U;\Cal I_E\varOmega^1(\lE))=0$, so there is $\eta\in \hat\Cal A^{1,1}(U)=\Cal A^1(\tilde U;\Cal I_E\varOmega^1(\lE))$ such that $\bar\partial\eta=\bar\partial\xi^{1,1}-\partial\xi^{0,2}$. Since the elements of $\Cal A^1(\tilde U;\Cal I_E\varOmega^1(\lE))$ are smooth ($\C I_E\varOmega^1(\lE)\subseteq\varOmega^1$), $\partial\eta\in \Cal A^1(\tilde U;\varOmega^2)=\hat\Cal A^{2,1}(U)$, so $\eta\in F^1\hat\Cal A^2(U)$. Hence, replacing $\xi$ with $\xi-\eta$, we may assume that $(d\xi)^{1,2}=0$. But now we have $\bar\partial(d\xi)^{2,1}=0$, so by [GR] there is $\tau\in \Cal A^0(\tilde U;\varOmega^2)=\Cal A^0(\tilde U;\Cal I_E\varOmega^2(\lE))=\hat\Cal A^{1,1}(U)$ such that $\bar\partial\tau=(d\xi)^{2,1}$, As above we may assume that $d\xi=0$. Hence (4.11) is surjective. 
\smallskip
\indent$\bullet$ ${\Cal S}\spdot=\hat\Cal N\spdot\,\, \text{and}\,\, (p,q)=(0,2)$: By definition, $H^{2}({\hat\Cal N}^{0,\,\cdot}_v)=H^{2}({\hat\Cal A}^{0,\,\cdot}_v)$, which was just shown to vanish. Again it follows easily from the definition of $\hat\Cal N\spdot$ that $H^2(\hat\Cal N^{\cdot}_v)\hookrightarrow H^2(\hat\Cal A^{\cdot}_v)$, which vanishes, as we saw above. So to show that $H^2(gr^0\hat\Cal N^{\cdot}_v)=0$ it suffices to show that
$$
H^2(\hat\Cal N^{\cdot}_v)\to  H^2(gr^0\hat\Cal N^{\cdot}_v)
$$
is surjective. But this is immediate from the fact that $\hat\Cal N^3=0$ (so $\xi\in \hat\C N^2(U)$ representing an element of $H^2(gr^0\hat\Cal N^{\cdot}_v)$ is automatically closed).
\smallskip
\indent$\bullet$ ${\Cal S}\spdot=\hat\Cal L\spdot \,\,\text{and}\,\, (p,q)=(1,0)$: By definition, $\ker \bar\partial_N=\varOmega^1_{(2)}(\tilde U)$ and so by (2.43) and the definition (3.11) of $\C N$, $\ker \bar\partial_D=\C N(N)(\tilde U)$. But by (3.24b) $\C N(N)(\tilde U)=\varOmega^1(\lE)(\tilde U)$, which equals $\varOmega^1(\tilde U)$ by (3.29a).  Hence, $\ker \bar\partial_D= \varOmega^1(\tilde U)$ and it follows easily from this that
${\kappa}^{0,1}_{{\Cal L},v}$ is an isomorphism: surjectivity is the only issue, and we need only show that if $\omega\in\dom \bar\partial^{1,0}_D$, then $\omega\in\dom \partial^{1,0}_D$. For this, the proof above that $\phi\in \dom\,\partial_D$, in case ${\Cal S}\spdot=\hat\Cal L\spdot \,\,\text{and}\,\, (p,q)=(0,0)$, applies essentially verbatim. The only thing to add is that the identity used there becomes
$$
\align
<\bar\partial(\eta{\phi}_i-\eta{\phi}_j),\bar\partial(\eta{\phi}_i-\eta{\phi}_j)>&=<\partial(\eta{\phi}_i-\eta{\phi}_j),\partial(\eta{\phi}_i-\eta{\phi}_j)>+<\bar\vartheta(\eta{\phi}_i-\eta{\phi}_j),\bar\vartheta(\eta{\phi}_i-\eta{\phi}_j)>\\&\geq <\partial(\eta{\phi}_i-\eta{\phi}_j),\partial(\eta{\phi}_i-\eta{\phi}_j)>
\endalign
$$
where $\bar\vartheta$ denotes as usual the formal adjoint of $\partial$.
\smallskip
\indent$\bullet$ ${\Cal S}\spdot=\hat\Cal A\spdot \,\,\text{and}\,\, (p,q)=(1,0)$: By definition, $H^0(U;\hat{\Cal A}^{1,\cdot})={\varOmega}^1(\lE)(\tilde U)$ which equals ${\varOmega}^1(\tilde U)$ by (ref). It follows from this that ${\kappa}^{0,1}_{{\Cal A},v}$ is an isomorphism.
\smallskip
\indent$\bullet$ ${\Cal S}\spdot=\hat\Cal N\spdot\,\, \text{and}\,\, (p,q)=(1,0)$: The proof here is identical to that of the previous case.
\smallskip
\indent$\bullet$ ${\Cal S}\spdot=\hat\Cal L\spdot \,\,\text{and}\,\, (p,q)=(1,1)$: Let $[\xi]\in H^2(U;gr^1\hat\Cal L\spdot)$ be in the kernel of ${\kappa}^{1,1}_{{\hat\Cal L},v}$. Then given the decomposition $\xi=\xi^{2,0}+\xi^{1,1}$ of $\xi$ into type, we have $d_N\xi\in F^2\hat\Cal L^3(U)$ and there is $\phi\in \hat\Cal L^{1,0}(U)$ such that $\bar\partial_D\phi=\xi^{1,1}$. We showed above in case ${\Cal S}\spdot=\hat\Cal L\spdot \,\,\text{and}\,\, (p,q)=(1,0)$ that this implies $\phi\in\dom d_D$, so we get $[\xi]=[\xi-d_D\phi]=0$ in $H^2(U;gr^1\hat\Cal L\spdot)$.

\indent To prove surjectivity of ${\kappa}^{1,1}_{{\hat\Cal L},v}$, suppose given $[\xi^{1,1}]\in H^1(\tilde U;\Cal L^{1,\cdot})$. We show below that we may take $\xi^{1,1}$ to be a $\C N(Z-E)$-valued (0,1)-form. (In fact, we show there is a surjection $H^1(\pi_*\C N(Z-E)_v)\to H^1(\hat\C L^{1,1}_v)$). If we assume this, then since $\C N(Z-E)\subseteq \varOmega^1$, $\partial\xi^{1,1}\in \Cal A^{2,1}(\tilde U)$, the smooth (2,1)-forms on $\tilde U$. Since $\bar\partial\partial\xi^{1,1}=-\partial\bar\partial\xi^{1,1}=0$, there is  ([GR]) $\xi^{2,0}\in \Cal A^{2,0}(\tilde U)$ such that $\bar\partial\xi^{2,0}=-\partial\xi^{1,1}$. Hence $\xi:= \xi^{2,0}+\xi^{1,1}\in F^1\Cal L\spdot(\tilde U)$, $d_N\xi\in F^2\Cal L\spdot(\tilde U)$ and ${\kappa}^{1,1}_{{\hat\Cal L},v}[\xi]=[\xi^{1,1}]$ as required.

\ind To justify the assumption just made,we show there is an isomorphism
$$
\im (H^1(\pi_*\C N(Z-E)_v)\to H^1(\pi_*\C N(N)_v))\cong H^1(\hat\C L^{1,\cdot}_v),\tag4.12
$$
which will also prove the isomorphism of part {\it b.} of the Theorem by (3.26a). 

To begin we show there is a  commutative diagram of sheaves on $\tilde V$
$$
\CD
\C M^0(\C N(Z-E)) @>\bar\partial^0>> \C M^1(\C N(Z-E)) @>\bar\partial^1>>
 \C M^2(\C N(Z-E)) \\
@VVV        @VVV         @AAiA\\
\C L^{1,0}_\gamma @>\bar \partial^0_D>> \C L^{1,1}_\gamma @>\bar \partial^1_N>> \C L^{1,2}_\gamma 
\endCD\tag4.13
$$
where $\C M\spdot(\C N(Z-E))$ denotes the $\bar\partial$-complex of sheaves of $\C N(Z-E))$-valued measurable forms of type $(0,\cdot)$ on $\tilde V$, $\C L^{p,q}_\gamma$ is the sheaf of measurable forms on $\tilde V$ which have locally finite $L_2$ norm with respect to the (degenerate) metric $\gamma$ pulled up from the induced Fubini-Study metric on $V$, and the vertical maps are inclusions. First of all, the two rightmost inclusions follow from (3.12d) and the middle one shows that $\C M^0(\C N(Z-E))\subseteq\dom\db^0_N$. To complete the justification of the diagram, we must show that $ \C M^0(\C N(Z-E))\subseteq \dom\bar\partial^0_D$. By (3.12d), $\C L^{1,0}_\gamma=\C M^0(\C N(D_{\gamma})):=\C M^0(\C N(Z-E+N))$. 
Let $\rho:\tilde U\to \Bbb R$ denote the composition of $\pi:\tilde U\to U$ with the distance map $r:\tilde U\to \Bbb R$ coming from the imbedding of $U$ into $\Bbb C^N$. Since the divisor $N\geq Z$ (by (3.3c)), we see from (3.4b) that if $\omega\in\C M^0(\C N(Z-E)))$, then $\omega/\rho\in \C L^{1,0}_\gamma$. We now conclude that $\omega\in\dom\db^0_D$ as in the case ${\Cal S}\spdot=\hat\Cal L\spdot \,\,\text{and}\,\, (p,q)=(0,0)$.

 Let us now apply $\Gamma(\tilde V;\frac{}{})$ to (4.13). Observe that the composite operator $i\bar\partial_N$ is compactly approximable in norm since $E$ is of real codimension two in $\tilde U$ (\cite{PS,, proof of (3.6)}.). This allows us to replace the diagram above with the commutative diagram of Hilbert spaces
$$
\CD
 M^0(\C N(Z-E)) @>\bar\partial^0>>  M^1(\C N(Z-E)) @>\bar\partial^1>>
  M^2(\C N(Z-E)) \\
@VVV        @VVV         @VV=V\\
 L^{1,0}_\gamma @>\bar \partial^0_D>>  L^{1,1}_\gamma  @>\bar \partial^1_D>>  M^2(\C N(Z-E)) 
\endCD\tag4.14
$$

To use (4.13), we introduce the following Lemma. It gives conditions on a map between complexes of Hilbert spaces under which one may conclude surjectivity of the induced map on cohomology; in effect, it gives conditions under which one may reverse the standard implication (\cite{PS,1.3(a)}) ``$H^k(V,\Cal L^{\cdot})$ finite dimensional $\Rightarrow$ range $d$ closed". 
\proclaim {\rm 4.15}{\bf Lemma} Let $i^{\cdot}:(M^{\cdot},D^{\cdot})\hookrightarrow (L^{\cdot},d^{\cdot})$ be a bounded inclusion of complexes of Hilbert spaces, $\cdot=0,1,2$, such that
\roster
\item the operators $D^{\cdot}$ and $d^{\cdot}$ are closed,
\item the  cohomology $H^1(M^{\cdot})$ is finite dimensional,
\item $i^2$ is an equality,
\item there is a subspace $L^{1}_c\subseteq M^{1}$ such that if $d^1_c:=d^1|L^1$, then the operator closure $\Bar{d^{1}_c}=d^1$ and
\item range $d^0$ and range $D^1$ are closed.
\endroster
Then if $\phi\in \dom d^1$, there exist $\psi\in \dom D^1$ and $\lambda\in L^0$ such that $\phi-d^0\lambda=\psi$. In particular, $H^1(M^{\cdot})\to H^1(L^{\cdot})$ is surjective.
\endproclaim
\demo
\ind {\bf Proof:} To aid in following the proof, we display $i^{\cdot}$:
$$
\CD
M^0 @>D^0>>  M^1 @>D^1>> M^2\\
@Vi^0VV      @Vi^1VV    @Vi^2V=V\\
L^0 @>d^0>> L^1 @>d^1>> L^2
\endCD\tag4.16
$$
An element of $M^i$ will be regarded when convenient as an element of $L^i$. 

\ind Let $\phi\in \dom d^1$. Assumption (4) says there is a sequence $\{\phi_j\}$ in $L^1_c$ such that $\phi_j @>L^1>> \phi$ and $\{d^1\phi_j=D^1\phi_j\}$ converges in $M^2=L^2$. (Here and below, $\xi_j @>H>> \xi$ means that the sequence $\xi_j$ converges in the Hilbert space $H$ to $\xi$.) By assumptions (1) and (2), $M^1$ has a ``Hodge Decomposition" (\cite{KK, Appendix})
$$
M^1=\im D^0\perp \im D^{1*}\perp \C H^1
$$
where $\C H^1:=\ker D^1\cap \ker D^{0*}$ is finite-dimensional. Using this, write 
$$
\phi_j=D^0\beta_j+D^{1*}\gamma_j+h_j
$$
Note that we do not know whether $\{\phi_j\}$ converges in $M^1$, so we can't conclude that any of these tems converge there. Since $D^1\phi_j=D^1D^{1*}\gamma_j$ converges in $M^2$ and the range of $D^1$ is closed, $D^{1*}\gamma_j$ converges in $M^1$. We replace $\phi$ with $\phi-\lim D^{1*}\gamma$ and denote it $\phi$ again. We now have $d^1\phi=0$ and a sequence $\{\phi_j=D^0\beta_j+h_j\}$ in $M^1$ (not necessarily in $L^1_c$) such that $\phi_j @>L^1>> \phi$ and $D^1\phi_j=d^1\phi_j @>M^2=L^2>> 0$.

\ind Now write
$$
h_j=h^0_j+h^1_j\quad\text{where}\quad h^0_j\in \ker (\C H^1\to H^1(L^{\cdot})),\,\,\text{and}\,\,h^1_j\in \ker (\C H^1\to H^1(L^{\cdot}))^{\perp},
$$
and, for each $j$, choose $\alpha_j\in L^0$ such that $d^0\alpha_j=D^0\beta_j+h^0_j$. 

Since $h^1_j/\|h^1_j\|_{L^1}$ is bounded in the finite-dimensional subspace $i^1\C H^1\subseteq L^1$, we may assume, perhaps after passing to a subsequence, that it converges, say to $h^1\in L^1$. We claim the sequence $\{\|h^1_j\|_{L^1}\}$ is bounded. If it were unbounded, then (passing again to a subsequence if necessary) $\phi/\|h^1_j\|_{L^1} @>L^1>> 0$, so $d^0(\alpha_j/\|h^1_j\|_{L^1}) @>L^1>> -h^1$. But then, since the range of $d^0$ is closed, $h^1\in \im d^0$, which contradicts the fact that $h^1\in\ker (\C H^1\to H^1(L^{\cdot}))^{\perp}$. 

\ind Now since $\|h^1_j\|_{L^1}$ is bounded, $\|h^1_j\|_{M^1}$ is bounded too, because the $L^1$- and $M^1$-norms are equivalent on $\C H^1$. So there is a convergent subsequence, $h^1_j @>M^1>> h^*$, and from 
$$
\phi_j=D^0\beta_j+h^0_j+h^1_j=d^0\alpha_j+h^1_j
$$
we get $\phi-h^*=\lim d^0\alpha_j$, which equals $d^0\lambda$, for some $\lambda\in L^0$, since $d^0$ is closed.
\enddemo

We have verified that all the hypotheses of (4.15) are satisfied in (4.14) except one, namely that the range of $\bar \partial^0_D$ be closed, which we now verify. First the range of $\bar \partial^0_D$ is closed if and only if that of its Hilbert space adjoint $(\bar \partial^0_D)^*$ is; and this would follow from the $\db_D$-Hodge decomposition if we knew $H^{1,0}_D(V)$ were finite. Because $H^{1,0}_D(V)\subseteq H^{1,0}_N(V)$, this follows from the computation made in (3.12d):

\proclaim{\rm 4.17}{\bf Lemma} $H^{1,0}_N(V)=H^0(\tilde V;\C N(D_\gamma))$
\endproclaim

It now follows from (4.15) that $H^1(\C N(Z-E)_v)\to H^1(\C L^{1,\cdot}_v)$ is surjective:
take $\phi\in\C L^{1,1}(\tilde U)$for some $\tilde U\supset E$, get a global form $\eta\phi\in L_\gamma$ using a cut-off $\eta$ as in the case ${\Cal S}\spdot=\hat\Cal L\spdot \,\,\text{and}\,\, (p,q)=(0,0)$ above and use (4.15) to get $\psi\in M^1(\C N(Z-E))$ and $\lambda\in L^{1,0}_\gamma$ such that $\phi-\db_D\lambda=\psi$. Then on some $\tilde U'\subseteq \tilde U$, $[\phi]$ is in the image of $H^1(\tilde U';\C N(Z-E))\to H^1(\tilde U';\C L^{1,\cdot}_v)$. 

Now we claim that the map $H^1(\tilde U;\C N(Z-E)) \to H^1(\tilde U;\C L^{1,\cdot})$ passes to a map
$$
\im (H^1(\tilde U;\C N(Z-E))\to H^1(\tilde U;\C N(N))) @>>> H^1(\tilde U;\C L^{1,\cdot}_\gamma).
$$
To prove this, suppose that $\db\psi=\phi\in\C M^1(\tilde U;\C N(Z-E))$, where $\psi\in\C M^0(\tilde U;\C N(N))$. Then in $\C L^{1,\cdot}(\tilde U)$, $\db_N\psi=\phi$, and we claim that we can replace $N$ by $D$. To do this we appeal to the argument in \cite{PS, (3.6)}. Namely, we need the 
"trace estimate'' \cite{loc. cit.,(3.7)}, which follows since the equation $\db\psi=\phi$ holds in $\C M^\cdot(\C N(D_\gamma))$, where we can appeal to the same elliptic regularity (\cite{H2, (4.2.3)}) as was used in the proof of \cite{PS, (3.6)}. Finally, we claim (4.12) is an isomorphism. To see this, suppose $\phi\in\C M^1(\C N(Z-E))$ and $\db_D\psi=\phi$, where $\psi\in\C L^{1,0}(\tilde U)=\C M^0(\tilde U;\C N(D_\gamma))$. Then by (2.47), we can arrange that $\psi/\rho\log\rho\in\C L^{1,0}(\tilde U)$. It now follows easily from (3.4b) that, since $N=D_\gamma-Z+E$, $\psi\in\C M^0(\tilde U;\C N(N))$. Hence, (4.12) is injective.
\smallskip 
{\rm 4.18} {\bf Remark} The argument in the last paragraph can be used together with (4.17) to show ($D_\gamma-(Z-E)=N$)
$$
H^{1,0}_D(V)=H^0(V;\C N(N))
$$
\noindent Details are left to the reader.
\smallskip
\indent \indent$\bullet$ ${\Cal S}\spdot=\hat\Cal A\spdot \,\,\text{and}\,\, (p,q)=(1,1)$: Let $\xi=\xi^{2,0}+\xi^{1,1}\in \Cal A^0(\tilde U;\varOmega^2)\oplus \Cal A^1(\tilde U;\Cal I_E\varOmega^1(\lE))$ represent an element of the kernel of ${\kappa}^{1,1}_{{\hat\Cal A},v}$: hence $\xi^{1,1}=\bar\partial\phi$ where $\phi\in \Cal A^0(\tilde U;\varOmega^1(\lE))$. This implies that $[\xi]$ is in the kernel of the natural map $H^1(\tilde U;\varOmega^1(\lE))\to H^1(\tilde U;\Cal I_E\varOmega^1(\lE))$. This map factors through $H^1(\tilde U;\varOmega^1)$ and(3.29b) shows $[\xi]$ vanishes there, so there is $\phi\in \Cal A^0(\tilde U;\varOmega^1)$ such that $\bar\partial\phi=\xi^{1,1}$. Since $\hat\C A^{2,0}(U):=\C A(\tilde U;\C I_E\varOmega^2(log\,E))=\C A(\tilde U;\varOmega^2)$, $\partial\phi\in \hat\C A^{2,0}(U)$, so we have $[\xi]=[\xi-\xi^{1,1}]=[\xi^{2,0}]=0$ in $H^2(gr^1\hat\Cal A_v)$. So ${\kappa}^{1,1}_{{\hat\Cal A},v}$ is injective. The proof of surjectivity is essentially the same as that in the immediately previous case.
\smallskip
\indent$\bullet$ ${\Cal S}\spdot=\hat\Cal N\spdot \,\,\text{and}\,\, (p,q)=(1,1)$: That ${\kappa}^{1,1}_{{\hat\Cal N},v}$ is injective is proved in the same way as for ${\kappa}^{1,1}_{{\hat\Cal A},v}$; we only need to add that $N^*(Z-E)\subseteq \Cal I_E\varOmega^1(\lE)$. Surjectivity again uses the same argument as for that in case $(p,q)=(1,1)$ and $\Cal S\spdot=\hat\Cal L\spdot$.

To prove the isomorphism of part {\it b.}, notice first that 
$$
H^1(\hat\C N^{1,\cdot}_v)=\im \left( H^1(\tilde U;\C N(Z-E))\to H^1(\tilde U;\varOmega^1(\lE))\right)
$$
which, using the surjection in (3.24a), is isomorphic to
$$
\im \left( H^1(\tilde U;\C I_E\varOmega^1(\lE)))\to H^1(\tilde U;\varOmega^1(\lE))\right)
$$
as desired.
\smallskip
\indent$\bullet$ ${\Cal S}\spdot=\hat\Cal L\spdot \,\,\text{and}\,\, p+q>2$: By [O1], $H^q(\hat\Cal L^{p,\cdot}_v)=0$. To show $H^{p+q}(gr^p\hat\Cal L\spdot_v)=0$, we must strengthen Ohsawa's argument to include $\partial$-control. Specifically, we need the following result.

\proclaim{\bf Lemma} Let $\phi\in \hat\Cal L^{p,q}(U)$, $p+q>2$. Then there exists $\nu\in \hat\Cal L^{p,q-1}(U)$ such that $\bar\partial_N\nu=\phi$ and $\nu\in \dom\partial_N$.
\endproclaim
The proof is easily adapted from the careful exposition of Ohsawa's argument in \cite{PS, 2.3}
\smallskip
Now we can argue as we have several times above: given any $[\psi]\in H^{p+q}(U;gr^p\hat\Cal L\spdot)$, we may assume the (1,1)-component of $\psi$ is zero, which means $[\psi]=0$.
\smallskip
\indent$\bullet$ ${\Cal S}\spdot=\hat\Cal L\spdot,\,\hat\Cal A\spdot\,\text{or}\,\,\hat\Cal N\spdot \,\,\text{and}\,\, p=2$: In these cases, $gr^2{\Cal S}^{\cdot}={\Cal S}^{2,\,\cdot}$. The calculation of $H^*({\Cal S}^{2,\,\cdot}_v)$ is clear in case $\Cal S\spdot=\hat\Cal A\spdot\,\text{or}\,\,\hat\Cal N\spdot$ and follows from the main theorem of [PS] in case $\Cal S\spdot=\hat\Cal L\spdot$. 
\bigskip
We now compute the Neumann and Dirichlet $L_2-\db$-cohomology groups $H^{p,q}_N(V)$ and $H^{p,q}_D(V)$ of an algebraic surface $V$. Most of this has already been done: from \cite{PS} we have 
$$
H^{2,q}_N(V)\cong H^q(\tilde V;\varOmega^2),\qquad  H^{2,q}_D(V)\cong H^q(\tilde V;\varOmega^2(E-Z))\tag4.19
$$ 
and by duality
$$
H^{0,q}_N(V)\cong H^q(\tilde V;\C O(Z-E)),\qquad  H^{0,q}_D(V)\cong H^q(\tilde V;\C O);\tag4.20
$$
and from part {\it b.} of Theorem (4.8) and (3.29a) we have
$$
H^{1,2}_N(V)\cong H^2(\tilde V;\C I_E\varOmega^1(\log E)),\qquad  H^{1,0}_D(V)\cong H^0(\tilde V;\varOmega^1(\log E))\cong H^0(\tilde V;\varOmega^1)\tag4.21
$$

\proclaim {\rm 4.22}{\bf Theorem} Let $V$ be a complex projective surface. Then
$$
H^{1,2}_D(V)\cong H^2(\tilde V;\C I_E\varOmega^1(\log E)\otimes\C O(E-Z),\qquad  H^{1,0}_N(V)\cong H^0(\tilde V;\varOmega^1(\log E)\otimes\C O(Z-E)),
$$ 
$$
H^{1,1}_D(V)\cong\im (H^1(\tilde U;\C I_E\varOmega^1(\log E)\otimes\C O(E-Z))\to H^1(\tilde U;\varOmega^1(\log E)))
$$
and
$$  
H^{1,1}_N(V)\cong \im (H^1(\tilde U;\C I_E\varOmega^1(\log E))\to H^1(\tilde U;\varOmega^1(\log E)\otimes\C O(Z-E)))
$$
\endproclaim

{\bf Proof:} We begin with the Neumann groups and follow the proof of Theorem (4.1) in the case ${\Cal S}\spdot=\hat\Cal L\spdot \,\,\text{and}\,\, (p,q)=(1,1)$. There is a  commutative diagram of sheaves on $\tilde V$
$$
\CD
\hat\C M^0(\C N(D_\gamma)) @>\bar\partial^0>> \C M^1(\C N(Z-E)) @>\bar\partial^1>>
 \C M^2(\C N(Z-E)) \\
@VVV        @VVV         @AAiA\\
\C L^{1,0}_\gamma @>\bar \partial^0_N>> \C L^{1,1}_\gamma @>\bar \partial^1_N>> \C L^{1,2}_\gamma 
\endCD
$$
where $\hat\C M^0(\C N(D_\gamma))$ consists of those $\omega\in\C M(\C N(D_\gamma))$ such that $\db^0\omega\in\C M^1(\C N(Z-E))$; in particular, the middle cohomology of the top complex is $\im (H^1(\tilde U;\C N(Z-E))\to H^1(\tilde U;\C N(D_\gamma)))$.  We again get a second diagram
$$
\CD
 \hat M^0(\C N(D_\gamma)) @>\bar\partial^0>>  M^1(\C N(Z-E)) @>\bar\partial^1>>
  M^2(\C N(Z-E)) \\
@VVV        @VVV         @VV=V\\
 L^{1,0}_\gamma @>\bar \partial^0_N>>  L^{1,1}_\gamma @>\bar \partial^1_D>>  M^2(\C N(Z-E)) 
\endCD
$$
and to apply (4.15) we need to know that the range of $\db^0_N$ is closed. This is equivalent to the range of $(\db^0_N)^*$ being closed, which in turn follows from the finiteness of $H^{1,0}_N(V)$ (Lemma (4.17)). Now following the argument above following gives us a surjection 
$$
\im (H^1(V;\C N(Z-E))\to H^1(V;\C N(D_\gamma)))\to H^{1,1}_N(V)\tag4.23
$$
Moreover, if $\ker \db^1\owns\psi=\db^0_N\omega\in L^{1,1}_N(V)$, then since $\C M^0(\C N(D_\gamma))=\C L^{1,0}$ and $\C M^1(\C N(Z-E))=\C L^{1,1}$, we get $\psi=\db^0\omega$. Hence (4.23) is an isomorphism. Hence the claimed computation of $H^{1,1}_N(V)$ follows from (3.26b). 

Now notice that since $\Lambda^2\C N(D_\gamma)=\varOmega^2$ ((3.12d)), we have nonsingular pairings
$$
\C N\times\C N(D_\gamma)\to\varOmega^2
$$
and
$$
\C N(Z-E)\times\C N(N)\to\varOmega^2
$$
It now follows from Serre duality and the duality between Dirichlet and Neumann cohomology (\cite{PS, (1.3c)}) that there is an isomorphism
$$
\im (H^1(\tilde U;\C N)\to H^1(\tilde U;\C N(N)))\to H^{1,1}_D(V)
$$
(This could also have been proved as the isomorhism (4.23) was.)

Finally, the computation of $H^{1,0}_N(V)$ follows from (4.18) and (3.12d).

\enddocument